\pdfoutput=1
\RequirePackage{ifpdf}
\ifpdf 
\documentclass[pdftex]{sigma}
\else
\documentclass{sigma}
\fi

\usepackage{tabularx,mdwtab}
\usepackage{bbm}

\newcommand{\eps}{\varepsilon}
\newcommand{\One}{\mathbbm{1}}
\newcommand{\un}{\underline{N}}
\newcommand {\Cee}{{\mathbb C}}
\newcommand {\Ree}{{\mathbb R}}

\newcommand{\Kee}{{\mathbb K}}
\newcommand{\Zee}{{\mathbb Z}}
\newcommand{\fc}{\mathfrak{c}}
\newcommand{\fe}{\mathfrak{e}}
\newcommand{\fg}{\mathfrak{g}}
\newcommand{\fh}{\mathfrak{h}}
\newcommand{\fq}{\mathfrak{q}}
\newcommand{\fa}{\mathfrak{a}}
\newcommand{\fb}{\mathfrak{b}}
\newcommand{\fk}{\mathfrak{k}}
\newcommand{\fm}{\mathfrak{m}}
\newcommand{\fp}{\mathfrak{p}}
\newcommand{\fo}{\mathfrak{o}}
\newcommand{\fgl}{\mathfrak{gl}}
\newcommand{\fab}{\mathfrak{ab}}

\newcommand{\fpsl}{\mathfrak{psl}}
\newcommand{\fsl}{\mathfrak{sl}}
\newcommand{\fhei}{\mathfrak{hei}}
\newcommand{\fosp}{\mathfrak{osp}}
\newcommand {\ev}{{\bar0}}
\newcommand {\od} {{\bar1}}
\newcommand {\bcdot}{\mathbin{\hbox{\raise.4ex\hbox{\bf.}}}}
\newcommand {\tto}{\longrightarrow}

\newcommand {\fbj}{{\mathfrak{bj}}}

\newcommand {\fbr}{{\mathfrak{br}}}

\newcommand {\fel}{{\mathfrak{el}}}

\newcommand {\fme}{{\mathfrak{me}}}

\newcommand {\fbgl}{{\mathfrak{bgl}}}

\newcommand {\fout}{\operatorname{\mathfrak{out}}}
\newcommand {\fder}{\operatorname{\mathfrak{der}}}
\newcommand {\fwk}{{\mathfrak{wk}}}

\newcommand {\htr}{{\text{htr}}}

\newcommand{\bo}[2]{o^{\{#1,#2\}}}

\numberwithin{equation}{section}

\newtheorem{Theorem}{Theorem}[section]
\newtheorem{Statement}[Theorem]{Statement}

\newtheorem{Lemma}[Theorem]{Lemma}

\newtheorem{Conjecture}[Theorem]{Conjecture}
 { \theoremstyle{definition}

\newtheorem{Fact}[Theorem]{Fact}

\newtheorem{Remark}[Theorem]{Remark}
\newtheorem{Remarks}[Theorem]{Remarks}

\newtheorem*{OpenProblems}{Open Problems}
\newtheorem{Comment}[Theorem]{Comment}

 }

\begin{document}

\allowdisplaybreaks

\newcommand{\arXivNumber}{1307.1858}

\renewcommand{\PaperNumber}{032}

\FirstPageHeading

\ShortArticleName{Derivations and Central Extensions of Symmetric Modular Lie Algebras and Superalgebras}

\ArticleName{Derivations and Central Extensions of Symmetric\\ Modular Lie Algebras and Superalgebras\\
(with an Appendix by Andrey Krutov)}

\Author{Sofiane BOUARROUDJ~$^{\rm a}$, Pavel GROZMAN~$^{\rm b}$, Alexei LEBEDEV~$^{\rm c}$ and
Dimitry LEITES~$^{\rm ad}$}

\AuthorNameForHeading{S.~Bouarroudj, P.~Grozman, A.~Lebedev and D.~Leites}

\Address{$^{\rm a)}$~New York University Abu Dhabi, Division of Science and Mathematics,\\
\hphantom{$^{\rm a)}$}~P.O.~Box 129188, United Arab Emirates}
\EmailD{\href{mailto:sofiane.bouarroudj@nyu.edu}{sofiane.bouarroudj@nyu.edu}, \href{mailto:dimleites@gmail.com}{dimleites@gmail.com}}

\Address{$^{\rm b)}$~Deceased}

\Address{$^{\rm c)}$~Equa Simulation AB, Stockholm, Sweden}
\EmailD{\href{mailto:alexeyvleb@gmail.com}{alexeyvleb@gmail.com}}

\Address{$^{\rm d)}$~Department of Mathematics, University of Stockholm, SE-106 91 Stockholm, Sweden}

\ArticleDates{Received November 16, 2016, in final form February 02, 2023; Published online May 29, 2023}

\Abstract{Over algebraically closed fields of positive characteristic, for simple Lie (super)\-alge\-bras, and certain Lie (super)algebras close to simple ones, with symmetric root systems (such that for each root, there is minus it of the same multiplicity) and of ranks less than or equal to~8---most needed in an approach to the classification of simple vectorial Lie superalgebras (i.e., Lie superalgebras realized by means of vector fields on a~supermanifold),---we list the outer derivations and nontrivial central extensions. When the conjectural answer is clear for the infinite series, it is given for any rank. We also list the outer derivations and nontrivial central extensions of one series of non-symmetric (except when considered in characteristic~2), namely periplectic, Lie superalgebras---the one that preserves the nondegene\-ra\-te symmetric odd bilinear form, and of the Lie algebras obtained from them by desuperization. We also list the outer derivations and nontrivial central extensions of an analog of the rank 2 exceptional Lie algebra discovered by Shen Guangyu. Several results indigenous to positive characteristic are of particular interest being unlike known theorems for characteristic 0, some results are, moreover, counterintuitive.}


\Classification{17B50; 17B55; 17B56; 17B20; 17B40}

\setcounter{tocdepth}{3}

{\small \tableofcontents}

\section{Introduction}\label{Si}

Hereafter,
$\Kee$ is an algebraically closed field of characteristic
$p>0$ and $\fg$ is a~finite-dimensional Lie (super)algebra. This paper is a~sequel to \cite{BGL2} which contains several vital definitions, e.g., of the codifferentials, and both \cite{BGL2} and this paper are sequels to \cite{BGL1}, a~step towards the classification of simple Lie (super)algebras over $\Kee$. Laborious computations were performed with the aid of Grozman's \textsc{SuperLie} package, see \cite{Gr}.

There are
two major types of Lie (super)algebras: having selected a~maximal torus, we say that the Lie (super)algebra is ``symmetric" if
with every root $\alpha$ it has a~root $-\alpha$ of the same
multiplicity as that of $\alpha$; the Lie (super)algebras without
this property are said to be ``non-symmetric'' (or ``lopsided"). In this paper we consider $\Zee$-graded symmetric simple Lie
(super)algebras; a~review of the situation with non-symmetric
simple Lie (super)algebras, together with new results, will be given elsewhere.

Appendix contains lists of derivations and central extensions of the \textit{deforms} (results of deformations, same as transforms
are results of transformations) of the symmetric Lie (super)algebras considered in the main text.

\subsection{Background and notation we use}\label{Notat} In the literature, it is customary to call Lie (super)algebras over $\Kee$ \textit{modular}; we
find this term overused, but use it for brevity.

The description of derivations and central extensions of modular Lie \textit{super}algebras reveals several surprises, even to experts, we will dwell a~bit on these phenomena, especially abundant when $p=2$ and~3.

For the most comprehensive background on specifics of Lie (super)algebras for $p=2$, see~\cite{BGLLS}. For classification of isomorphism classes of orthogonal Lie (super)algebras without
Cartan matrix, and their periplectic analogs preserving
non-degenerate odd symmetric bilinear form for any $p$, and its numerous versions for $p=2$, see \cite{LeP} and Section~\ref{Speripl} below.

One of the main types of our examples are Lie (super)algebras with
Cartan matrix, or their ``relatives"'' such as the simple subquotients, non-trivial central extensions, and algebras of derivations. For notions and notation related to Cartan matrix, see~\cite{BGL1} to which we gladly add a~recent paper \cite{BLLoS} where the adequate definition of the notion of ``root'' over fields of positive characteristic is suggested, see also~\cite{KLS}. Observe that all indecomposable Cartan matrices of finite-dimensional modular Lie (super)algebras are symmetrizable; for their classification, see~\cite{BGL1}.

Of course, the weights of
cocycles do not depend on a~realization (the choice of Cartan
matrix), but the form of cocycles does, and hence, for each Lie
(super)algebra, we give the result for one (``simplest" in some
sense) incarnation, if there are several. We indicate the Cartan
matrix to which our cocycles correspond.

If $p=2$, we give only the
``non-super'' version of the Cartan matrix; for Lie superalgebras that superize them, we take
the ``same'' Cartan matrix (but with 0 instead of $\ev$ on the main
diagonal).

The answer for the simple subquotient $\fg(A)^{(k)}/\fc$, where $\fc$ is the center and the unspecified superscript $k$ is the smallest $k$ for which $\fg^{(k)}=\fg^{(k+1)}$, where $\fg^{(k)}$ is the $k$th derived algebra, of the Lie
(super)algebra of the form $\fg(A)$ is given with respect to the
same basis denoted by $x_i$ for positive root vectors and $y_i$ for negative root
vectors, same as for $\fg(A)$ itself. The elements $h_i:=[x_i,y_i]$ belong to the
maximal torus for every pair of Chevalley generators $x_i$ and $y_i$ (do not confuse the
\textit{generators} with the other elements of the Chevalley \textit{basis}). For the noninvertible Cartan
matrices, we denote the derivations representing the classes of outer derivations, denoted by $d_j$ in \cite{BGL1}, by
symbols $D_j$, as in comments after the proof of Statement~\ref{Statement}, in order not to confuse them with the codifferential $d$ in cohomology.

We consider algebras as rank grows, but if the general answer for a~
series considered is clear, we give it at the first instance. For
each serial Lie (super)algebra, we give the answer for the smallest
value of $p$ with which the domain of stabilization of the answer begins.

In Sections~\ref{Ssymm} and~\ref{SymmWithoutCM}, we consider Lie (super)algebras with $n\times n$ Cartan matrix $\fg(A)$ and their simple ``relatives" of the form $\fg(A)^{(k)}/\fc$; by ``rank'' of these Lie (super)algebras we mean the size $n$ of the $n\times n$ matrix $A$.

In Section~\ref{Speripl}, we consider symmetric Lie (super)algebras without Cartan matrix---$\mathfrak{o}_{\rm gen}(2n)$ and
$\mathfrak{pe}_{\rm gen}(n)$, their projectivisations---$\mathfrak{p}(\mathfrak{o})(2n):=\mathfrak{o}(2n)/\Kee 1_{2n}$ and $\mathfrak{p}(\mathfrak{pe})(n):=\mathfrak{pe}(n)/\Kee 1_{2n}$, and simple ``relatives'' of all these algebras.

Thanks to symmetry of the root system of these Lie (super) algebras, it suffices to give cocycles
of only nonpositive degree, and therefore we convene:
\begin{equation}\label{conv}
\text{\begin{minipage}[c]{14cm}
\textit{when we give cocycles representing a~basis of $H^i$, the cocycles of positive degree symmetric to those of negative degree are assumed}.
In the repeatedly used phrases ``or a
basis of $H^1(\fg;\fg)$ we can take \textit{the classes of} derivations represented by the following cocycles'' we skip ``the classes of''.
The odd cocycles are underlined.
\end{minipage}}
\end{equation}

By $[c]$ we denote the cohomology class of the cocycle $c$, but often, by abuse
of notation we write, for example, $H^2(\fg)=\operatorname{Span}(c_0, c_{-2})$ meaning $H^2(\fg)=\operatorname{Span}([c_0], [c_{-2}], [c_{2}])$, thanks to convention~\eqref{conv} that saves a~lot of paper in multi-dimensional cases.

Let the Lie (super)algebra
of ``outer derivations'' be
$\fout \fg:=\fder\ (\fg)/(\fg/\fc(\fg))$, where $\fc(\fg)$ is the center of $\fg$.

When a~Lie (super)algebra $\fg$ has incarnations for various $p$, our claim of the form ``\textit{For any $p>2$}''
means that we have checked this conjecture for several values of $p$ ($=3,5,7$ and sometimes~11 to be sure).

\subsection[Goal: classify simple finite-dimensional modular Lie algebras. What is done]{Goal: classify simple finite-dimensional modular Lie algebras.\\ What is done}\label{goal}

\subsubsection[Arbitrary $p$, but only algebras with indecomposable Cartan matrix]{Arbitrary $\boldsymbol{p}$, but only algebras with indecomposable Cartan matrix}

Weisfeiler and Kac \cite{WK}
offered a~classification of finite-dimensional Lie algebras $\fg(A)$
with indecomposable Cartan matrix $A$ for any $p>0$. For a
verification of the corrections of \cite{WK} contained in \cite{KWK} and \cite{Sk29}, precise definitions of various related notions, and superization thereof, see \cite{BGL1}.

\subsubsection[All algebras, but $p\geq 5$]{All algebras, but $\boldsymbol{p\geq 5}$}

The Kostrikin--Shafarevich method
conjecturally produced all \textit{restricted} simple finite-di\-men\-sion\-al modular Lie
algebras over algebraically closed fields of characteristic $p>7$. Block and Wilson proved the restricted version of the KSh-conjecture, see~\cite{BW}. From~\cite{BLLS}: ``This classification is implicit to this day when dealing with deforms: simple deforms of the divergence-free algebras~\cite{W} and of Hamiltonian type algebras \cite{Sk1, Sk0} were classified only several years after~\cite{BW} was published. In the divergence-free case, explicit formulas of the $p$-structure were obtained only recently, see~\cite{BKLLS}. (\lq\lq The problem of restrictedness is approached. \dots\ [But] the family of Hamiltonian algebras \dots\ is not yet handable'', see \cite[p.~357]{S}.'')

The \textit{generalized} conjecture considered \textit{all} simple algebras; it consisted of two parts: one part clearly described $\Zee$-graded Lie algebras; this part turned out to be correct for $p\geq7$.

The bigger and vague part of the generalized conjecture was more and more lucidly formulated later, thanks to ideas of Kac, who suggested to reduce the description of all deformations (infinitesimally, of $H^2(\fg;\fg)$) to the much simpler task of the description of the tensors preserved by $\fg$ which boils down to the description of filtered deformations (infinitesimally, of $H^2(\fg_-;\fg)$). The latter are covered by results due to Wilson~\cite{W} (and earlier, but with an isomorphism missed, by S.~Tyurin \cite{Tyu}),
who described classes of volume forms, and Skryabin who described classes of symplectic and contact forms, see \cite{Sk1, Sk0}. Elsewhere, we intend to review these results together with new ones, for non-symmetric Lie \textit{super}algebras. Interestingly, the well-known example of ``quantization'' of the Poisson Lie algebra was not taken into account, although it manifestly shows that the bold idea to replace computation of $H^2(\fg;\fg)$ by that of $H^2(\fg_-;\fg)$ does not always work. For one more, very interesting exceptional example of Shchepochkina, see~\cite{BGLLS}.


In \cite{KD}, it was suggested to improve the generalized KSh-procedure for non-restricted algebras by a
carefully study of deforms of certain ``standard'' examples having
added to simple ``standard'' objects one nontrivial central
extension. This improvement works for $p\geq 5$, see \cite{BGP,
KD, S}. (The proof in~\cite{KD} is absolutely correct, but its English is a~bit broken and
the Poisson Lie algebra is called Hamiltonian. However, the result was
double-checked for the case of smallest dimension, see~\cite{MeZu}.)

\subsubsection[A conjectural method for producing all simple Lie algebras for $p\geq 3$]{A conjectural method for producing all simple Lie algebras for $\boldsymbol{p\geq 3}$}

This
modification of the KSh-procedure reduce the stock of
``standard" examples to Lie algebras with indecomposable Cartan matrices.
In the improved version of the KSh-procedure, we should%
\begin{equation*}
\text{\begin{minipage}[c]{14.5cm}
\phantom{xx}(a) take the nonpositive part of
Lie algebras with Cartan matrix over $\Kee$, construct its
generalized \textit{Cartan--Tanaka--Shchepochkina prolong} (complete or \textit{partial}), and
called \textit{CTS-prolong} in what follows. (Most lucidly the CTS-prolong is defined in \cite{Shch} with examples for $p=0$ and $p>0$, see also \cite{BGLLS}.)
Proceed inductively, as described in
\cite{Ltow};\\
\phantom{xx}(b) deform the simple Lie algebras obtained at step (a), perhaps, after passing to a~first
or second derived algebra of the prolong;\\
\phantom{xx}(c) select non-isomorphic examples among the above (\textit{true} deforms).
\end{minipage}}
\end{equation*}
For $p>3$, one thus gets all simple examples.
We
\textit{conjecture that one thus gets all simple examples for $p=3$} as well, see interpretations of various Skryabin algebras, as well as Ermolayev and Frank algebras,
described in \cite{GL4}.

\subsection[Phenomena indigenous to $p=2$]{Phenomena indigenous to $\boldsymbol{p=2}$} Recall the definition of a~Lie superalgebra in characteristic~2, see~\cite{BGL2}. In~\cite{BLLSq}, we
suggested several versions of restrictedness for $p=2$ initiated in
\cite{BGL1, LeD}. Moreover, \textit{we reduced classification of simple Lie
superalgebras to that of simple Lie algebras and the $\Zee/2$-gradings of the latter}. (To classify $\Zee/2$-gradings might be rather tough, cf.~\cite{KrLe}.)

For $p=2$, the stock of algebras for the input for the CTS-prolongation is wider than for $p\geq 3$. The process is inductive, see~\cite{Ltow}.
In~\cite{SkT1}, Skryabin showed that certain semisimple Lie
algebras also have to be added to the list of ``standard'' examples, see also~\cite{GZ}. In addition, there are Eick algebras~\cite{Ei}, not explained yet.

Recently, Skryabin's approach to classification of symplectic forms was extended in \cite{KKCh,Kon} to ``non-alternate" forms, see also~\cite{LeP}. (Observe that Skryabin applied the adjective ``Hamiltonian'' not only to algebras preserving these forms, which is the standard usage of the term, but also to the forms themselves, which is at variance with the usage of the term in theoretical mechanics and differential geometry, cf.~\cite{Sk0} with earlier~\cite{Sk1}.)

Additionally, there might appear new CTS-prolongs of certain pairs $(\fg_{-}, \fg_0)$ with exceptional $ \fg_0$-modules $\fg_{-1}$.

There are more types of symmetric simple Lie
(super)algebras: the ``non-alternate'' orthogonal ones, see
\cite{LeP}, and queerification of symmetric Lie algebras, see
\cite{BLLSq}. For details and examples, see \cite{BGLLS1, BGLLS}.

\subsection[Phenomena indigenous to $p=3$. Pre-Lie superalgebras]{Phenomena indigenous to $\boldsymbol{p=3}$. Pre-Lie superalgebras}\label{ssDerpIs3}

Recall the definition of a~Lie superalgebra in characteristic 3 and a~pre-Lie superalgebra (indigenous to $p=3$), see \cite{BGL2}.

Any~Lie (super)algebra with a~non-degenerate invariant symmetric bilinear form $B$ will be briefly called a~\textit{NIS-(super)algebra}. The \textit{double extension} $\mathfrak{g}$ of a~NIS-(super)algebra $\mathfrak{a}$ is the result of simultaneous adding to $\mathfrak{a}$ a~central element $c$ and a~derivation $D$ so that the direct sum \textit{as spaces} $\fg:={\mathcal K} \oplus \mathfrak{a} \oplus {\mathcal K} ^*$, where ${\mathcal K} :=\Kee c$ and ${\mathcal K}^* :=\Kee D$, is a~NIS-algebra. Most known examples of double extensions are affine Kac--Moody algebras (over $\Cee$ or $\Ree$).

Computing double extensions of restricted Lie superalgebras the authors of \cite{BeBou2} discovered an interesting phenomenon that helped us to interpret certain derivations and central extensions.

As we know from \cite{BeBou1} succinctly summarized in \cite{BLS}, the double extension is \textit{indecomposable}, i.e., not isomorphic to the direct sum of two ideals $\mathfrak{a}$ and ${\mathcal K} \oplus {\mathcal K} ^*$, if the following conditions are satisfied:
\begin{enumerate}\itemsep=0pt
\item[a)] the derivation $D$ of $\mathfrak{a}$ \textit{must be outer} for any $p$; if $p=2$, and $D\in(\mathfrak{out}\, \mathfrak{a})_\od$, then, moreover, the condition $D^2=0$ is a~must, see \cite{BKLS};
\item[b)] the central extension has to be non-trivial, see \cite[Section~8]{BLS}.
\end{enumerate}

If $p=3$ and a~derivation $D$ preserves an \textit{invariant and symmetric} (IS for short) bilinear form~$B$, then $\omega(a,b) := B(Da, b)$ does not have to describe a~central extension of a~Lie superalgebra~$\fg$, because it may happen that
\[
\omega\big(x, x^2\big)\neq 0 \qquad \text{for an element $x\in\fg_\od$.}
\]
So if we try to construct a~central extension using this $\omega$, the resulting superalgebra wouldn't necessarily satisfy the indigenous for $p=3$ part of the Jacobi identity for 3 odd elements in the form $\big[x, x^2\big]=0$; though it would still satisfy
\begin{gather}\label{JI3}
[u, [v,w]]+[v,[w,u]]+[w,[u,w]]=0 \qquad \text{for any $u,v,w\in\fg_\od$.}
\end{gather}

 Whether this $\omega$ is a~cocycle or not, depends on how we define cochains and the action of the exterior derivation on them; accordingly, we calculate the cocycles correctly or not (either by pen on paper or by means of the code \textsc{SuperLie}).

 If we use the \textit{used to be} standard definition according to which the cochains are anti-supersymmetric functions $\omega$, then the condition $d\omega=0$ assures us that $\omega$ is a~cocycle but not of a~Lie superalgebra. By definition, the \textit{pre-Lie superalgebra} satisfies the Jacobi identity for different triples of arguments, e.g., \eqref{JI3}, but not $\big[x, x^2\big]=0$ for $x$ odd.

If $p\neq 2,3$, then a~function $\omega$ such that $\omega\big(x,x^2\big)\neq 0$ for some odd $x$ cannot be a~cocycle,
since $d\omega(x,x,x) = 3\omega(x,[x,x])$,
therefore \textsc{SuperLie}, and we, consider $d\omega\neq 0$ if $\omega\big(x,x^2\big)\neq 0$ even for $p=3$.

As we have discussed in \cite{BGL2}, the \textit{used to be} standard definition of cochains is not good for computing infinitesimal deformations of Lie superalgebras in characteristic $2$. It turns out that it is not good for finding central extensions in characteristic~$3$ (and, perhaps,~$2$) either.

The definition of cochains can be interpreted in different ways, depending on what we mean by cochains and by a~derivation preserving a~symmetric bilinear form.

In the standard way, one defines cochains, as super anti-symmetric functions, then $B(Dx, y)$ \textit{is} a~cocycle, it is just that not every 2-cocycle with coefficients in the trivial module corresponds to a~central extension; \textsc{SuperLie} selects only those cocycle which do correspond to central extensions.

Corrections to the standard definition: we define cochains as it is implemented in \textsc{SuperLie}:
\begin{itemize}\itemsep=0pt
\item $B(Dx, y)$ is not necessarily a~cocycle.
\item We change the definition of ``a derivation preserving an IS form''. Namely, in addition to just requiring
\[
 B(Dx,y) = (-1)^{p(D)p(x)}B(x,Dy) \qquad \text{for any $x,y\in\fg$}
 \]
 \textit{we should add the following condition if $p=2$ or $3$}:
\end{itemize}

If $p=2$ then
\[
B(Dx,x)=0 \qquad \text{for all $x\in\fg_\od$}.
\]
Since, for any Lie superalgebra $\fg$, we want $\fder \fg$ be a~Lie superalgebra, one has to add the condition
\[
D\big(x^{2}\big) = [Dx,x] \qquad \text{for any $x\in\fg_\od$ and $D\in\fder \fg$.}
\]

If $p=3$, then \textit{we should add the condition}
\begin{gather}\label{der3}
B\big(Dx,x^2\big) = 0\qquad \text{for all $x\in\fg_\od$}.
\end{gather}

If $p\neq 2,3$ or if $p=3$ and $D$ is even, then the condition \eqref{der3} is satisfied automatically, since
\begin{align*}
B\big(Dx, x^2\big) &= \tfrac12B(Dx, [x,x]) = \tfrac12B([Dx,x], x) = \tfrac14B(D[x,x], x)\\
&=
\tfrac12B\big(Dx^2, x\big) = \tfrac12B\big(x, Dx^2\big) = (-1)^{p(D)} \tfrac12B\big(Dx, x^2\big)= \tfrac12B\big(Dx, x^2\big),
\end{align*}
so $B\big(Dx, x^2\big) = 0$.

\subsection{Derivations and central extensions of simple Lie (super)algebras}\label{ssOutDer} The classification
of \textit{derivations and
central extensions of simple Lie (super)algebras} is of interest \textit{per se},
but the knowledge of the result is also needed in one of the powerful
methods of classification of simple Lie (super)algebras: see~\cite{LSh}.

\textbf{On importance of restrictedness.} At first, Shafarevich and Kostrikin
considered in their conjecture only \textit{restricted} Lie
algebras. Recently, on a~different occasion, Deligne gave us an
advice to look, if $p>0$, at the groups (geometry) rather than Lie
algebras, see Deligne's appendix in \cite{LL}. Deligne advised us
to restrict the classification problem of simple Lie
(super)algebras, and modules over them, to restricted ones---at least, to begin with---because
only restricted Lie algebras correspond to geometry.
Conceding that the study of restricted
Lie (super)algebras is a~natural first step, we consider in this paper naturally arising nonrestricted Lie (super)algebras as well; they are often needed to describe restricted ones, see also~\cite{BLLSq}. The referees remind us also that non-restricted Lie algebras are applied in the study of
$p$-groups, and naturally appear in representations of quantum groups (algebras) $U_q(\fg)$ for the parameter $q$ equal to a~$p$th roots of unity.

The lack of outer derivations is a~sufficient condition for the
algebra to be restricted. The complete inventory of central
extensions is of interest \textit{per se}, but becomes indispensable in
classification of simple Lie (super)algebras if $p=0$, see \cite{LSh}. To
classify (at least certain types of) simple \textit{vectorial} Lie (super)algebras (i.e., Lie superalgebras realized by means of vector fields on a~supermanifold) \`a~la~\cite{LSh} if $p>0$ is
another possible application of the results of this paper.

In this subsection, $\operatorname{char} \Kee=0$ and $\fg$ is a~Lie algebra, unless otherwise stated.
For the Lie
algebra with a~non-degenerate invariant symmetric
bilinear form (NIS for short)~$b$, the problems of describing outer derivations and
nontrivial central extensions are practically equivalent. Indeed, the
outer derivations are described by cocycles representing the classes
of $H^1(\fg; \fg)$ whereas the nontrivial central extensions are
defined by cocycles representing the classes of $H^2(\fg)$, i.e.,
with trivial coefficients. In presence of a~NIS $b$, we have
$\fg\simeq\fg^*$, and hence $H^1(\fg; \fg)\simeq H^1(\fg; \fg^*)$.

For any cocycle $\omega$ representing the class $[\omega]\in H^2(\fg)$, there is an embedding
\begin{gather}
j\colon \ H^2(\fg)\tto H^1(\fg; \fg^*),\nonumber\\
j(\omega)\colon \ X\longmapsto
j(\omega)(X), \qquad \text{where} \ (j(\omega)(X))(Y)=\omega(X, Y) \ \text{for any $X,Y\in \fg$}.\label{mapJ}
\end{gather}
The
following sequence is exact:
\begin{gather}\label{exact}
0 \tto H^2(\fg) \stackrel{j}{\tto} H^1(\fg; \fg^*)
\stackrel{v}{\tto} (S^2(\fg))^\fg \stackrel{K}{\tto} H^3(\fg)\tto
H^2(\fg; \fg^*) \tto H^1\big(\fg;S^2(\fg)\big),
\end{gather}
where $\big(S^2(\fg)\big)^\fg$ denotes the space of $\fg$-invariant
symmetric bilinear forms on $\fg$, and the maps $v$ and $K$ are
defined as follows:
\begin{gather*}
v(f)(X, Y) = f(X)(Y) + f(Y)(X)\qquad \text{for any $X,Y\in \fg$ and $f\in Z^1(\fg; \fg^*)$},\\
K(b)(X, Y, Z) = b([X, Y], Z) \qquad \text{for any $X,Y, Z\in \fg$ and $b\in \big(S^2(\fg)\big)^\fg$}.
\end{gather*}
The beginning of the sequence~\eqref{exact} was discovered by
Koszul, who introduced the map $K$; then the sequence was further
extended to the right, see~\cite{NW}. The fact ``the map $j$
is an embedding'' is proved in \cite{Dz5}; the exactness of the sequence
\eqref{exact} is proved in \cite[Proposition~7.2]{NW}.

In 1992, in a~letter to D.L., Dzhumadildaev wrote about \textit{symmetric $($co$)$homology} $HS^{\bcdot}$, later elaborated in \cite{DzS,DzZ}. The term $HS^0$ is $\big(S^2(\fg)\big)^\fg$ and the exactness of the sequence \eqref{exact} means, for example, that for complex simple finite-dimensional Lie algebras, the 3-cocycle $K(b)(-, -, -)$ is nontrivial, i.e., $\operatorname{Ker} K = 0$, hence, $H^2(\fg) \simeq H^1(\fg; \fg^{\ast}) = 0$.

If $\fg$ is a~vectorial Lie algebra, it can happen that $\operatorname{Ker} K$ is nontrivial, e.g., for a~modular~$\mathfrak{svect}(3)$, see \cite{BKLS}, and for Hamiltonian Lie algebras $\mathfrak{h}$ on tori where there are invariant forms on $\mathfrak{h}$. Then, ``Leibniz central extensions'' arise, i.e., a~central extension of the Lie algebra, the result of which is not a~Lie algebra, but a~\textit{Leibniz algebra} which satisfies only the Jacobi identity, but not the anti-symmetry one.

Thus, a~part of $H^1(\fg; \fg^{\ast})$ describes Lie central extensions, while the other part describes Leibniz central extensions.

\textbf{Superization of the sequence \eqref{exact}}. For the case
where all $\fg$-invariant symmetric bilinear forms $(S^2(\fg))^\fg$
are \emph{even}, the arguments of \cite{NW} are literally applicable
to Lie \textit{super}algebras and help us to verify the results if $p\neq 2$. (Note that the title of \cite{BKLS} is misleading: there are symmetric and anti-symmetric bilinear forms on superspaces, but there are no "super symmetric" forms. For an elucidation of this, see \cite[Section~2.4.2]{KLLS}.)

If $p=2$, there are subtleties we'll discuss in what follows.

It is an \textit{open problem} to construct an analog of the sequence~\eqref{exact} if NIS on $\fg$ is odd.

\subsubsection{What can be said a~priori: ``extra derivations''}\label{apriori1} Let $\fg=\fg(A)$ be a~Lie (super)algebra with a~noninvertible
Cartan matrix. Then, $\fg$ has a~central element $z$ and a~grading
element $D$, both lying in the maximal torus $\fh$ and
$D\not\in\fg^{(1)}$, where $\fg^{(1)}$ is the derived algebra.

Let $\{X_1, \dots, X_n\}$ be a~basis of $\fg$; let the $X_i^*\in \fg^*$ be the dual basis.
Let $\widehat X_i=\Pi(X_i^*)\in \Pi(\fg^*)$ \index{$\widehat X_i=\Pi(X_i^*)\in \Pi(\fg^*)$}
when cochains of degree $>1$ are considered; otherwise, the change of parity $\Pi$ can be ignored.

\begin{Statement}[unexpected]\label{Statement} Let
$\operatorname{rk}(A)=\operatorname{size}(A)-1$, so $\fg = \fg^{(1)}\ltimes \Kee\cdot D$ is a~semi-direct sum. The operator $M$ on $\fg$ corresponding
to the cochain $z\otimes \widehat{D}$, i.e., such that
\begin{equation}\label{AL} MD=z \qquad \text{and} \qquad M\big(\fg^{(1)}\big)=0, \end{equation}
is an outer derivation of $\fg$.
\end{Statement}

\begin{proof} To show that $M$ is a~derivation, we need to show that
\[
M[x,y]=[Mx,y]+[x,My] \qquad \text{for any $x,y\in\fg$}.
\]
Indeed, the left-hand side vanishes thanks to \eqref{AL}, and the
right-hand side vanishes due to the fact that $\operatorname{Im} M\subset \fc$, where
$\fc$ is the center of $\fg$.

There is no $x\in\fg$ such that $M=\operatorname{ad}_x$. Indeed, otherwise
$[x,D]=z$, but $D\in\fh$, so $[D,\fh]=0$ and
$[D,\fg_\pm]\subset\fg_\pm$. Thus, $z\not\in [D,\fg]$ since
$z\in\fh$.
\end{proof}

\begin{Remark}\label{RemD} If $\operatorname{size}(A)-\operatorname{rk}(A)=k$, then, as is known from \cite{BGL1}, there are~$k$ linearly independent
central elements $z_1,\dots,z_k$, and $k$ grading elements
$D_1,\dots,D_k$. Similarly to the above, there are~$k^2$ linearly
independent outer derivations of~$\fg$, i.e., operators $M_{ij}$
corresponding to the cochains $z_i\otimes \widehat{D}_j$.
\end{Remark}

\subsubsection{On central extensions and derivations of simple Lie algebras for $p>0$}\label{known} The results of Ibraev \cite{Ib, Ib2} and Permyakov \cite{Per} imply a~complete description of central extensions and derivations of simple Lie algebras of ``classical type'' (i.e., of the form $\fg(A)$ or $\fg(A)^{(i)}/\fc$, except for $\fwk$ and $\fbr$ cases) for $p>0$. Therefore, \textit{we will not consider these cases}, except for an occasional illustration (but we do consider $\fwk$ and $\fbr$ cases---new results).

\subsubsection{What can not be said a~priori and we have to compute}\label{apriori2}
It is clear that if $\big(S^2(\fg)\big)^\fg\neq 0$, there is an isomorphism
$H^1(\fg; \fg^*)\simeq H^2(\fg)$, if $p\neq 2$. This is wonderful, but we
are currently interested in $H^1(\fg; \fg)$ and not in $H^1(\fg; \fg^*)$, and therefore
\[
\begin{array}{l}
\text{if $\big(S^2(\fg)\big)^\fg=0$ or if NIS on $\fg$ is odd,}\\
\text{we have to compute both $H^1(\fg; \fg)$ and $H^2(\fg)$}.
\end{array}
\]

\subsubsection{What can be said a~posteriori}\label{aposteriori} If there is a~NIS on $\fg$, then $\dim \big(S^2(\fg)\big)^\fg=1$, unless $p=2$; for a~proof, see \cite{KLLS}. Let us consider examples.

For $p\neq 2$, the fact that $j$ is an embedding implies that if there is an even
NIS $b$ on $\fg$ (e.g., $\fg$ has an invertible Cartan matrix)
and $H^1(\fg;
\fg)=0$, then $H^2(\fg)=0$.

For $p=2$, several results of our computations contradicted the fact that $j$ is embedding. But, indeed,
\[
\text{if $p=2$, then the map $j$, see \eqref{mapJ}, is not necessarily an embedding!}
\]

By the standard definition, $n$-cochains of the Lie superalgebra $\fg$ with values in the module $V$ are super anti-symmetric $n$-linear functions on $\fg$ with values in $V$. In particular, 1-cochains are just linear functions on $\fg$ with values in $V$. The standard exterior differential (codifferential in cohomological terms) is defined on the space of 1-cochains $c_1$ as follows (the signs are irrelevant if $p=2$):{\samepage
\[
 dc_1(x,y) = c_1([x,y]) - xc_1(y) + yc_1(x)\qquad \text{for any $x,y\in\fg$}.
\]
Thus defined, the codifferential $d$ does not take squaring into account if $p=2$.}

\textit{The definition for $p=2$}: The definition of $n$-cochains of Lie superalgebras used in \textsc{SuperLie} is not equivalent to the above definition of $n$-cochains as super anti-symmetric $n$-linear functions when $n\geq p>0$. In particular, in the case of $p=2$, the \textsc{SuperLie} definition of $C^2(\fg; M)$ is equivalent to the space of all pairs $(c,\fq)$, where $c$ is an anti-symmetric\footnote{Note that it is antisymmetric, not super anti-symmetric; i.e., $c(x,x)=0$ for all $x\in\fg$, independently of the parity of $x$.} bilinear function on $\fg$ with values in $M$ and $\fq$ is a quadratic form on $\fg_\od$ with values in $M$ such that its polar form is the restriction of $c$ to $\fg_\od$.

Recall the definition of $(c,\fq)\in C^2(\fg)$ and ${\mathfrak d}^1$, see \cite[Section ``Cohomology of Lie superalgebras in degrees 1 and 2 for $p=2$'']{BGL2}. \textit{In what follows, we use the notation $a\wedge b$ to represent the pair $a,b\in\fg^*$ such that }
\[
c(x,y) := a(x)b(y) + a(y)b(x) \qquad \text{for all $x,y\in\fg$ and $\fq(x):=a(x)b(x)$ for all $x\in\fg_\od$.}
\]
The codifferential ${\mathfrak d}^1$ maps a 1-cochain $c_1\in C^1(\fg; M)$ (the definition of $C^1(\fg; M)$ is equivalent to the classical one, i.e., linear functions from $\fg$ to $M$) to the pair $(dc_1, \fq)$ where $dc_1$ is as in the classical definition and
\[
\text{$\fq(x)=c_1(x^2) + xc_1(x)$ for any $x\in\fg_\od$.}
\]

In particular, if $x\in\fg_\od$ and $c_1\in C^1(\fg)$ is such that $c_1(x)=0$ and $c_1(x^2)=1$, then ${\mathfrak d}^1c_1$ is non-zero, since $\fq(x)=1$.

For several Lie superalgebras $\fg$ in characteristic 2 we got
\begin{gather}\label{inP=2}
\dim H^2(\fg)\geq\dim H^1(\fg; \fg^*).
\end{gather}
This can indeed be so if $p=2$, but this can not be true for the classical cochains (i.e., super anti-symmetric functions on which the codifferential does not take squaring into account) for $p\neq 2$. Observe that for $p=3$, the definition of 3-cochains is also different from the classical, so a~new phenomena concerning $H^2$ may arise.

\begin{Remarks} 1) $\fout(\fh)\simeq H^1(\fh;\fh)$ is not a~must only for $p=2$. Since the polar form is always anti-symmetric (i.e., $c(x,x) = 0$) for $p=2$, we could equivalently define $c\in C^2(\fg; M)$ to be anti-symmetric whatever the parity of $x$.

In characteristic $2$, in all examples computed in this paper, the space $\fout(\fg)$ of outer derivations is isomorphic to $H^1(\fg;\fg)$. Although the Chevalley--Eilenberg codifferential does not consider the squaring, \textsc{SuperLie} does take squaring into account when computing Lie superalgebra cohomology.

The isomorphism $\fout(\fg)\simeq H^1(\fg;\fg)$ does not take place in general. Consider a~$1|1$-dimensional Lie superalgebra~$\fg$ spanned by an odd basis element $x$ and even $x^2$ (with ${\big[x,x^2\big]=0}$ by the Jacobi identity). Since the bracket in $\fg$ identically vanishes, it follows that~$d$
 is the zero operator on~$C^{\bcdot}(\fg;\fg)$; this means that $Z^1(\fg; \fg)=C^1(\fg; \fg)\simeq\fgl(\fg)$ and $B^1(\fg; \fg)=0$, and hence $H^1(\fg; \fg)\simeq\fgl(\fg)$, so $\operatorname{sdim} H^1(\fg; \fg)=2|2$. However, any derivation $f$ of $\fg$ must satisfy
 \begin{gather}\label{(2)}
 f\big(x^2\big) = [x,f(x)] \qquad \text{for any $x\in \fg_\od$}
\end{gather}
since $[x,y]=0$ for any $y\in \fg$, so $\fout(\fg) = \fder \fg$ is $1|1$-dimensional. (The idea is that to be a~derivation, a~linear map $f\colon \fg\to\fg$ has to satisfy the only condition~\eqref{(2)}. So we can freely choose $f(x)$, which gives $1|1$-dimensional space of degrees of freedom, and $f(x^2)$ is determined by the condition \eqref{(2)}, so there is no additional freedom.)

2) On the other hand, we have the following lemma.

\begin{Lemma} If $p=2$, then $\fout(\fh)\simeq H^1(\fh;\fh)$ for any centerless Lie superalgebra~$\fh$. \end{Lemma}

\begin{proof}
A linear map $f\colon \fh\tto \fh$ belongs to $Z^1(\fh;\fh)$ if and only if it satisfies the condition (we skip signs since $p=2$)
\begin{gather}\label{(1)}
 f([x,y]) = [f(x),y] + [x,f(y)] \qquad \text{for any $x,y\in \fh$}.
\end{gather}
Observe that $f$ is a~derivation of $\fh$ if and only if it satisfies conditions~\eqref{(1)} and~\eqref{(2)}.

 But if $\fh$ has no center, then \eqref{(2)} follows from \eqref{(1)}: if $f$ satisfies \eqref{(1)}, then for $x$ odd and any $y$, we have
\[
 f\big(\big[x^2,y\big]\big) = \big[f\big(x^2\big),y\big] + \big[x^2,f(y)\big].
\]
 On the other hand,
\begin{align*} f\big(\big[x^2,y\big]\big)& = f([x,[x,y]]) = [f(x),[x,y]] + [x,f([x,y])] \\
& =[f(x),[x,y]] + [x,[f(x),y]] + [x,[x,f(y)]] \\
&= [[x,f(x)],y] + \big[x^2,f(y)\big].
\end{align*}
 So we see that
 \[
 \big[f\big(x^2\big),y\big] = [[x,f(x)],y] \qquad \text{for any $y$}.
 \]
 Since $\fh$ has no center, this means that $f\big(x^2\big) = [x,f(x)]$.
 So, $Z^1(\fh;\fh) \simeq \fder \fh$. And since $B^1(\fh;\fh)$ is the space of inner derivations of $\fh$ in any case, we have $H^1(\fh;\fh) \simeq \fout(\fh)$. \end{proof}
\end{Remarks}

\begin{Conjecture}[on inequalities between $\dim H^1(\fg; \fg^*)$ and $\dim H^2(\fg)$] \label{ineq} For $p=2$, the inequality \eqref{inP=2} should be modified to
\begin{gather*}
\dim H^1(\fg; \fg^*)\leq \dim H^2(\fg) \leq \dim H^1(\fg; \fg^*) + \dim \fg_\od - \operatorname{codim} \fg^{[1]},
\end{gather*}
where $\fg^{[1]} := \operatorname{Span}([x,y] \mid x,y\in \fg)$.
Do not confuse $\fg^{[1]}$ with the first derived algebra $\fg^{(1)}$.
\end{Conjecture}

Since $\operatorname{codim} \fg^{[1]}$ is the dimension of the space of $2$-co\-boun\-daries of the form $a\wedge a$, we have
\[
\dim \fg_\od - \operatorname{codim} \fg^{[1]}=\dim\operatorname{Ker} j,
\]
see formula~\eqref{mapJ}.

\subsection{Super goal: classify simple finite-dimensional modular Lie superalgebras}

First, let us perform an inventory.
For the classification of simple Lie superalgebras of the form $\fg(A)$ for $A$ indecomposable and invertible, and simple subquotients of $\fg(A)$ for $A$ indecomposable but non-invertible, see~\cite{BGL1}. All other types of simple Lie
superalgebras are not classified yet.

For $p>5$, we have only
conjectural methods for producing all of them (see the list of examples in Sections~\ref{1.6.1} and~\ref{1.6.2}), but we are sure these methods yield the complete classification.

For $p=3$ and $5$, we only hope our methods produce all simple Lie superalgebras.

For $p=2$, expectedly the most difficult case, we offered two methods that pro\-duce all simple Lie superalgebras out of simple Lie algebras, and proved that every simple Lie superalgebra is obtained from a~simple Lie algebra by one of these two methods, see~\cite{BLLSq}. This result is astonishing: in the most difficult case we have a~complete classification! But there is a~catch: the classification is modulo the classification of the simple Lie algebras which seems to be far out of reach at the moment. Moreover, one of these two methods requires classification of $\Zee/2$-gradings of simple Lie algebras, which is a~very tough, if not wild, problem, see~\cite{KrLe} where it is solved in a~particular case.

\subsubsection{Symmetric simple Lie superalgebras}\label{1.6.1} Here is the list. Deforms of non-symmetric algebras can be symmetric. We are not aware of the opposite examples.
\begin{enumerate}\itemsep=0pt
\item[A)] Lie superalgebras with indecomposable Cartan matrices, their
simple subquotients. These are classified in \cite{BGL1}.
\item[B)] Queer Lie superalgebras (any $p$) and queerified symmetric simple Lie algebras ($p=2$).
\end{enumerate}

The deforms of the symmetric Lie (super)algebras (except for numerous queerifications
and superizations of Chebochko's examples, see~\cite{Ch1}) are completely
described in~\cite{BGL2}, where, apart from simple Lie
(super)algebras with indecomposable Cartan matrices and their simple relatives, there are considered other types of simple Lie
(super)algebras, but only of rank $\leq 8$.

\subsubsection{Non-symmetric simple Lie superalgebras}\label{1.6.2} These simple Lie superalgebras are usually realized as vectorial, i.e., by means of vector fields (unlike matrix or linear superalgebras realized by matrices or linear operators). In this paper we consider symmetric simple Lie superalgebras, the non-symmetric simple Lie superalgebras will be considered separately \textit{with one exception}: periplectic Lie superalgebras become symmetric for $p=2$, so we consider periplectic cases here, for any $p$ for consistency.

For $p> 5$, one can consider the direct modular
versions of the simple vectorial Lie superalgebras over $\Cee$, see~\cite{LSh}, and their filtered deforms analogous to the Lie algebra
case, see, e.g.,~\cite{Kos}. For the input of CTS-procedure (most lucidly described in~\cite{Shch}) take
\begin{itemize}\itemsep=0pt
\item nonpositive parts of Lie superalgebras of the form $\fg(A)$, or
\item nonpositive parts of $\fg^{(i)}(A)/\fc$, see
\cite{BGLL, BGL1}; or
\item pairs $(\fg_-,\fg_0)$ producing simple
exceptional vectorial Lie superalgebras, cf. \cite{LSh}.
\end{itemize}
Then,
conjecturally (but undoubtedly), we get as the result of the CTS-procedure all simple non-symmetric $\Zee$-graded Lie superalgebras. Together with their deforms, we get all simple superalgebras.

For $p= 5$,
new examples (not discovered yet) could be added to these inputs. For $p< 5$, examples already found should be added. In particular,

For $p=3$, see \cite{BGL3, BL}; two super analogs of Melikyan
examples appear; unlike the Melikyan algebras, both these super
analogs seem to qualify as ``standard" examples, hopefully
exhausting them, together with the modular versions of simple
vectorial Lie superalgebras over $\Cee$, see~\cite{LSh}.

For $p=2$, the candidates for the role of ``standard" examples are numerous and diverse, see \cite{BGLLS1,BGLLS,BLLS, BLLSq,GL4,
ILL,LeP}.

\subsection{Symmetric Lie algebras (known results)}\label{ss1.5} The derivations and central
extensions are computed at the moment only for certain of
$\Zee$-graded simple Lie (super)algebras.

For $p>2$, the spaces of outer
derivations of Lie algebras with indecomposable Cartan matrix and
their simple relatives of the form $\fg^{(i)}(A)/\fc$ are computed
in \cite{Ib}; \textit{they are at most $1$-dimensional, except for $\fpsl(3)$
for $p=3$}.

For $p=2$, the dimensions of the spaces of outer derivations
of ``symmetric'' simple Lie algebras with Cartan matrix and simple
relatives of the form $\fg^{(i)}(A)/\fc$ of Lie algebras with
indecomposable Cartan matrix (except for $\fwk$ and $\fbr$ cases) are described in~\cite{Per}.
In~\cite{Per}, there are also described dimensions
of the spaces of outer derivations of the Lie algebras
$\fg_\Kee:=\fg_\Zee\otimes_\Zee\Kee$ obtained from the integer form
$\fg_\Zee$ spanned by the Chevalley basis of the simple complex Lie
algebra $\fg$ and their quotients modulo center, $\fg_\Kee/\fc$.

Observe that Lie algebras with roots of series $B$ and $C$, and the exception $F_4$ yield
non-simple Lie algebras\footnote{It is unclear to us where such non-simple Lie algebras and
their outer derivations might be of interest; however, see Section~\ref{sssUnclear}.} of the form $\fg_\Kee/\fc$, while the exceptional
Lie algebra $\fg(2)$ (often denoted by its root system: $G_2$) turns under the passage to the simple subquotient into $\fpsl(4)$ if $p=2$.

\subsection{Double extension: an interesting and important notion}\label{DE} Recently, there was distinguished a~notion in which a~central extension, and a~derivation, \textit{and a~NIS} on a~given Lie (super)algebra are considered \textit{simultaneously}. Among the most interesting and most known examples of double extensions we encounter affine Kac--Moody Lie (super)algebras with Cartan matrix over $\Cee$ and $\fgl(np)$ in characteristic $p>0$; for a~succinct overview of double extensions together with several new results explaining certain known facts, see \cite{BLS}.

Let $\fh$ be a~Lie superalgebra over any field with a~NIS $B_\fh$ on it, and $D\in\fder\ \fh$ a~derivation such that $B_\fh$ is $D$-invariant, i.e.,
\[
B_\fh(D(a),b)+B_\fh(a,D(b))=0 \qquad \text{for any $a,b\in \fh$.}
\]
Then, if $\fg$ has a~non-trivial central extension, there exists a~Lie superalgebra structure on the space $\fg:={\mathcal K} \oplus \fh\oplus {\mathcal K}^*$, called the \textit{double extension} (or just $D$-\textit{extension}) of $\fh$, where ${\mathcal K}=\Kee c$ and ${\mathcal K}^*=\Kee c^*$, where $c^*=D$, is the dual space, defined as follows (for any $a,b\in\fh$ and $D:=c^*$):
\begin{equation*}
[c,\fg]=0, \qquad [a,b]_\fg=[a,b]_\fh+B_\fh (D(a),b)c, \qquad [D,a]=D(a).
\end{equation*}

If $p\neq 2$, the notion of D-extension has a~natural super version provided $D^2=0$ if $D$ is odd, see~\cite{ABB, ABBQ}.

If $p=2$, the definition of double extension is non-trivial, see \cite{BeBou1, BeBou2}; it helps us to interpret one example in~Lemma~\ref{gl4_22} below.

\subsection{Our results} Over $\Kee$ for $p>0$, we listed the outer derivations and
nontrivial central extensions of simple finite-dimensional Lie
(super)algebras of rank $\leq 8$ for algebras with symmetric root
systems.

If $p=2$, we considered several types of algebras admitting both symmetric and non-symmetric root system---periplectic Lie superalgebras, and their desuperizations.

When the pattern was clear we gave the answer for any rank.
Computations were performed with the aid of \textsc{SuperLie} package,
see~\cite{Gr}.

These computations confirmed the results of Ibraev~\cite{Ib} and Permyakov~\cite{Per}
concerning the simple Lie algebras; moreover, we exhibit the
explicit cocycles; they are sometimes needed. In certain cases we prove our claim concerning the general rank thanks to Permyakov's result~\cite{Per}.
He gives arguments to prove that \lq\lq there are that many cocycles"; here we exhibit exactly that many cochains and it is easy to check that they are (a) cocycles, (b) correspond to the outer derivations, (c) are linearly independent.

Observe that for $p=2$,
there are other simple Lie algebras of ``symmetric'' type in addition
to those considered in \cite{Per}; we considered these algebras as
well. For all cases where the description of central extensions did
not follow from the general theory of formula~\eqref{exact}, we gave the description
separately.
\begin{itemize}\itemsep=0pt
\item Certain results concerning outer derivations of Lie algebras with
noninvertible Cartan matrix are counterintuitive, e.g., see Remark~\ref{2.2.2aa}; they are based on the innocent-looking
Statement~\ref{Statement}.

\item Answering the remarks of the referees we added several new results related to the latest studies
of \textit{double extensions}.

\item Andrey Krutov wrote an appendix classifying derivations and central extensions of the deforms with even parameter (listed in~\cite{BGL2}) of the Lie (super)algebras considered in the main text.
\end{itemize}

\textit{In what follows ``rank'' refers to the size of Cartan matrix, or $m+n$ for $\mathfrak{oo}_{B}(m|n)$ which has no Cartan matrix}.

\begin{OpenProblems}\label{SSop}\quad
\begin{enumerate}\itemsep=0pt
\item[1a)] What is the analog of the exact sequence
\eqref{exact} for Lie superalgebras with an even symmetric $\fg$-invariant bilinear
form on it for any $p$?
\item[1b)] What is the analog of the exact sequence
\eqref{exact} for Lie superalgebras with an \emph{odd} symmetric $\fg$-invariant bilinear
form on it for any $p$?
\item[1c)] Does the sequence
\eqref{exact} remain exact if $p>2$? If not, what happens?
\item[2)] For $p=2$, what are the rules according to which the spaces $H^1(\fg;\fg)$, where $\fg=\fg(A)$,
vary as $\ev$'s turn into $0$'s on the diagonal of ``the same" Cartan matrix, see \cite{BGL1}?
\item[3)] Describe (super)groups of automorphisms of simple
Lie (super)algebras not considered in~\cite{FG}.
\end{enumerate}
\end{OpenProblems}

\section{Symmetric Lie (super)algebras}\label{Ssymm}

In this section, we consider ``symmetric'' Lie (super)algebras, except for---if $p=2$---the orthogonal series and queerifications of symmetric Lie algebras. For further notions and notation related to Cartan matrix, see \cite{BGL1}. We let the positive Chevalley \textit{generators} be of degree~1, and the
elements of Chevalley \textit{basis} they generate, by $x_i$, the corresponding negative
basis elements by $y_i$; we set $h_i:=[x_i, y_i]$ for the \textit{generators} $x_i$ and $y_i$ of degree $\pm 1$ only---this grading is referred to as \textit{standard} or \textit{principal}.
The cocycles are indexed by their degree induced by the $\Zee$-grading of $\fg$; the superscript labels independent cocycles of the same degree or weight.

If $\{X_1, \dots, X_n\}$ is a~basis of $\fg$, we set $\widehat X_i:=X_i^*\in \fg^*$ or $\widehat X_i:=\Pi(X_i^*)\in \Pi(\fg^*)$
when cochains of degree $>1$ are considered.

\textit{A tough choice: what to compute?} The distinction
between the Lie algebra with Cartan matrix and its derived or quotient or simple
subquotient (such as $\fgl(pn)$, $\fsl(pn)$, $\fp\fgl(pn)$, and
$\fpsl(pn)$) was often disregarded. We insist on careful distinction, but are puzzled by the question: if $A$ is not invertible, so there are $\fg(A)$, $\fg(A)^{(i)}$, $\fg(A)/\fc$, $\fg(A)^{(i)}/\fc$,
\[
\text{should we compute $H^1(\fg;\fg)$ and $H^2(\fg)$ for all four relatives of $\fg(A)$?}
\]
We think not all these cases are usually interesting. Our choice is as follows: since the most often encountered are $\fg(A)$ and its simple subquotient $\fg(A)^{(i)}/\fc$, we should definitely consider these two cases. If some other relative proved to be especially nice, e.g., $\fsl(pn)$, we consider it as well; otherwise, we ignore the other two cases.

\subsection[Rank 1: Lie (super)algebras with Cartan matrices $(\od)$ or $(2)$ and $(\ev)$ (resp. $(1)$ and $(0)$). One more superization]{Rank 1: Lie (super)algebras with Cartan matrices $\boldsymbol{(\od)}$ or $\boldsymbol{(2)}$ and $\boldsymbol{(\ev)}$\\ (resp.\ $\boldsymbol{(1)}$ and $\boldsymbol{(0)}$). One more superization}\label{SSrk1Ort}

The case where size of the Cartan matrix $A$ is equal to $1$ is exceptional, some of its subcases have no analogs of Cartan matrices of larger size.

For $p=2$, these are
\begin{itemize}\itemsep=0pt
\item the simple Lie algebra
$\fo^{(1)}(3)$ with Cartan matrix $(\bar 1)$, and its superizations $\mathfrak{oo}^{(1)}_{I\Pi}(1|2)$ with Cartan matrix $(1)$, and $\mathfrak{oo}^{(1)}_{II}(1|2)$ without any Cartan matrix;
\item the Lie (super)algebra with Cartan matrix $(\ev)$ (resp.\ $(0)$) is $\fgl(2)$ (resp.\ $\fgl(1|1)$); the derived of this algebra is $\fhei(2|0)\simeq\fsl(2)$ (resp.\ $\fhei(0|2)\simeq \fsl(1|1)$);
\item there is also a~totally different superization of $\fhei(2|0)$, namely $\fb\fa(1|1)$, the \textit{anti-bracket} analog of the Heisenberg algebra with an odd center spanned by $\zeta$ and two more elements: even $a$ and odd $a^+$ with the only non-zero bracket $[a^+, a]=\zeta$. Same as $\fhei(2|0)$ (resp.\ $\fhei(0|2)$) can be realized as the negative part $\fg_-:=\oplus_{i<0}\fg_i$ of the Lie superalgebra of contact vector fields $\fk(3|0)$ (resp. $\fk(1|2)$), the Lie superalgebra $\fb\fa(1|1)$ can be realized as the negative part of the vectorial Lie superalgebra $\fg=\fm(1)$---the ``odd'' (a.k.a.\ pericontact) version of the contact algebra in its standard $\Zee$-grading, see~\eqref{st}. For details, see~\cite{BGLLS}; we recall only the necessary formulas.
\end{itemize}

For any $f\in \Kee[\tau, q, \xi]$, where $\tau$ and $\xi$ are odd, and $q$ is even, let
\begin{equation*}
M_f:=(2-E)(f)\frac{\partial}{\partial \tau}- (-1)^{p(f)} \frac{\partial f}{\partial \tau} E-\frac{\partial f}{\partial q}
\frac{\partial}{\partial \xi}-(-1)^{p(f)} \frac{\partial f}{\partial\xi} \frac{\partial}{\partial q},
\end{equation*}
where $E:=\xi \frac{\partial}{\partial \xi}+q\frac{\partial}{\partial q}$. Set
\begin{gather}\label{st}
\deg \tau=2,\qquad \deg q=\deg\xi=1,\qquad \deg M_f=\deg f-2, \qquad p(M_f)=p(f)+\od.
\end{gather}
Clearly,
\[
\fb\fa(1|1)=\operatorname{Span}\big(\zeta=M_1, a=M_\xi, a^+=M_q\big).
\]

\begin{Lemma}\label{ba}\quad
\begin{enumerate}\itemsep=0pt
\item[$1)$] For $\fg=\fb\fa(1|1)$, for any $p$ we have: $H^1(\fg;\fg)$ is spanned by
\[
 a^+\otimes\widehat a,\qquad
 a\otimes\widehat a~+ a^+\otimes\widehat{a^+},\qquad
 a\otimes\widehat a~- \zeta\otimes\widehat\zeta.
\]
Indeed: $\fout \fb\fa(1|1)=\operatorname{Span}(M_\tau, M_{q\xi}, M_{q^2})$.
\item[$2)$] The space $H^2(\fg)$ is spanned by
\begin{gather*}
 c_2 = \widehat{a^+} \wedge \widehat{a^+}, \qquad
 \underline{c_3} = {} \widehat a\wedge \widehat\zeta \quad \text{$($both for any $p)$},\nonumber\\
 c_4 = \begin{cases}
 0 & \text{for $p\neq 2$,}\\
 \widehat\zeta\wedge\widehat\zeta & \text{for $p=2$}. 
 \end{cases}
\end{gather*}
\end{enumerate}
\end{Lemma}

\begin{Remark} Assuming that $\fb\fa(1|1)=\fm(1)_-$, the negative part of $\fm(1)$, is an analog of $\fsl(2)$, or rather $\fhei(2|0)$, we can say that the non-positive part $\fm(1)_{\leq 0}$ of $\fm(1)$ is a~(far fetched and ad hoc) analog of $\fgl(2)$. What are the derivations and central extensions of $\fm(1)_{\leq 0}$?

Answer: For $p=2$, the space $H^2(\fg)$ is spanned by
\[
c_{0}^1={} \widehat{ M_\tau}\wedge \widehat{ M_{q\xi}},
\qquad c_{0}^2 ={} \widehat{ M_{q^{(2)}}}\wedge\widehat{ M_{q^{(2)}}}
\qquad c_{2} ={} \widehat{ M_q}\wedge \widehat{ M_q},
\qquad c_{4} ={} \widehat{ M_1}\wedge \widehat{ M_1}.\]
For $p\neq2$, the space $H^2(\fg)$ is spanned by
\[
c_{0}= {} \widehat{ M_\tau}\wedge \widehat{ M_{q\xi}}.\]

The space of $H^1(\fg;\fg)$ is spanned by $M_1\otimes \widehat{M_{q^{(2)}}}$ if $p=2$, otherwise $H^1(\fg;\fg)=0$.
\end{Remark}

For $p> 2$, the simple Lie algebra $\fo(3)\simeq \fsl(2)$ has Cartan matrix $(2)$ and its super versions $\fsl(1|1)$, or rather $\fgl(1|1)$, has Cartan matrix $(0)$, whereas $\fosp(1|2)$ has Cartan matrix $(1)$.

\begin{Lemma}\label{2.1.1b} For any $p$ and $\fg=\fhei(2|0)$:
\begin{enumerate}\itemsep=0pt
\item[$1)$] We have $\fout\ \fg\simeq\fgl(2)$; for a~basis of $H^1(\fg;\fg)$ we can take the following derivations $($recall convention~\eqref{conv}$)$
\begin{equation*}
c_{-2}= y_1\otimes\widehat x_1,\qquad
 c_{0}^1= h_1\otimes\widehat h_1+x_1\otimes\widehat x_1,\qquad
c_{0}^2= -x_1\otimes\widehat x_1+y_1\otimes\widehat y_1.
\end{equation*}
\item[$2)$] We have $H^2(\fg)=\operatorname{Span}\big(\widehat x_1\wedge\widehat h_1, \widehat y_1\wedge\widehat h_1\big)$.
\end{enumerate}
\end{Lemma}

\begin{Comment}\label{2.1.1c} Let us show that the Lie algebra
$\fout\ \fg$, is not ``too big'',
i.e., Lemma~\ref{2.1.1b} is correct. The basis of the space of
0-degree 1-cochains is $\big\{h_1\otimes\widehat h_1, x_1\otimes\widehat x_1,
y_1\otimes\widehat y_1\big\}$. We have
\begin{equation}\label{comb}
d\big(a{h}_1\otimes\widehat h_1+bx_1\otimes\widehat x_1+c y_1\otimes\widehat y_1\big)=(a+b+c)
h_1\otimes\widehat x_1\wedge\widehat y_1 .
\end{equation}
For the linear combination to be a~cocycle, the expression
\eqref{comb} should vanish, i.e., $a+b+c=0$. Besides, there are
coboundaries of degree 0 because $d({h}_1)=0$. Now, \textsc{SuperLie} is
taking a~basis with $a=b=1$, $c=0$ and $a=0$, $b=c=1$. To have a
symmetric answer, we'd take a~basis with $a=b=1$, $c=0$ and $a=c=1$,
$b=0$, but the code \textsc{SuperLie} has different esthetic
criteria.
\end{Comment}

\begin{Lemma}\label{2.1.1d} For any $p$
and $\fg=\fhei(0|2)$, for a~basis of
$H^1(\fg;\fg)$ we can take the following derivations:
\begin{equation*}
c_{0}^1=h_1\otimes\widehat h_1-x_1\otimes\widehat x_1,\qquad
c_{0}^2=-x_1\otimes\widehat x_1+y_1\otimes\widehat y_1,
\end{equation*}
and hence
$\fout \fg\simeq\begin{cases}\fo_\Pi(2)\oplus\Kee z&\text{for $p\neq 2$,}\\
\fo_\Pi^{(1)}(2)\oplus\Kee z&\text{for $p= 2$.}\end{cases}$

We have $H^2(\fg)=\operatorname{Span}\big(\widehat x_1\wedge\widehat h_1, \widehat y_1\wedge\widehat h_1\big)$.
\end{Lemma}

\begin{Lemma}\label{2.1.2b} Let $\fg=\fosp(1|2)$. Let $p>3$. Then, $H^1(\fg;\fg)=0$ and $H^2(\fg)=0$.

Let $p=3$.
\begin{enumerate}\itemsep=0pt
\item[$a)$] For $\fg=\fosp(1|2)$, for a~basis
of $H^1(\fg;\fg)$ we can take the following odd derivations $($recall convention~\eqref{conv} and that $x_2:=x_1^2$, $y_2:=y_1^2)$
\[
\underline{c_{-3}}= - y_1\otimes \widehat x_2+y_2\otimes\widehat x_1.
\]
\item[$b)$] $H^2(\fg)=0$.
\end{enumerate}
\end{Lemma}

Let $p=2$.
(Recall that for $p=2$, there is no $\fosp(1|2)$, there are two non-isomorphic Lie superalgebras $\fo\fo_B(1|2)$ for $B=I\Pi\sim\Pi\Pi$ and $B=II\sim \Pi I$.)

\begin{Lemma}\label{2.1.2c}
For $\fg=\fo\fo_{I\Pi}^{(1)}(1|2)$ and $\fo\fo_{I\Pi}(1|2)$, as well as for $\fo\fo_{II}^{(1)}(1|2)$ and $\fo\fo_{II}(1|2)$, we have
$H^1(\fg;\fg)=0$ and $H^2(\fg)=0$.
\end{Lemma}

One might expect the same result for $\fo\fo_{I\Pi}^{(1)}(1|2)$ as for
$\fo^{(1)}_{\Pi}(3)$. Lemma~\ref{2.1.2c} is, however, correct: let us
prove it.

\begin{proof} There is one coboundary: $x_1\otimes\widehat x_1+y_1\otimes
\widehat y_1$. The space of cochains is 5-dimensional, that of cocycles is
equal to 2. For $\fo^{(1)}_{\Pi}(3)$, the situation is different
(first, the dimensions are smaller): the space of cochains is
3-dimensional, and that of cocycles is equal to 1. The coboundary,
which is also $x_1\otimes\widehat x_1+y_1\otimes\widehat y_1$, reduces the answer to 0.
\end{proof}

For $\fg=\fsl(2)$, we have $H^1(\fg;\fg)=0$ and $H^2(\fg)=0$ for any $p\neq 2$, see Section \ref{known}.

\subsection[Rank 2: $\fsl(3)$ and $\fgl(3)$ for $p=3$; $\fsl(1|2)$ for any $p$; $\fbr(2;\eps)$ and $\mathfrak{brj}(2;3)$ for $p=3$;
$\mathfrak{brj}(2;5)$ for $p=5$]{Rank 2: $\boldsymbol{\fsl(3)}$ and $\boldsymbol{\fgl(3)}$ for $\boldsymbol{p=3}$; $\boldsymbol{\fsl(1|2)}$ for any $\boldsymbol{p}$; \\ $\boldsymbol{\fbr(2;\eps)}$ and $\boldsymbol{\mathfrak{brj}(2;3)}$ for $\boldsymbol{p=3}$;
$\boldsymbol{\mathfrak{brj}(2;5)}$ for $\boldsymbol{p=5}$}

\begin{Lemma}\label{2.2.1b} Let $p\neq 2$ and $\fg=\fsl(1|2)$. We have
$H^1(\fg;\fg)=0$ and $H^2(\fg)=0$.

For $p=2$, we have
$H^1(\fg;\fg)=0$ and the space $H^2(\fg)$ is spanned by $\widehat{x}_1\wedge\widehat{x}_1$ and
 $\widehat{x}_3\wedge\widehat{x}_3$.
\end{Lemma}

\begin{proof} For $p\neq 2$: same arguments as in \cite{Dz}; for $p=2$ we used \textsc{SuperLie}. \end{proof}

\begin{Conjecture} For $p=2$, if $x\in\fg_\od$ and $x^2=0$, then $\widehat{x}\wedge\widehat{x}$
is a~nontrivial cocycle, see cases $\fgl(a|a+pk)$ in Section~$\ref{sssAK}$.
\end{Conjecture}

\begin{Lemma}\label{2.2.2a} Let $p=3$.
\begin{enumerate}\itemsep=0pt
\item[$(a1)$] For
$\fg=\fgl(3)$, for a~basis of $H^1(\fg;\fg)$ we can take the
following derivation:
\[
c_{0}=2 h_1\otimes \widehat D_3
 +h_2\otimes \widehat D_3,
 \]
where $D_3$ denotes the grading operator, see Remark $\ref{RemD}$.
\item[$(a2)$] We have $H^2(\fgl(3))=0$.
\item[$(b1)$] See {\rm \cite{DzI,Ib,Ib2}}. For $\fg=\fpsl(3)$, for a~basis of $H^1(\fg;\fg)$ we can
take the following $7$ derivations $($recall convention~\eqref{conv}$)$
\begin{alignat*}{3}
& \deg=-3\colon \quad && c_{-3}^1= y_1\otimes \widehat x_3+y_3\otimes \widehat x_1,\qquad
c_{-3}^2= y_2\otimes \widehat x_3+y_3\otimes \widehat x_2,& \\
& \deg=0\colon \quad && c_{0}^1= 2 x_2\otimes \widehat x_1+y_1\otimes \widehat y_2,\qquad c_{0}^2= 2 x_1\otimes \widehat x_2+y_2\otimes \widehat y_1,&\\
&&& c_{0}^3= x_1\otimes \widehat x_1 +2 x_2\otimes \widehat x_2 +2
y_1\otimes \widehat y_1 +y_2\otimes \widehat y_2.&
\end{alignat*}
\item[$(b2)$] See {\rm \cite{Ib,Ib2}}. We have $\dim H^2(\fpsl(3))=7$, for a~basis we can take the following cocycles $($recall convention~\eqref{conv}$)$
\begin{alignat*}{3}
&\deg=-3\colon \quad &&c_{-3}^1=\widehat x_2\wedge\widehat x_3, \qquad c_{-3}^2=\widehat x_1\wedge\widehat x_3;&\\
&\deg=0\colon \quad &&c_{0}^1=\widehat x_1\wedge\widehat y_2, \qquad c_{0}^2=\widehat x_2\wedge\widehat y_2-\widehat x_1\wedge\widehat y_1, \qquad c_{0}^3=\widehat x_2\wedge\widehat y_1.
\end{alignat*}
\end{enumerate}
\end{Lemma}

\subsubsection{Clarifying an incredible result}\label{2.2.2aa}
We precede the proof of Lemma~\ref{2.2.2a} with these remarks: they are for those who, like the authors themselves at first, can not believe the
claims in headings a) of Lemma~\ref{2.2.2a} are true.

(a1) It is doubtful because ``$\fgl$ can not
have outer derivations by definition". The expression in
quotation marks is not true unless $p=0$, see Statement~
\ref{Statement}. Since the Cartan matrix of $\fgl(3)$ for $p=3$ is
not invertible, there is a~grading operator $D$, see \cite{BGL1}.
For such a~$D$ one can take, e.g., $D_3=2 E_{1,1}+E_{2,2}+E_{3,3}$, or simply $E_{1,1}$, both of
which, of course, \textit{belong to the algebra, but are outer derivations from the cohomological point of view}.

(a2) Having read the arXiv version of this paper, A.~Dzhumadildaev and P.~Zusmanovich wrote to us that they doubted claims of item (a2) of Lemma~\ref{2.2.2a} because of the following fact (due to J.~Dixmier, and in a~more general setting of a~\textit{subalgebra} $\mathfrak{i}$, not an \textit{ideal}, due to A.~Dzhumadildaev; for a~proof \textit{in characteristic $0$} and references, see~\cite{Zus}):

\begin{Fact}\label{Fact} Let $\fg=\mathfrak{i} \, {\subset\hspace{-3.5mm} +}\, \Kee x$. Then,
\begin{gather}\label{wrong}
H^m(\fg;V)=H^m(\mathfrak{i}; V) \oplus H^{m-1}(\mathfrak{i};V)\wedge x.
\end{gather}
\end{Fact}

Proof of equality~\eqref{wrong} in the particular case $\mathfrak{i}=\fsl(n)\subset \fg=\fgl(n)$ is supposed to be as follows. Let~$i_x$ be the inner multiplication of cochains by $x\in\fgl(n)$, and $\psi$ a~degree $m-1$ cocycle of~$\fsl(n)$ with values in $V$. Since the trace $\operatorname{tr}$ is a~1-cocycle of~$\fgl(n)$, it follows that $\operatorname{tr}\wedge\psi$ is a~cocycle of~$\fgl(n)$ of degree~$m$.
Therefore, (since the codifferential is a~derivation) any degree~$m$ cocycle~$\alpha$ of~$\fgl(n)$ can be decomposed into a~sum of cocycles
\[
\alpha=\beta+\gamma, \qquad \text{where $i_x\beta=0$ and $i_x\gamma$ induces a~degree $m-1$ cocycle of $\fsl(n)$.}
\]
Therefore, $\alpha$ induces a~degree $m$ cocycle of $\fsl(n)$.
Hence, $\alpha = \operatorname{tr} \wedge i_{1_n}\gamma$.

However, this proof works if $x=1_n$ complements $\fsl(n)$ to $\fgl(n)$, \textit{but fails if $p$ divides $n$}.

\begin{proof}[Proof of Lemma~$\ref{2.2.2a}$] Recall that $p=3$.

(a1) The grading $\deg(\text{Chevalley generator})=\pm1$ of both $\fsl$ and $\fgl$ yields a~grading of depth~4 of $C^1(\fsl(3); \fgl(3))$. No 1-cocycle exists in degrees $-4$ and $-3$. In degree $-2$, we have just one 1-cocycle (of weight $(-1,-1,-1)$):
\[
2y_2\otimes \widehat x_1+h_1\otimes \widehat x_3+h_2\otimes x_3 +y_3\otimes \widehat h_1+y_3\otimes \widehat h_2+y_1\otimes \widehat x_2.
\]
This 1-cocycle is trivial, because it is just $d(y_3)$, where $y_3=[y_1, y_2]$. In degree $-1$ we have the following two cocycles
\begin{alignat*}{3}
& (-2,1,-1)\colon \quad && 2 x_2\otimes \widehat x_3 +2 y_3\otimes \widehat y_2 +h_1\otimes \widehat x_1+y_1\otimes \widehat h_1+y_1\otimes \widehat h_2,&\\
& (1,-2,0) \colon \quad & & h_2\otimes\widehat x_2 +
y_2\otimes \widehat h_1+y_2\otimes \widehat h_2 +x_1\otimes \widehat x_3 +y_3\otimes \widehat y_1.&
 \end{alignat*}
These 1-cocycles are trivial because they are the coboundaries of $y_1$ and $y_2$, respectively. In degree 0, we have the following 1-cocycles
\begin{alignat*}{3}
& (0,0,0)\colon \quad &&
2 x_2\otimes \widehat x_2 +2y_1\otimes \widehat y_1 +x_1\otimes \widehat x_1+y_2\otimes \widehat y_2,& \\
& (0,0,0)\colon \quad && 2x_3\otimes \widehat x_3 +2x_1\otimes \widehat x_1+y_1\otimes \widehat y_1+y_3\otimes \widehat y_3.&
\end{alignat*}
These cocycles are coboundaries of $h_2-h_3$ and $2h_3$, respectively.
Hence, $H^1(\fsl(3); \fgl(3))=0$.

(a2) Observe that $H^2(\fsl(3))$ is spanned by (the classes of) the following cocycles:
\begin{alignat*}{3}
& \deg=-3\colon \quad &&\widehat x_1\wedge\widehat x_3, \qquad \widehat x_2\wedge\widehat x_3,&\\
&\deg=0\colon \quad &&\widehat x_1\wedge\widehat y_2, \ \ \ \widehat x_2\wedge\widehat y_1.&
\end{alignat*}

Let us verify that $H^2(\fgl(3))=0$, although $H^2(\fsl(3))\neq 0$. It is easy to check manually that the following list exhausts all 2-cochains of $\fsl(3)$ and of $\fgl(3)$
with trivial coefficients:
\begin{alignat*}{5}
& (-2,1,1)\colon \quad && \widehat E_{1,2}\wedge \widehat E_{1,3},\qquad && (1,-2,1)\colon \quad && \widehat E_{2,3}\wedge \widehat E_{2,1},& \\
& (-1,-1,2)\colon \quad &&\widehat E_{2,3}\wedge \widehat E_{1,3},\qquad && (1,1,-2)\colon \quad && \widehat E_{3,2}\wedge \widehat E_{3,1}, &\\
& (-1,2,-1)\colon \quad && \widehat E_{1,2}\wedge \widehat E_{3,2},\qquad && (2,-1,-1)\colon \quad && \widehat E_{2,1}\wedge \widehat E_{3,1}.&
\end{alignat*}
For $\fgl(3)$, there are no 1-cochains of these weights, so it suffices to check if the above cochains are closed cocycles. For example
\[
d\big(\widehat E_{1,2}\wedge \widehat E_{1,3}\big) =-\widehat 1_{3}\wedge
\widehat E_{1,2}\wedge \widehat E_{1,3}.
\]
Evidently, this is not zero on $\fgl(3)$, but vanishes on $\fsl(3)$.
\end{proof}

\subsection[Lie superalgebra $\fgl(a|a+kp)$ and its simple relatives (A.~Krutov)]{Lie superalgebra $\boldsymbol{\fgl(a|a+kp)}$ and its simple relatives (A.~Krutov)}\label{sssAK} Let us collect here all these cases of rank $\leq 8$.
 For the elements of Chevalley basis, see \cite[Section~9.2]{BGL2}.

Let $p=5$:
For $\fg=\fgl(5)$, the space $H^1(\fg;\fg)$ is spanned by
\[
c_0 = (4 h_1 + 3 h_2 + 2 h_3 + h_4)\otimes\widehat D
\]
and $H^2(\fg)=0$.

For $\fg=\fpsl(5)$, the space $H^1(\fg;\fg)$ is spanned by
\begin{gather*}
 c_0 =
 4 x_2\otimes\widehat x_2+x_3\otimes\widehat x_3+4 x_5\otimes\widehat x_5+x_7\otimes\widehat x_7+y_2\otimes\widehat y_2+4 y_3\otimes\widehat y_3+y_5\otimes\widehat y_5+4 y_7\otimes\widehat y_7
\end{gather*}
and the space $H^2(\fg)$ is spanned by
\[
 c_0 = \widehat x_2\wedge \widehat y_2+4 \widehat x_3\wedge \widehat y_3+4 \widehat x_5\wedge \widehat y_5+\widehat x_7\wedge \widehat y_7.
\]

For $\fpsl(1|6)$, we see that $H^2(\fg)$ is spanned by
\[
 c_0 = \widehat x_4\wedge \widehat y_4+4 \widehat x_5\wedge \widehat y_5+4 \widehat x_9\wedge \widehat y_9+\widehat x_{11}\wedge \widehat y_{11}+\widehat x_{13}\wedge \widehat y_{13}+\widehat x_{16}\wedge \widehat y_{16}
\]
and $H^1(\fg;\fg)$ is spanned by
\begin{gather*} c_0 = x_2\otimes\widehat x_2+4 x_3\otimes\widehat x_3+4 x_7\otimes\widehat x_7+4 x_9\otimes\widehat x_9+4 x_{14}\otimes\widehat x_{14}+4 x_{18}\otimes\widehat x_{18}\\
\hphantom{c_0 =}{}
+4 y_2\otimes\widehat y_2+y_3\otimes\widehat y_3+y_7\otimes\widehat y_7+y_9\otimes\widehat y_9+y_{14}\otimes\widehat y_{14}+y_{18}\otimes\widehat y_{18}.
\end{gather*}

For $\fg=\fgl(2|7)$, we see that $H^2(\fg)=0$ and $H^1(\fg;\fg)$ is spanned by
\[
c_0 = (4h_1+3h_2+h_3+4h_5+3h_6+2h_7+h_8)\otimes\widehat D.
\]

For $\fg=\fpsl(2|7)$, we see that $H^2(\fg)$ is spanned by
\begin{gather*}
 c_0 = 4 \widehat x_5\wedge \widehat y_5+\widehat x_6\wedge \widehat y_6+\widehat x_{12}\wedge \widehat y_{12}+4 \widehat x_{14}\wedge \widehat y_{14}+4 \widehat x_{18}\wedge \widehat y_{18}+\widehat x_{21}\wedge \widehat y_{21}\\
 \hphantom{c_0 =}{} +4 \widehat x_{23}\wedge \widehat y_{23}+\widehat x_{27}\wedge \widehat y_{27},
\end{gather*}
and $H^1(\fg;\fg)$ is spanned by
\begin{gather*} c_0 = x_4\otimes\widehat x_4+4 x_5\otimes\widehat x_5+x_{11}\otimes\widehat x_{11}+4 x_{13}\otimes\widehat x_{13}+4 x_{17}\otimes\widehat x_{17}+4 x_{20}\otimes\widehat x_{20}\\
\hphantom{c_0 =}{}
+4 x_{22}\otimes\widehat x_{22}+4 x_{26}\otimes\widehat x_{26}+4 y_4\otimes\widehat y_4+y_5\otimes\widehat y_5+4 y_{11}\otimes\widehat y_{11}+y_{13}\otimes\widehat y_{13}
\\
\hphantom{c_0 =}{}+y_{17}\otimes\widehat y_{17}+y_{20}\otimes\widehat y_{20}+y_{22}\otimes\widehat y_{22}+y_{26}\otimes\widehat y_{26}.
\end{gather*}

Let $p=3$:

For $\fgl(3)$ and $\fpsl(3)$, see Section \ref{2.2.2a}.

For $\fgl(6)$ and $\fpsl(6)$, see Section \ref{gl4_22}.

For $\fg=\fgl(1|4)$, the space $H^1(\fg;\fg)$ is spanned by
\[
c_0 = 2 h_1\otimes\widehat D+2 h_3\otimes\widehat
D+h_4\otimes\widehat D
\]
and $H^2(\fg) = 0$.

For $\fg=\fpsl(1|4)$, the space $H^1(\fg;\fg)$ is spanned by (see Lemma~\ref{2.4.2a})
\begin{gather*}
 c_{0} = x_2\otimes\widehat x_2+2 x_3\otimes\widehat x_3+2 x_5\otimes\widehat x_5+2 x_7\otimes\widehat x_7+2 y_2\otimes\widehat y_2+y_3\otimes\widehat y_3+y_5\otimes\widehat y_5+y_7\otimes\widehat y_7
\end{gather*}
and $H^2(\fg)$ is spanned by
\[
c_0 = 2 \widehat x_3\wedge \widehat y_3+\widehat x_4\wedge \widehat y_4+\widehat x_6\wedge \widehat y_6+\widehat x_8\wedge \widehat y_8.
\]

For $\fg=\fgl(1|7)$, we see that $H^2(\fg)=0$ and $H^1(\fg;\fg)$ is spanned by
\[
c_0 = (2h_1+2h_3+h_4+2h_6+h_7)\otimes\widehat D.
\]

For $\fg=\fpsl(1|7)$, we see that $H^2(\fg)$ is spanned by
\begin{gather*}
 c_{0} =
 2 \widehat x_4\wedge \widehat y_4+\widehat x_5\wedge \widehat y_5+\widehat x_{10}\wedge \widehat y_{10}+2 \widehat x_{12}\wedge \widehat y_{12}+2 \widehat x_{15}\wedge \widehat y_{15}+\widehat x_{18}\wedge \widehat y_{18}+2 \widehat x_{19}\wedge \widehat y_{19},
\end{gather*}
and $H^1(\fg;\fg)$ is spanned by
\begin{gather*}
 c_{0} = x_3\otimes\widehat x_3+2 x_4\otimes\widehat x_4+x_9\otimes\widehat x_9+2 x_{11}\otimes\widehat x_{11}+2 x_{14}\otimes\widehat x_{14}+2 x_{17}\otimes\widehat x_{17}+2 x_{22}\otimes\widehat x_{22}\\
\hphantom{c_{0} =}{}
+2 y_3\otimes\widehat y_3+y_4\otimes\widehat y_4+2 y_9\otimes\widehat y_9+y_{11}\otimes\widehat y_{11}+y_{14}\otimes\widehat y_{14}+y_{17}\otimes\widehat y_{17}+y_{22}\otimes\widehat y_{22}.
 \end{gather*}

For $\fg=\fgl(2|5)$, we see that $H^2(\fg)=0$ and $H^1(\fg;\fg)$ is spanned by
\[
 c_0 = (2h_1+h_2+h_3+2h_5+h_6)\otimes\widehat D.
 \]

For $\fg=\fpsl(2|5)$, we see that $H^2(\fg)$ is spanned by
\[
 c_0 = \widehat x_4\wedge \widehat y_4+2 \widehat x_5\wedge \widehat y_5+2 \widehat x_9\wedge \widehat y_9+\widehat x_{11}\wedge \widehat y_{11}+2 \widehat x_{13}\wedge \widehat y_{13}+\widehat x_{16}\wedge \widehat y_{16},
\]
and $H^1(\fg;\fg)$ is spanned by
\begin{gather*}
 c_{0} = x_3\otimes\widehat x_3+2 x_4\otimes\widehat x_4+2 x_8\otimes\widehat x_8+2 x_{10}\otimes\widehat x_{10}+2 x_{12}\otimes\widehat x_{12}+2 x_{15}\otimes\widehat x_{15}\\
\hphantom{c_{0} =}{}
+2 y_3\otimes\widehat y_3+y_4\otimes\widehat y_4+y_8\otimes\widehat y_8+y_{10}\otimes\widehat y_{10}+y_{12}\otimes\widehat y_{12}+y_{15}\otimes\widehat y_{15}.
\end{gather*}

For $\fgl(3|3)$ and $\fpsl(3|3)$, see Section \ref{gl4_22}.

For $\fg=\fgl(3|6)$, we see that $H^2(\fg)=0$, and $H^1(\fg;\fg)$ is spanned by
\[
 c_0 = (2h_1+h_2+2h_4+h_5+2h_7+h_8)\otimes\widehat D.
 \]

For $\fg=\fpsl(3|6)$, we see that $H^2(\fg)$ is spanned by
\begin{gather*}
 c_{0} =
 2 \widehat x_1\wedge \widehat y_1+\widehat x_2\wedge \widehat y_2+2 \widehat x_{10}\wedge \widehat y_{10}+\widehat x_{17}\wedge \widehat y_{17}+2 \widehat x_{23}\wedge \widehat y_{23}+\widehat x_{28}\wedge \widehat y_{28}\\
 \hphantom{c_{0} = }{}
 +2 \widehat x_{32}\wedge \widehat y_{32}+\widehat x_{35}\wedge \widehat y_{35},
 \end{gather*}
and $H^1(\fg;\fg)$ is spanned by
\begin{gather*}
 c_{0} = 2 x_3\otimes\widehat x_3+2 x_4\otimes\widehat x_4+2 x_{10}\otimes\widehat x_{10}+2 x_{12}\otimes\widehat x_{12}+2 x_{16}\otimes\widehat x_{16}+2 x_{19}\otimes\widehat x_{19}\\
\hphantom{c_{0} =}{}
+2 x_{25}\otimes\widehat x_{25}+2 x_{30}\otimes\widehat x_{30}+y_3\otimes\widehat y_3+y_4\otimes\widehat y_4+y_{10}\otimes\widehat y_{10}+y_{12}\otimes\widehat y_{12}\\
\hphantom{c_{0} =}{}+y_{16}\otimes\widehat y_{16}+y_{19}\otimes\widehat y_{19}+y_{25}\otimes\widehat y_{25}+y_{30}\otimes\widehat y_{30}.
\end{gather*}

Let $p=2$:

For $\fgl(4)$ and $\fpsl(4)\simeq \textbf{F}(\fh_\Pi^{(1)}(0|4))$ (for this exceptional case, see Sections \ref{gl4_22} and~\ref{psl22}).

For $\fgl(6)$ and $\fpsl(6)$, see Section~\ref{2.5.1b}.

For $\fg=\fpsl(8)$, the space $H^2(\fg)$ is spanned by
\[
 c_0 = \widehat x_4\wedge \widehat y_4+\widehat x_5\wedge \widehat y_5+\widehat x_{10}\wedge \widehat
 y_{10}+\widehat x_{12}\wedge \widehat y_{12}+\widehat x_{15}\wedge \widehat y_{15}+\widehat x_{18}\wedge \widehat y_{18}+\widehat x_{19}\wedge \widehat y_{19},
\]
and $H^1(\fg;\fg)$ is spanned by
\begin{gather*}
 c_{0} =
 x_4\otimes\widehat x_4+x_5\otimes\widehat x_5+x_{10}\otimes\widehat x_{10}+x_{12}\otimes\widehat x_{12}+x_{15}\otimes\widehat x_{15}+x_{18}\otimes\widehat x_{18}+x_{19}\otimes\widehat x_{19}\\
\hphantom{c_{0} =}{}
+y_4\otimes\widehat y_4+y_5\otimes\widehat y_5+y_{10}\otimes\widehat y_{10}+y_{12}\otimes\widehat y_{12}+y_{15}\otimes\widehat y_{15}+y_{18}\otimes\widehat y_{18}+y_{19}\otimes\widehat y_{19}.
\end{gather*}

For $\fg=\fgl(1|3)$, we see that $H^2(\fg)$ is spanned by $\widehat x_i\wedge\widehat x_i$ for all $i$ such that $x_i$ is odd, and
\[
H^1(\fg;\fg)=\Kee\big[(h_1+h_3)\otimes\widehat D\big].
\]

For $\fg = \fpsl(1|3)$, we see that $H^2(\fg)$ is spanned by (in the standard format of $\fgl(1|3)$, the only simple root, the one corresponding to $x_1$, is odd; the other odd positive root vectors are~$x_4$ and~$x_6$:
\begin{gather*}
c_{-6} = \widehat x_6\wedge \widehat x_6,\qquad
c_{-4,1} = \widehat x_4\wedge \widehat x_4,\qquad
c_{-4,2} = \widehat x_2\wedge \widehat x_6+\widehat x_4\wedge \widehat x_5,\qquad
c_{-2,1} = \widehat x_1\wedge \widehat x_1,\\
c_{-2,2} = \widehat x_6\wedge \widehat y_2+\widehat x_1\wedge \widehat x_3,\qquad
c_{0,1} = \widehat x_1\wedge \widehat y_3+\widehat x_4\wedge \widehat y_5,\qquad
c_{0,2} = \widehat x_3\wedge \widehat y_1+\widehat x_5\wedge \widehat y_4,\\
c_{0,3} = \widehat x_2\wedge \widehat y_2+\widehat x_3\wedge \widehat y_3+\widehat x_4\wedge \widehat y_4.
\end{gather*}

For $\fgl(1|5)$ and $\fpsl(1|5)$, see Section \ref{2.5.1b}.

For $\fg=\fgl(1|7)$, we see that $H^2(\fg)$ is spanned by $\widehat x_i\otimes\widehat x_i$ for all $i$ such
that $x_i$ is odd, and
\[
H^1(\fg;\fg)=\Kee\big[(h_1+h_3+h_5+h_7)\otimes\widehat D\big].
\]

For $\fg=\fpsl(1|7)$, we see that $H^2(\fg)$ is spanned by
\begin{gather*}
\widehat x_i\wedge \widehat x_i\qquad\text{for all $i$ such that $x_i$ is odd},\\
 c_{0} = \widehat x_4\wedge \widehat y_4+\widehat x_5\wedge \widehat y_5+\widehat x_{10}\wedge \widehat y_{10}+\widehat x_{12}\wedge \widehat y_{12}+\widehat x_{15}\wedge \widehat y_{15}+\widehat x_{18}\wedge \widehat y_{18}+\widehat x_{19}\wedge \widehat y_{19},
\end{gather*}
and $H^1(\fg;\fg)$ is spanned by
\begin{gather*}
c_0 = x_4\otimes\widehat x_4+x_5\otimes\widehat x_5+x_{10}\otimes\widehat x_{10}+x_{12}\otimes\widehat x_{12}+x_{15}\otimes\widehat x_{15}+x_{18}\otimes\widehat x_{18}+x_{19}\otimes\widehat x_{19}\\
\hphantom{c_0 =}{}
+y_4\otimes\widehat y_4+y_5\otimes\widehat y_5+y_{10}\otimes\widehat y_{10}+y_{12}\otimes\widehat y_{12}+y_{15}\otimes\widehat y_{15}+y_{18}\otimes\widehat y_{18}+y_{19}\otimes\widehat y_{19}.
\end{gather*}

For $\fgl(2|2)$ and $\fpsl(2|2)\simeq\fh_\Pi^{(1)}(0|4)$ (for this exceptional case, see Sections~\ref{gl4_22} and~\ref{psl22}).

For $\fgl(2|4)$ and $\fpsl(2|4)$, see Section~\ref{2.5.1b}.

For $\fg=\fgl(2|6)$, we see that $H^2(\fg)$ is spanned by $\widehat x_i\otimes\widehat x_i$ for all $i$ such
that $x_i$ is odd, and
\[
H^1(\fg;\fg)=\Kee\big[(h_1+h_3+h_5+h_7)\otimes\widehat D\big].
\]

For $\fg= \fpsl(2|6)$, we see that $H^2(\fg)$ is spanned by
\begin{gather*}
 \widehat x_i\wedge \widehat x_i\qquad\text{for all $i$ such that $x_i$ is odd},\\
 c_{0} = \widehat x_4\wedge \widehat y_4+\widehat x_5\wedge \widehat y_5+\widehat x_{10}\wedge \widehat y_{10}+\widehat x_{12}\wedge \widehat y_{12}+\widehat x_{15}\wedge \widehat y_{15}+\widehat x_{18}\wedge \widehat y_{18}+\widehat x_{19}\wedge \widehat y_{19},
\end{gather*}
and $H^1(\fg;\fg)$ is spanned by
\begin{gather*}
 c_{0} =
 x_4\otimes\widehat x_4+x_5\otimes\widehat x_5+x_{10}\otimes\widehat x_{10}+x_{12}\otimes\widehat x_{12}+x_{15}\otimes\widehat x_{15}+x_{18}\otimes\widehat x_{18}+x_{19}\otimes\widehat x_{19}\\
\hphantom{c_{0} =}{}
 +y_4\otimes\widehat y_4+y_5\otimes\widehat y_5+y_{10}\otimes\widehat y_{10}+y_{12}\otimes\widehat y_{12}+y_{15}\otimes\widehat y_{15}+y_{18}\otimes\widehat y_{18}+y_{19}\otimes\widehat y_{19}.
 \end{gather*}

For $\fg=\fgl(3|5)$, we see that $H^2(\fg)$ is spanned by $\widehat x_i\otimes\widehat x_i$ for all $i$ such
that $x_i$ is odd, and
\[
H^1(\fg;\fg)=\Kee\big[(h_1+h_3+h_5+h_7)\otimes\widehat D\big].
\]

For $\fg = \fpsl(3|5)$, we see that $H^2(\fg)$ is spanned by
\begin{gather*}
\widehat x_i\wedge \widehat x_i\qquad\text{for all $i$ such that $x_i$ is odd},\\
 c_{0} = \widehat x_4\wedge \widehat y_4+\widehat x_5\wedge \widehat y_5+\widehat x_{10}\wedge \widehat y_{10}+\widehat x_{12}\wedge \widehat y_{12}+\widehat x_{15}\wedge \widehat y_{15}+\widehat x_{18}\wedge \widehat y_{18}+\widehat x_{19}\wedge \widehat y_{19},
\end{gather*}
and $H^1(\fg;\fg)$ is spanned by
\begin{gather*}
 c_{0} =
 x_4\otimes\widehat x_4+x_5\otimes\widehat x_5+x_{10}\otimes\widehat x_{10}+x_{12}\otimes\widehat x_{12}+x_{15}\otimes\widehat x_{15}+x_{18}\otimes\widehat x_{18}+x_{19}\otimes\widehat x_{19}\\
\hphantom{c_{0} =}{}
+y_4\otimes\widehat y_4+y_5\otimes\widehat y_5+y_{10}\otimes\widehat y_{10}+y_{12}\otimes\widehat y_{12}+y_{15}\otimes\widehat y_{15}+y_{18}\otimes\widehat y_{18}+y_{19}\otimes\widehat y_{19}.
\end{gather*}

For $\fg=\fgl(4|4)$, we see that $H^2(\fg)$ is spanned by $\widehat x_i\otimes\widehat x_i$ for all $i$ such
that $x_i$ is odd and
\[
H^1(\fg;\fg)=\Kee\big[(h_1+h_3+h_5+h_7)\otimes\widehat D\big].
\]

For $\fg=\fpsl(4|4)$, we see that $H^2(\fg)$ is spanned by
\begin{gather*}
 \widehat x_i\wedge \widehat x_i\qquad\text{for all $i$ such that $x_i$ is odd},\\
 c_{0} = \widehat x_4\wedge \widehat y_4+\widehat x_5\wedge \widehat y_5+\widehat x_{10}\wedge \widehat y_{10}+\widehat x_{12}\wedge \widehat y_{12}+\widehat x_{15}\wedge \widehat y_{15}+\widehat x_{18}\wedge \widehat y_{18}+\widehat x_{19}\wedge \widehat y_{19},
\end{gather*}
and $H^1(\fg;\fg)$ is spanned by
\begin{gather*}
 c_0 = x_4\otimes\widehat x_4+x_5\otimes\widehat x_5+x_{10}\otimes\widehat x_{10}+x_{12}\otimes\widehat x_{12}+x_{15}\otimes\widehat x_{15}+x_{18}\otimes\widehat x_{18}+x_{19}\otimes\widehat x_{19}\\
\hphantom{c_0 =}{}
+y_4\otimes\widehat y_4+y_5\otimes\widehat y_5+y_{10}\otimes\widehat y_{10}+y_{12}\otimes\widehat y_{12}+y_{15}\otimes\widehat y_{15}+y_{18}\otimes\widehat y_{18}+y_{19}\otimes\widehat y_{19}.
\end{gather*}

\begin{Lemma}\label{2.2.2b} Let $p=3$. Let $\fg=\mathfrak{brj}(2;3)$ with Cartan matrix
\[
\begin{pmatrix}
 \hphantom{-}0& -1 \\
 -2 & \hphantom{-}1 \end{pmatrix}
 \]
 and basis
\begin{gather*}
x_1,\ x_2,\
x_3=[x_1, x_2],\ x_4=[x_2, x_2],\ x_5=[x_2,[x_1, x_2]],\ x_6=[[x_1, x_2], [x_2, x_2]],\nonumber\\
x_7=[[x_2, x_2], [x_2, [x_1, x_2]]], \
x_8=[[x_2, [x_1, x_2]], [x_2, [x_1, x_2]]]. 
\end{gather*}
\begin{enumerate}\itemsep=0pt
\item[$(a1)$] For a~basis of $H^1(\fg;\fg)$ we can take the
following derivations $($recall convention~\eqref{conv}$)$
\[
c_{-3}= x_1\otimes\widehat x_6+x_3\otimes\widehat x_7+2 y_2\otimes
\widehat x_4+y_4\otimes\widehat x_2+2 y_6\otimes\widehat y_1+y_7\otimes\widehat y_3.
\]
\item[$(a2)$] We have $H^2(\fg)=0$.
\end{enumerate}
\end{Lemma}

\begin{Lemma}\label{br2eps}
Let $p=3$. Let $\fg=\mathfrak{br}(2;\eps)$, where $\eps\neq 0$, with Cartan matrix
\[
\begin{pmatrix}
\hphantom{-}2&-1\\
-2&1-\eps\end{pmatrix}
\]
and basis
\begin{equation*}
 x_1, \ x_2, \
x_3=[x_1, x_2],\ x_4=-\operatorname{ad}_{x_2}^2(x_1),
\end{equation*}
Then, $H^1(\fg;\fg)=0$ and $H^2(\fg)=0$.
\end{Lemma}

\begin{Lemma}\label{2.3.1a} Let $p> 3$. For $\fg=\fsl(1|2)$, $\fosp(3|2)$, and $\fosp(1|4)$, we have $H^1(\fg;\fg)=0$ and $H^2(\fg)=0$.
\end{Lemma}

\begin{proof} We use the arguments of \cite{Dz}.\end{proof}

\begin{Lemma}\label{2.3.1b5} Let $p=5$. For $\fg=\mathfrak{brj} (2;5)$, we have $H^1(\fg;\fg)=0$ and $H^2(\fg)=0$.
\end{Lemma}

\subsection[Rank 3: $\fpsl(4)$ and $\fgl(4)$ for
$p=2$; $\mathfrak{ag}(2)$ for $p>3$; $\fbr(3)$ and
$\fg(1,6)$, $\fg(2,3)$ and $\fg^{(1)}(2,3)/\fc$ for $p=3$; $\fwk(3;a)$
and $\fbgl(3;a)$ for $p=2$]{Rank 3: $\boldsymbol{\fpsl(4)}$ and $\boldsymbol{\fgl(4)}$ for
$\boldsymbol{p=2}$; $\boldsymbol{\mathfrak{ag}(2)}$ for $\boldsymbol{p>3}$; $\boldsymbol{\fbr(3)}$ and
$\boldsymbol{\fg(1,6)}$,\\
 $\boldsymbol{\fg(2,3)}$ and $\boldsymbol{\fg^{(1)}(2,3)/\fc}$ for $\boldsymbol{p=3}$; $\boldsymbol{\fwk(3;a)}$
and $\boldsymbol{\fbgl(3;a)}$ for $\boldsymbol{p=2}$}\label{ssR3}

\begin{Lemma}\label{2.3.1b} Let $p=2$.
Let $\fg=\fwk(3;\alpha)$ and
$\fg=\fbgl(3;\alpha)$ with Cartan matrix
\[
\begin{pmatrix}
 \ev & \alpha & 0 \\
 \alpha & \ev & 1 \\
 0 & 1 & \ev\end{pmatrix}
 \]
 and basis
\begin{gather}\label{wk3a}
 x_1,\ x_2,\ x_3,\
 x_4=[x_1, x_2],\ x_5=[x_2, x_3],\ x_6=[x_3,[x_1, x_2]],\ x_7=[[x_1, x_2], [x_2,x_3]].
 \end{gather}
\begin{enumerate}\itemsep=0pt
\item[$(a1)$] For a~basis of $H^1(\fg;\fg)$ we can take
the following derivation, where $D_4$ is an outer derivation $($see Remark $\ref{RemD})$:
\[
c_{0}=h_1\otimes \widehat D_4+\alpha h_3\otimes \widehat D_4.
\]
\item[$(a2)$] We have $H^2(\fg)=0$.
\item[$(b1)$] For $\fg=\fwk(3;a)/\fc$ and $\fbgl(3;a)/\fc$, we
have $\dim H^1\big(\fg^{(1)};\fg^{(1)}\big)=1$ and $\fder \fg^{(1)}=\fg$.
\item[$(b2)$] The space $H^2(\fg^{(1)})$ is spanned by
\[
c_0=\widehat x_1\wedge\widehat y_1+\alpha\widehat x_4\wedge\widehat y_4+\alpha\widehat x_6\wedge\widehat y_6+\alpha(\alpha+1)\widehat x_7\wedge\widehat y_7.
\]
\end{enumerate}
\end{Lemma}

\begin{Lemma}\label{gl4_22} \quad
\begin{enumerate}\itemsep=0pt
\item[$(a1)$] Let $p=2$. For $\fg=\fgl(4)$, for a
basis of $H^1(\fg;\fg)$ we can take the following derivation where~$D_4$ is an outer derivation $($see Remark $\ref{RemD})$:
\[
c_{0}=h_1\otimes \widehat D_4+h_3\otimes \widehat D_4.
\]
\item[$(a2)$] We have $H^2(\fgl(4))=0$.
\item[$(b1)$] Let $p>2$. For $\fg=\fgl(2|2)$, for a~basis of $H^1(\fg;\fg)$ we can take
the following derivation, where $D_4$ is an outer derivation $($see Remark $\ref{RemD})$:
\[
c_{0}=- h_1\otimes \widehat D_4 +(p-2) h_2\otimes \widehat D_4 + h_3\otimes \widehat D_4.
\]
\item[$(b2)$] Let $p>2$. We have $H^2(\fgl(2|2))=0$.
\item[$(c1)$] For any $p>2$, for
$\fg=\fgl(1|3)$, for a~basis of $H^1(\fg;\fg)$ we can take the following derivatives, where $D_4$ is an outer
derivation $($see Remark $\ref{RemD})$:
\[
 c_{0} = (p-3) h_1\otimes\widehat D_4 + 2h_2\otimes\widehat D_4 + h_3\otimes\widehat D_4
 + 2 D_4\otimes\widehat D_4.
\]
\item[$(c2)$] For any $p>2$, we have $H^2(\fgl(1|3))=0$.
\end{enumerate}
\end{Lemma}

\subsection[$\fpsl(2|2)\simeq\fh_\Pi^{(1)}(0|4)$]{$\boldsymbol{\fpsl(2|2)\simeq\fh_\Pi^{(1)}(0|4)}$}\label{psl22}

Recall that the space of \textit{Poisson superalgebra} $\mathfrak{po}_\Pi(0|2n)$ is the Grassmann algebra generated by $\xi=(\xi_1, \dots, \xi_n)$ and $\eta=(\eta_1, \dots, \eta_n)$ with the Poisson bracket
\[
\{f,g\}:=-(-1)^{p(f)}\sum\left(\frac{\partial f}{\partial\xi_i}\frac{\partial g}{\partial\eta_i}+ \frac{\partial f}{\partial\eta_i}\frac{\partial g}{\partial\xi_i}\right).
\]
The \textit{Hamiltonian Lie superalgebra} is $\fh_\Pi(0|2n):=\mathfrak{po}_\Pi(0|2n)/\fc$, where $\fc$ is spanned by constants.

An isomorphism $\fpsl(2|2)\simeq\fh_\Pi^{(1)}(0|4)$ is explicitly given on generators by the following correspondences (other elements are obtained by bracketing): take all basis monomials, except for 1 and $\xi_1\xi_2 \eta_1\eta_2$:
\begin{equation*}
\begin{array}{@{}l@{\,}c@{\,}ll@{\,}c@{\,}ll@{\,}c@{\,}ll@{\,}c@{\,}ll@{\,}c@{\,}l}
\eta_1&\longleftrightarrow &E_{3,2}, &\xi_1\eta_2 &\longleftrightarrow &E_{2,1}, & \xi_1\xi_2& \longleftrightarrow &E_{4,3}, &\xi_1\xi_2 \eta_2 & \longleftrightarrow & E_{2,3},\\
\xi_2\eta_1 &\longleftrightarrow & E_{1,2},& \eta_1\eta_2 & \longleftrightarrow & E_{3,4}, &
1 &\longleftrightarrow & I, & \xi_1\xi_2 \eta_1\eta_2 & \longleftrightarrow & E_{2,2}.
\end{array}
\end{equation*}

\begin{Lemma}\label{gl4_22c} In this Lemma we use convention~\eqref{conv} and $\textbf{F}$ denotes the \textit{desuperization} functor which forgets the squaring if $p=2$.
Let $p=2$.
\begin{enumerate}\itemsep=0pt
\item[$(Aa)$] The space $\text{H}^1\big(\fh_\Pi^{(1)}(0|4); \fh_\Pi^{(1)}(0|4)\big)$ is spanned by the $7$ cocycles $($in formula~\eqref{cocy} the superscript denotes the weight and the subscript the degree$)$:
\begin{gather}
D_{-1}= \xi _1\otimes\big(\widehat{\xi _1 \xi _2 \eta _2}\big)+\xi
 _2\otimes\big(\widehat{\xi _1 \xi _2 \eta _1}\big)+\eta _1\otimes\big(\widehat{\xi _2
 \eta _1 \eta _2}\big)+\eta _2\otimes\big(\widehat{\xi _1 \eta _1 \eta _2}\big), \nonumber\\
D_0^{-2,0}= \eta _1\otimes\big(\widehat{\xi _1}\big)+ \xi _2 \eta_1\otimes \big(\widehat{\xi _1\xi _2}\big)+\eta _1
 \eta _2 \otimes\big(\widehat{\xi _1 \eta _2}\big)+ \xi_2\eta _1 \eta _2\otimes \big(\widehat{\xi _1\xi _2 \eta _2}\big),\nonumber\\
D_0^{0, -2}= \eta _2\otimes\big(\widehat{\xi _2}\big)+ \xi _1 \eta_2\otimes \big(\widehat{\xi _1 \xi _2}\big)+ \eta _1
 \eta _2 \otimes\big(\widehat{\xi _2 \eta _1}\big)+ \xi_1 \eta _1 \eta _2 \otimes\big(\widehat{\xi _1 \xi _2 \eta _1}\big),\nonumber\\
D_0^{0,0}= \xi _2\otimes\big(\widehat{\xi _2}\big)+ \eta_1\otimes\big(\widehat{\eta _1}\big)+ \xi _1 \eta_2\otimes \big(\widehat{\xi _1 \eta _2}\big)+ \xi _2
 \eta _1\otimes \big(\widehat{\xi _2 \eta _1}\big)+ \xi_1 \xi _2 \eta _2 \otimes\big(\widehat{\xi _1 \xi _2 \eta _2}\big)\nonumber\\
 \hphantom{D_0^{0,0}=}{}
 + \xi _1 \eta _1 \eta _2 \otimes\big(\widehat{\xi_1 \eta_1 \eta _2}\big).\label{cocy}
\end{gather}

\item[$(Ab)$] The space $\text{H}^2\big(\fh_\Pi^{(1)}(0|4)\big)$ is spanned by the $15$ cocycles~\eqref{cocyder} and \eqref{cocydery}
\begin{gather}
 c_{0,-2}^1= \big(\widehat{\xi _2}\big)\wedge \big(\widehat{\xi _1 \xi _2 \eta _1}\big)+ \big(\widehat{\xi _1 \xi _2}\big)\wedge \big(\widehat{\xi _2 \eta_1}\big),\nonumber\\
c_{-2,0}^1= \big(\widehat{\xi _1}\big)\wedge \big(\widehat{\xi _1 \xi _2 \eta _2}\big)+ \big(\widehat{\xi _1 \xi _2}\big)\wedge \big(\widehat{\xi _1 \eta_2}\big),\nonumber\\
c_{0,0}^1= \big(\widehat{\xi_1}\big)\wedge \big(\widehat{\eta _1}\big)+\big(\widehat{\xi_2}\big) \wedge \big(\widehat{\eta _2}\big),\nonumber\\
c_{0,0}^2= \big(\widehat{\xi _1}\big)\wedge \big(\widehat{\xi _2 \eta _1 \eta _2}\big)+\big(\widehat{\xi _2}\big)\wedge \big(\widehat{\xi _1 \eta _1 \eta_2}\big)+\big(\widehat{\xi _1 \xi _2}\big)\wedge \big(\widehat{\eta _1
\eta _2}\big),\nonumber\\
c_{0,0}^3= \big(\widehat{\xi _1 \xi _2 \eta _1}\big)\wedge \big(\widehat{\xi _1 \eta _1 \eta
_2}\big)+ \big(\widehat{\xi _1 \xi _2 \eta _2}\big)\wedge
\big(\widehat{\xi _2 \eta _1 \eta _2}\big),\label{cocyder}
\\
c_{-2,0}^2= \big(\widehat{\xi _1}\big)^2, \qquad
c_{-2,0}^3= \big(\widehat{\xi _1 \xi _2 \eta _2}\big)^2,\nonumber\\
c_{0,-2}^2= \big(\widehat{\xi_2}\big)^2,\qquad
c_{0,-2}^3= \big(\widehat{\xi _1 \xi _2 \eta _1}\big)^2.\label{cocydery}
\end{gather}

\item[$(B)$] The non-isomorphic double extensions corresponding to the above outer derivations and central extensions are the three Lie superalgebras described in {\rm \cite{BeBou1}}: $\fgl(2|2)$, $\mathfrak{po}(0|4)$, and one more, denoted $\widetilde{\mathfrak{po}}(0|4)$; for further clarifications, see~{\rm \cite{BLS}}.

\item[$(Ca)$] Let $\fg=\fpsl(4)\simeq \textbf{F}\big(\fh_\Pi^{(1)}(0|4)\big)$. For a~basis of $H^1(\fg;\fg)$ we can
take the following derivations and the corresponding subalgebra of
$\fder\ \fg$ generated by adding just one outer derivation to the
inner ones is indicated on the right if we can identify it $($recall that
$\fa\ltimes\fb$ denotes a~semi-direct sum in which $\fa$ is an ideal$)$:
\[
\begin{array}{@{}l|l}
c_{-4}= y_2\otimes\widehat x_6+y_4\otimes \widehat x_5+y_5\otimes \widehat x_4+y_6
\otimes \widehat
x_2&\\[2mm]
c_{-2}= x_2\otimes \widehat x_6+y_1\otimes \widehat x_3+y_3\otimes \widehat
x_1+y_6\otimes \widehat
y_2&\textbf{F}(\fh_\Pi(0|4))\\
c_{0}^1=x_3\otimes \widehat x_1 +x_5\otimes \widehat x_4
 +y_1\otimes \widehat y_3 +y_4\otimes \widehat y_5&\\
c_{0}^2=x_1\otimes \widehat x_3 +x_4\otimes \widehat x_5
 +y_3\otimes \widehat y_1 +y_5\otimes \widehat y_4&\\
c_{0}^3=x_2\otimes \widehat x_2 +x_3\otimes \widehat x_3
 +x_4\otimes \widehat x_4 +y_2\otimes \widehat y_2&\\
\hphantom{c_{0}^3=}{} +y_3\otimes \widehat y_3
 +y_4\otimes \widehat y_4&\mathfrak{pgl}(4)\simeq \textbf{F}\big(\fh_\Pi^{(1)}(0|4)\ltimes \Kee E\big)
 \end{array}
 \]

\item[$(Cb)$] The nontrivial central extensions of $\fpsl(4)$: the right-most
column contains the description of the result when we can identify
the algebra; the independent cocycles are as follows:
\begin{equation}\label{centextrpsl4}
\begin{array}{@{}l|l}
c_{-4}= \widehat x_2 \wedge \widehat x_6 + \widehat x_4 \wedge \widehat x_5& \\
c_{-2}= \widehat x_1 \wedge \widehat x_3 + \widehat x_6\wedge\widehat y_2&\textbf{F}\big(\mathfrak{po}^{(1)}(0|4)\big)=\mathfrak{po}^{(1)}(4; \One)\\
c_{0}^1= \widehat x_1 \wedge \widehat y_3 + \widehat x_4 \wedge \widehat y_5& \\
c_{0}^2= \widehat x_3 \wedge \widehat y_1 + \widehat x_5 \wedge \widehat y_4& \\
c_{0}^3= \widehat x_2 \wedge \widehat y_2 + \widehat x_3 \wedge \widehat y_3
 + \widehat x_4\wedge \widehat y_4&\fsl(4)
\end{array}
\end{equation}
\end{enumerate}
\end{Lemma}

\begin{Lemma}\label{psl22a} For any $p\neq 2$ $($for $p=2$, see Section~$\ref{sssUnclear})$, we have the following $($recall convention~\eqref{conv}$)$.
\begin{enumerate}\itemsep=0pt
\item[$(a1)$] Let $\fg=\fpsl(2|2)$. For a~basis of $H^1(\fg;\fg)$ we can take
the following derivations \emph{(recall convention~\eqref{conv})}
\begin{gather*}
c_{-4} = -y_2 \otimes \widehat x_6 +
 y_4 \otimes \widehat x_5+ y_5 \otimes \widehat x_4 - y_6 \otimes \widehat x_2,\\
c_{0} = -x_2 \otimes \widehat x_2- x_3 \otimes \widehat x_3 - x_4\otimes \widehat x_4 +
 y_2\otimes \widehat y_2+
 y_3\otimes \widehat y_3+
 y_4 \otimes \widehat y_4
 \end{gather*}
computed for the Cartan matrix
\[
\begin{pmatrix}
 \hphantom{-}2 & -1 & 0 \\
 -1 & \hphantom{-}0 & 1 \\
 \hphantom{-}0 & -1 & 2\end{pmatrix}
 \]
and the basis
\[
 x_1,\ x_2,\ x_3,\
 x_4=[x_1, x_2],\ x_5=[x_2, x_3],\ x_6=[x_3,[x_1, x_2]].
\]

So $\fder \fg\simeq\fosp_{-1}(4|2)$, see formula~\eqref{exSeqOsp4,2,a}.
\item[$(a2)$] For a~basis of $H^2(\fpsl(2|2))$ we can take
\[
\widehat x_2 \wedge \widehat y_2 +\widehat x_3 \wedge \widehat y_3 - \widehat x_4 \wedge \widehat y_4, \qquad \widehat y_2\wedge \widehat y_6+\widehat y_4\wedge \widehat y_5, \qquad \widehat x_2\wedge \widehat x_6+\widehat x_4\wedge \widehat x_5;
\]
the central extension is, clearly, $\fosp_0(4|2)$, see formula~\eqref{exSeqOsp4,2,a}.
\end{enumerate}
\end{Lemma}

\begin{Lemma}\label{2.3.2b} Let $p=3$.
\begin{enumerate}\itemsep=0pt
\item[$(a1)$] For $\fg=\fg(2,3)$ with Cartan matrix
\[
\begin{pmatrix}
 \hphantom{-}2 & -1 & -1 \\
 -1 & \hphantom{-}2 & -1 \\
 -1 & -1 & \hphantom{-}0\end{pmatrix}
 \]
 and basis
 \begin{gather*}
x_1,\ x_2,\ x_3,\
x_4=[x_1, x_2],\
x_5=[x_1, x_3],\ x_6=[x_2, x_3],\
x_7=[x_3, [x_1, x_2]],\\
x_8=[[x_1, x_2], [x_1, x_3]],\ x_9=[[x_1,
x_2], [x_2, x_3]],\
x_{10}=[[x_1, x_2], [x_3, [x_1, x_2]]],\\
x_{11}=[[x_3, [x_1, x_2]], [x_3, [x_1, x_2]]]
\end{gather*}
for a~basis of $H^1(\fg;\fg)$ we can take
the following derivation, \emph{where $D_4$ is an outer derivation, see Remark~$\ref{RemD}$}:
\[
c_{0}=- h_1\otimes \widehat D_4 + h_2\otimes \widehat D_4.
\]

\item[$(a2)$] We have $H^2(\fg)=0$.

\item[$(b1)$] For $\fg=\fbj:=\fg(2,3)^{(1)}/\fc$, for a~basis of
$H^1(\fg;\fg)$ we can take the following derivations $($recall convention~\eqref{conv}$)$
\begin{alignat*}{3}
& \deg=-3\colon \quad &&c_{-3}^1= 2
 x_3\otimes \widehat x_8+x_6\otimes \widehat x_{10}+y_1\otimes \widehat x_4+
 y_4\otimes \widehat x_1+2 y_8\otimes \widehat y_3+y_{10}\otimes \widehat y_6,&\\
&&& c_{-3}^2=x_3\otimes \widehat x_9+x_5\otimes \widehat x_{10}+2 y_2\otimes \widehat x_4+2
 y_4\otimes \widehat x_2+y_9\otimes \widehat y_3+y_{10}\otimes \widehat y_5;&\\
& \deg=0\colon \quad &&c_{0}^1=x_2\otimes \widehat x_1+2 x_6\otimes \widehat x_5+2 x_9\otimes d
x_8+2 y_1\otimes \widehat y_2+y_5\otimes \widehat y_6+y_8\otimes \widehat y_9,&\\
&&&c_{0}^2=x_1\otimes \widehat x_1+2 x_2\otimes \widehat x_2+2 x_5\otimes d
x_5+x_6\otimes \widehat x_6+2
x_8\otimes \widehat x_8+x_9\otimes \widehat x_9& \\
&&&\hphantom{c_{0}^2=}{}
 +2 y_1\otimes \widehat y_1+y_2\otimes \widehat y_2+y_5\otimes \widehat y_5+2 y_6\otimes
\widehat y_6+y_8\otimes \widehat y_8+2
y_9\otimes \widehat y_9,&\\
&&& c_{0}^3=x_1\otimes \widehat x_2+2 x_5\otimes \widehat x_6+2 x_8\otimes \widehat x_9+2
y_2\otimes \widehat y_1+y_6\otimes \widehat y_5+y_9\otimes \widehat y_8.
 \end{alignat*}

\item[$(b2)$] We have $\dim H^2(\fg)=7$; for a~basis we can take the following cocycles $($recall convention~\eqref{conv}$)$
\begin{alignat*}{3}
& \deg=-3\colon \quad && c_{-3}^1=2 \widehat x_1 \wedge \widehat x_4+ \widehat x_8 \wedge \widehat y_3+ \widehat x_{10} \wedge\widehat y_6,&\\
&&&c_{-3}^2=\widehat x_2 \wedge \widehat x_4+2 \widehat x_9 \wedge \widehat y_3+\widehat x_{10} \wedge \widehat y_5, &\\
& \deg=0\colon \quad &&c_{0}^1=2 \widehat x_1 \wedge \widehat y_2 + \widehat x_5 \wedge \widehat y_6 + \widehat x_8 \wedge \widehat y_9,& \\
&&& c_{0}^2=2 \widehat x_1 \wedge \widehat y_1+ x_2 \wedge \widehat y_2+ \widehat x_5 \wedge \widehat y_5 + 2 \widehat x_6 \wedge \widehat y_6+ \widehat x_8 \wedge \widehat y_8+ 2 \widehat x_9 \wedge \widehat y_9,& \\
&&& c_{0}^3=\widehat x_2 \wedge \widehat y_1+2 \widehat x_6 \wedge \widehat y_5 +
2 \widehat x_9 \wedge \widehat y_8.&
 \end{alignat*}
\end{enumerate}
\end{Lemma}

\begin{Lemma}\label{2.3.2c} \quad
\begin{enumerate}\itemsep=0pt
\item[$(a)$] For $p\geq 5$ and $\fg=\mathfrak{ag}(2)$, we have $H^1(\fg;\fg)=0$ and $H^2(\fg)=0$.
\item[$(b)$] For $p=3$ and $\fg=\mathfrak{br}(3)$, we have $H^1(\fg;\fg)=0$ and $H^2(\fg)=0$.
\item[$(c)$] For $p\geq3$ and $\fg=\fg(1,6)$, $\fosp(1|6)$, $\fosp(2|4)$, $\fosp(3|4)$, $\fosp(4|2)$, $\fosp(5|2)$, $\fsl(1|3)$, we have $H^1(\fg;\fg)=0$ and $H^2(\fg)=0$.
\end{enumerate}
\end{Lemma}

\subsection[Rank 4: $\fsl(5)$ and $\fgl(5)$ for $p=5$;
 $\fab(3)$ for $p>3$;
$\fg(3,6)$, $\fg(4,3)$, $\fg(3,3)$ and $\fg^{(1)}(3,3)/\fc$ for $p=3$;
$\fwk(4;a)$ and $\fbgl(4;a)$ for $p=2$]{Rank 4: $\boldsymbol{\fsl(5)}$ and $\boldsymbol{\fgl(5)}$ for $\boldsymbol{p=5}$;
 $\boldsymbol{\fab(3)}$ for $\boldsymbol{p>3}$;
$\boldsymbol{\fg(3,6)}$, $\boldsymbol{\fg(4,3)}$, $\boldsymbol{\fg(3,3)}$\\ and $\boldsymbol{\fg^{(1)}(3,3)/\fc}$ for $\boldsymbol{p=3}$;
$\boldsymbol{\fwk(4;a)}$ and $\boldsymbol{\fbgl(4;a)}$ for $\boldsymbol{p=2}$}

\begin{Lemma}\label{2.4.2a}\quad
\begin{enumerate}\itemsep=0pt
\item[$(a)$] For $p=2$ and $\fg= \fsl(1|4)$,
$\fsl(3|2)$, $\fwk(4;a)$, we have $H^1(\fg;\fg)=0$ and $H^2(\fg)=0$.
\item[$(b)$] For $p=3$ and $\fg=
\fg(4,3)$, $\fg(3,6)$, we have $H^1(\fg;\fg)=0$ and $H^2(\fg)=0$.
\item[$(c)$] For $p=3$ and $\fg=\fosp(1|8)$, $\fosp(2|6)$, $\fosp(3|6)$, $\fosp(4|4)$, $\fosp(5|2)$, we have $H^1(\fg;\fg)=0$ and $H^2(\fg)=0$.
\end{enumerate}
For $p=3$ and $\fg=\fsl(1|4)$, we see that $H^1(\fg;\fg)$ is spanned by
\[
x_2\otimes \widehat x_2 -x_3\otimes \widehat x_3 -x_5\otimes \widehat x_5 -
 x_7\otimes \widehat x_7 -y_2\otimes \widehat y_2 + y_3\otimes \widehat y_3 +y_5\otimes \widehat y_5+ y_7\otimes \widehat y_7,
 \]
where
\[
\begin{pmatrix}
 \hphantom{-}0 & -1 & \hphantom{-}0 &\hphantom{-}0 \\
 -1 & \hphantom{-}2 & -1 &\hphantom{-}0\\
 \hphantom{-}0 & -1 & \hphantom{-}2&-1\\
 \hphantom{-}0& \hphantom{-}0& -1 &\hphantom{-}2
 \end{pmatrix}
 \]
and the basis
\begin{gather*}
 x_1,\ x_2,\ x_3, \ x_4,\
 x_5=[x_1, x_2],\ x_6=[x_2, x_3],\ x_7=[x_3, x_4],\\
 x_8=[x_3,[x_1, x_2]],\ x_9=[x_4,[x_2, x_3]],\ x_{10}=[[x_1,x_2],[x_3, x_4]].
 \end{gather*}

For $\fg=\fsl(2|3)$, we see that $H^2(\fg)=0$ and $H^1(\fg;\fg)$ is spanned by
\[
 - x_2\otimes\widehat x_2 - x_3\otimes\widehat x_3 + x_5\otimes\widehat x_5
 - x_7\otimes\widehat x_7 + y_2\otimes \widehat y_2 + y_3\otimes\widehat y_3
 - y_5\otimes\widehat y_5 + y_7\otimes\widehat y_7,
\]
where the following Cartan matrix was used for computations
\[
\begin{pmatrix}
 \hphantom{-}0 & \hphantom{-}1 & \hphantom{-}0 & \hphantom{-}0 \\
 -1 & \hphantom{-}0 & \hphantom{-}1 & \hphantom{-}0 \\
 \hphantom{-}0 & -1 & \hphantom{-}2 &-1 \\
 \hphantom{-}0 & \hphantom{-}0 & -1 & \hphantom{-}2
 \end{pmatrix}
\]
and the basis
\begin{gather*}
 x_1,\ x_2,\ x_3,\ x_4,\
 x_5=[x_1,x_2],\ x_6=[x_2,x_3],\ x_7=[x_3,x_4],\\
 x_8=[x_3,[x_1,x_2]],\ x_9=[x_4,[x_2,x_3]],\
 x_{10}=[[x_1,x_2],[x_3,x_4]].
\end{gather*}
\end{Lemma}

\begin{Lemma}\label{2.4.2b} Let $p=3$.
\begin{enumerate}\itemsep=0pt
\item[$(a1)$] Let $\fg=\fg(3,3)$ with the
following Cartan matrix and basis
\begin{gather*} 
\begin{pmatrix}
 \hphantom{-}2 & -1 & \hphantom{-}0 & \hphantom{-}0 \\
 -1 & \hphantom{-}2 & -1 & \hphantom{-}0 \\
 \hphantom{-}0 & -2 & \hphantom{-}2 & -1 \\
 \hphantom{-}0 & \hphantom{-}0 & \hphantom{-}1 & \hphantom{-}0\end{pmatrix},
\\
 x_1,\ x_2,\ x_3,\ x_4,\
x_5=[x_1, x_2],\ x_6=[x_2, x_3],\ x_7=[x_3, x_4],\\
x_8=[x_3, [x_1, x_2]],\ x_9=[x_3, [x_2, x_3]],\ x_{10}=[x_4, [x_2, x_3]],\\
x_{11}=[x_3, [x_3, [x_1, x_2]]], \ x_{12}=[[x_1, x_2], [x_3, x_4]], \
x_{13}=[[x_2,
x_3], [x_3, x_4]], \\
x_{14}=[[x_2, x_3], [x_3, [x_1, x_2]]], \ x_{15}=[[x_3, x_4], [x_3,
[x_1, x_2]]], \\
x_{16}=[[x_3, [x_1, x_2]], [x_4, [x_2, x_3]]], \
x_{17}=[[x_4, [x_2, x_3]], [x_3, [x_3, [x_1, x_2]]]].
\end{gather*}
For a~basis of $H^1(\fg;\fg)$ we can take
the following derivation, where $D_5$ is an outer derivation, see Remark $\ref{RemD}$:
\[
c_0=- h_1\otimes \widehat D_5 +h_2\otimes \widehat D_5 +h_4\otimes \widehat D_5.
\]

\item[$(a2)$] We have $H^2(\fg)=0$.

\item[$(b1)$] Let $\fg=\fg(3,3)^{(1)}/\fc$. For a~basis of
$H^1(\fg;\fg)$ we can take the following derivations $($recall convention~\eqref{conv}$)$
\begin{gather*}
c_{-8}= y_4\otimes \widehat x_{17} + y_7\otimes \widehat x_{16} + 2
y_{10}\otimes \widehat x_{15} +y_{12}\otimes \widehat x_{13} +y_{13}\otimes \widehat
x_{12} + 2 y_{15}\otimes \widehat
x_{10}\\
\hphantom{c_{-8}=}{}
+y_{16}\otimes \widehat x_7 +y_{17}\otimes \widehat x_4, \\
c_0= 2 x_4\otimes \widehat x_4 + 2 x_7\otimes \widehat x_7 + 2 x_{10}\otimes \widehat
x_{10} + 2 x_{12}\otimes \widehat x_{12} + 2
x_{13}\otimes \widehat x_{13} + 2 x_{15}\otimes \widehat x_{15} \\
\hphantom{c_0=}{}
 +2 x_{16}\otimes \widehat x_{16} + 2 x_{17}\otimes \widehat x_{17} + y_4\otimes
\widehat y_4 +y_7\otimes \widehat y_7 +y_{10}\otimes \widehat y_{10}
+y_{12}\otimes \widehat y_{12}\\
\hphantom{c_0=}{} +y_{13}\otimes \widehat y_{13}+y_{15}\otimes \widehat y_{15} +y_{16}\otimes \widehat y_{16} +y_{17}\otimes \widehat
y_{17}.
\end{gather*}

\item[$(b2)$] For a~basis of $H^2(\fg)$ we can take $($the classes of$)$ the following cocycles $($recall convention~\eqref{conv}$)$
\begin{gather*}
c_{-8}=\widehat x_4 \wedge
\widehat x_{17} +\widehat x_7 \wedge \widehat x_{16} +\widehat x_{10} \wedge \widehat x_{15} +\widehat x_{12} \wedge \widehat x_{13}, \\
c_{0}=\widehat x_4 \wedge \widehat y_4+ 2 \widehat x_7\wedge \widehat y_7+
 \widehat x_{10} \wedge \widehat y_{10}+ 2
 \widehat x_{12} \wedge
 \widehat y_{12} + 2
 \widehat x_{13} \wedge
 \widehat y_{13} + \widehat x_{15} \wedge
\widehat y_{15} \\
\hphantom{c_{0}=}{}
+ 2 \widehat x_{16} \wedge \widehat y_{16} +\widehat x_{17} \wedge \widehat y_{17}.
\end{gather*}
\end{enumerate}
\end{Lemma}

\begin{Lemma} Let $p=5$. For $\fg=\fosp(1|8)$, $\fosp(2|6)$, $\fosp(3|6)$, $\fosp(4|4)$, $\fosp(5|2)$, we have
$H^1(\fg;\fg)= H^2(\fg)=0$.
\end{Lemma}

\subsection[Rank 5: $\fsl(6)$ and $\fgl(6)$ for $p=2, 3$; $\fg(8,3)$, $\fg(2,6)$
and $\fg^{(1)}(2,6)/\fc$, $\fel(5;3)$ for $p=3$; $\fel(5;5)$ for $p=5$]{Rank 5: $\boldsymbol{\fsl(6)}$ and $\boldsymbol{\fgl(6)}$ for $\boldsymbol{p=2, 3}$; $\boldsymbol{\fg(8,3)}$, $\boldsymbol{\fg(2,6)}$\\
and $\boldsymbol{\fg^{(1)}(2,6)/\fc}$, $\boldsymbol{\fel(5;3)}$ for $\boldsymbol{p=3}$; $\boldsymbol{\fel(5;5)}$ for $\boldsymbol{p=5}$}

\begin{Lemma}\label{2.5.1b} Let $p=2$, $\fg=\fgl(6)$ or its simple relative, or a~superization.
\begin{enumerate}\itemsep=0pt
\item[$(a)$] For $\fg=\fgl(6)$, as well as $\fg=\fgl(5|1)$,
$\fgl(4|2)$, $\fgl(3|3)$, we have
\[
H^1(\fg;\fg)=\Kee\big[h_1\otimes \widehat D_6 + h_3\otimes \widehat D_6 + h_5\otimes \widehat D_6\big],
\]
where $D_6$ is an outer derivation, see Remark $\ref{RemD}$.

We have $H^2(\fg)=0$ for $\fgl(6)$, whereas for its superizations, for a~basis of $H^2(\fg)$ we can take the following cocycles $($recall convention~\eqref{conv}$)$
\begin{alignat*}{5}
&&&\fgl(5|1)&&\fgl(4|2)&&\fgl(3|3)&\\
&\deg=-10\colon \quad &&\widehat x_{15}\wedge\widehat x_{15}\quad &&\widehat x_{15}\wedge\widehat x_{15}\quad &&\widehat x_{15}\wedge\widehat x_{15},&\\
& \deg=-8\colon \quad && \widehat x_{13}\wedge\widehat x_{13}\quad && \widehat x_{13}\wedge\widehat x_{13}, \ \widehat x_{14}\wedge\widehat x_{14} \quad &&\widehat x_{13}\wedge\widehat x_{13}, \ \widehat x_{14}\wedge\widehat x_{14},&\\
&\deg=-6\colon \quad &&\widehat x_{10}\wedge\widehat x_{10}\quad &&\widehat x_{10}\wedge\widehat x_{10}, \ \widehat x_{11}\wedge\widehat x_{11}\quad &&\widehat x_{10}\wedge\widehat x_{10},\ \widehat x_{11}\wedge\widehat x_{11}, \ \widehat x_{12}\wedge\widehat x_{12},&\\
&\deg=-4\colon \quad && \widehat x_{6}\wedge\widehat x_{6}\quad &&\widehat x_{7}\wedge\widehat x_{7}, \ \widehat x_{6}\wedge\widehat x_{6} \quad &&\widehat x_{7}\wedge\widehat x_{7},\ \widehat x_{8}\wedge\widehat x_{8},&\\
&\deg=-2\colon \quad && \widehat x_{1}\wedge\widehat x_{1}\quad &&\widehat x_{2}\wedge\widehat x_{2}\quad &&\widehat x_{3}\wedge\widehat x_{3},&
\end{alignat*}
for the following Chevalley basis
\begin{gather*}
\begin{pmatrix}
 \star & 1 & 0 & 0 & 0\\
1 & \star & 1 & 0 & 0 \\
0 & 1 & \star & 1 & 0\\
0 & 0 & 1 & \star & 1\\
0 & 0 & 0 & 1 & \star
\end{pmatrix}
\\
 x_1,\ x_2,\ x_3,\ x_4, x_5, \
x_6=[x_1, x_2],\ x_7=[x_2, x_3],\ x_8=[x_3, x_4],\ x_9=[x_4, x_5],\\
x_{10}=[x_3, [x_1, x_2]],\ x_{11}=[x_4, [x_2, x_3]],\ x_{12}=[x_5, [x_3, x_4]],\
x_{13}=[[x_1,x_2], [x_3, x_4]],\\
 x_{14}=[[x_2, x_3], [x_4, x_5]], \
 x_{15}=[[x_4, x_5], [x_3, [x_1, x_2]]].
\end{gather*}

\textbf{Distribution of parities}:
For $\fpsl(5|1)$, we take $\{\od,\ev,\ev,\ev,\ev\}$; for $\fpsl(4|2)$, we take $\{\ev, \od,\ev,\ev,\ev\}$; for $\fpsl(3|3)$, we take $\{\ev,\ev,\od,\ev,\ev\}$.

\item[$(b)$] For $\fg=\fpsl(6)$, as well as for $\fpsl(5|1)$,
$\fpsl(4|2)$, $\fpsl(3|3)$, for a~basis of $H^1(\fg;\fg)$ we can
take the following derivation turning $\fpsl$ into $\mathfrak{pgl}$:
\begin{gather*}
c_{0}= x_3\otimes \widehat x_3+ x_4\otimes \widehat x_4+ x_7\otimes \widehat x_7+
x_9\otimes \widehat x_9+ x_{10}\otimes \widehat x_{10}+ y_3\otimes \widehat y_3+y_4\otimes \widehat y_4\\
\hphantom{c_{0}=}{}
+ y_7\otimes \widehat y_7+ y_9\otimes \widehat y_9+ y_{10}\otimes \widehat y_{10}.
\end{gather*}

The space $H^2(\fpsl(6))$ is spanned by the following cocycle:
\[
\widehat x_3\wedge\widehat y_3+\widehat x_4\wedge\widehat y_4+\widehat x_7\wedge\widehat y_7+\widehat x_9\wedge\widehat y_9+\widehat x_{10}\wedge\widehat y_{10}.
\]
For superizations of $\fpsl(6)$, i.e., for $\fg=\fpsl(5|1)$, $\fpsl(4|2)$, and $\fpsl(3|3)$, the space of $H^2(\fg)$ is spanned by the following cocycles $($recall convention~\eqref{conv}$)$
\begin{alignat*}{3}
& \deg=-10\colon \quad &&\widehat x_{15}\wedge\widehat x_{15},\\
&\deg=-8\colon \quad &&\begin{cases} \widehat x_{13}\wedge\widehat x_{13},&\text{for $\fg=\fpsl(5|1)$},\\
\widehat x_{13}\wedge\widehat x_{13}, \ \widehat x_{14}\wedge\widehat x_{14}&\text{for $\fg=\fpsl(4|2)$ and $\fpsl(3|3)$},
\end{cases}&\\
&\deg=-6\colon \quad &&\begin{cases}\widehat x_{10}\wedge\widehat x_{10},&\text{for $\fg=\fpsl(5|1)$},\\
\widehat x_{10}\wedge\widehat x_{10}, \ \widehat x_{11}\wedge\widehat x_{11}&\text{for $\fg=\fpsl(4|2)$},\\
\widehat x_{10}\wedge\widehat x_{10}, \ \widehat x_{11}\wedge\widehat x_{11}, \ \widehat x_{12}\wedge\widehat x_{12}&\text{for $\fg=\fpsl(3|3)$},
\end{cases}& \\
&\deg=-4\colon \quad &&\begin{cases}\widehat x_{6}\wedge\widehat x_{6} &\text{for $\fg=\fpsl(5|1)$},\\
\widehat x_{7}\wedge\widehat x_{7}, \ \widehat x_{6}\wedge\widehat x_{6} &\text{for $\fg=\fpsl(4|2)$ and $\fpsl(3|3)$},\\
\widehat x_{7}\wedge\widehat x_{7}, \ \widehat x_{8}\wedge\widehat x_{8} &\text{for $\fg=\fpsl(3|3)$},
\end{cases}& \\
& \deg=-2\colon \quad &&\begin{cases}\widehat x_{1}\wedge\widehat x_{1} &\text{for $\fg=\fpsl(5|1)$},\\
\widehat x_{2}\wedge\widehat x_{2}
&\text{for $\fg=\fpsl(4|2)$},\\
\widehat x_{3}\wedge\widehat x_{3}&\text{for $\fg=\fpsl(3|3)$},
\end{cases}&\\
& \deg=0\colon \quad &&\widehat x_3\wedge\widehat y_3+\widehat x_4\wedge\widehat y_4+\widehat x_7\wedge\widehat y_7+\widehat x_9\wedge\widehat y_9+\widehat x_{10}\wedge\widehat y_{10}.&
\end{alignat*}
\end{enumerate}
\end{Lemma}

\begin{Lemma}\label{2.5.2a} For $p=3$, consider $\fgl(6)$, its simple relative, and super versions.

\emph{(a1)} For $\fg=\fgl(6)$, and
$\fgl(3|3)$, for a~basis of $H^1(\fg;\fg)$ we can take the following
derivation:
\[
c_{0}=\begin{cases}2 h_1\otimes \widehat D_6 + h_2\otimes \widehat D_6 + 2 h_4\otimes \widehat D_6 +
h_5\otimes \widehat D_6&\text{for $p=3$},\\
4h_1\otimes \widehat D_6 + 3 h_2\otimes \widehat D_6 + 2 h_3\otimes \widehat D_6 +
2 h_4\otimes \widehat D_6 + h_5\otimes \widehat D_6&\text{for $p=5$},\\\end{cases}
\]
where $D_6$ is an outer derivation, see Remark $\ref{RemD}$.

\emph{(a2)} We have $H^2(\fg)=0$.

\emph{(b1)} For $\fg=\fpsl(6)$, we have $\dim H^1(\fg;\fg)=1$; the
outer derivation turns $\fpsl(6)$ into $\mathfrak{pgl}(6)$.

\emph{(b2)} The space $H^2(\fg)$ is spanned by the following cocycle:
\[
c_0=\widehat x_3\wedge\widehat y_3-\widehat x_4\wedge\widehat y_4-\widehat x_7\wedge\widehat y_7+\widehat x_9\wedge\widehat y_9+\widehat x_{10}\wedge\widehat y_{10}.
\]

\emph{(c1)} For $\fg=\fpsl(3|3)$, we have $\dim H^1(\fg;\fg)=1$; the
outer derivation turns $\fpsl(3|3)$ into $\mathfrak{pgl}(3|3)$.

\emph{(c2)} The space $H^2(\fg)$ is spanned by the following cocycle \emph{(for the parities of Chevalley generators distributed as $(\ev,\ev,\od,\ev,\ev)$)}:
\[
c_0=\widehat x_2\wedge\widehat y_2-\widehat x_1\wedge\widehat y_1-\widehat x_7\wedge\widehat y_7+
\widehat x_{11}\wedge\widehat y_{11}-\widehat x_{14}\wedge\widehat y_{14}.
\]
\end{Lemma}

\begin{Lemma}\label{2.5.2b} Let $p=3$.
\begin{enumerate}\itemsep=0pt
\item[$(a1)$] For $\fg=\fg(2,6)$ with Cartan matrix and basis
\begin{gather*}
\begin{pmatrix}
 \hphantom{-}2 & -1 & \hphantom{-}0 & \hphantom{-}0 & \hphantom{-}0 \\
 -1 & \hphantom{-}2 & -1 & \hphantom{-}0 & \hphantom{-}0 \\
 \hphantom{-}0 & -1 & \hphantom{-}2 & -1 & -1 \\
 \hphantom{-}0 & \hphantom{-}0 & -1 & \hphantom{-}2 & \hphantom{-}0 \\
 \hphantom{-}0 & \hphantom{-}0 & \hphantom{-}1 & \hphantom{-}0 & \hphantom{-}0
\end{pmatrix}
\\
x_1,\ x_2,\ x_3,\ x_4,\ x_5,\
x_6=[x_1, x_2],\ x_7=[x_2, x_3],\ x_8=[x_3, x_4],\ x_9=[x_3,
x_5],\\
x_{10}=[x_3, [x_1, x_2]],\ x_{11}=[x_4, [x_2, x_3]],\ x_{12}=[x_5,
[x_2, x_3]],\ x_{13}=[x_5, [x_3, x_4]],\\
x_{14}=[x_5, [x_4, [x_2, x_3]]],\ x_{15}=[[x_1, x_2], [x_3, x_4]],\
x_{16}=[[x_1, x_2], [x_3, x_5]], \\
 x_{17}=[[x_1, x_2], [x_5, [x_3,
x_4]]], \ x_{18}=[[x_3, x_5], [x_4, [x_2, x_3]]],\\
x_{19}=[[x_3, [x_1, x_2]], [x_5, [x_3, x_4]]], \ x_{20}=[[x_5, [x_2,
x_3]], [x_5, [x_3, x_4]]], \\
x_{21}= [[x_5, [x_2, x_3]], [[x_1, x_2], [x_3, x_4]]],\
[x_{22}=[[x_5, [x_3, x_4]], [[x_1, x_2], [x_3, x_5]]], \\
x_{23}= [[x_5, [x_4, [x_2, x_3]]], [[x_1, x_2], [x_3, x_5]]],\\
x_{24}=[[[x_1, x_2], [x_3, x_5]], [[x_3, x_5], [x_4, [x_2,
x_3]]]],\\
x_{25}=[[[x_1, x_2], [x_5, [x_3, x_4]]], [[x_3, x_5], [x_4, [x_2,
x_3]]]]
\end{gather*}
and grading operator $D_6=(1,0,0,0,0)$, see Statement~$\ref{Statement}$, we have
\[
H^1(\fg;\fg)=\Kee\big[2h_1\otimes \widehat D_6+ h_2\otimes \widehat D_6+ h_5\otimes \widehat D_6\big].
\]

\item[$(a2)$] We have $H^2(\fg)=0$.

\item[$(b1)$] Let $\fg=\fg(2,6)^{(1)}/\fc$. For a~basis of $H^1(\fg;\fg)$ we can take
the following derivation
\begin{gather*}
c_0= 2 x_4\otimes \widehat x_4 +2 x_5\otimes \widehat
x_5 +2 x_8\otimes \widehat x_8 +2 x_9\otimes \widehat x_9 +2 x_{11}\otimes \widehat
x_{11} +2 x_{12}\otimes \widehat x_{12}\\
\hphantom{c_0=}{}
+2x_{15}\otimes \widehat x_{15}+2 x_{16}\otimes \widehat x_{16} +x_{20}\otimes \widehat x_{20} +x_{22}\otimes \widehat
x_{22} +x_{23}\otimes \widehat x_{23} +x_{24}\otimes \widehat x_{24}\\
\hphantom{c_0=}{}
+y_4\otimes \widehat
y_4 +y_5\otimes \widehat y_5 +y_8\otimes \widehat y_8 +y_9\otimes \widehat y_9 +y_{11}\otimes \widehat y_{11} +y_{12}\otimes \widehat y_{12}
+y_{15}\otimes \widehat y_{15}\\
\hphantom{c_0=}{}+y_{16}\otimes \widehat y_{16} +2y_{20}\otimes \widehat
y_{20} +2 y_{22}\otimes \widehat y_{22} +2 y_{23}\otimes \widehat y_{23} +2 y_{24}\otimes \widehat y_{24}.
\end{gather*}

\item[$(b2)$] The space $H^2(\fg)$ is spanned by the following cocycle:
\begin{gather*} c_0= \widehat x_1\wedge\widehat y_1-\widehat x_5\wedge\widehat y_5-\widehat x_6\wedge\widehat y_6+\widehat x_9\wedge\widehat y_9+\widehat x_{10}\wedge\widehat y_{10}-\widehat x_{12}\wedge\widehat y_{12}-\widehat x_{13}\wedge\widehat y_{13}\\
 \hphantom{c_0=}{}
 +\widehat x_{14}\wedge\widehat y_{14}-\widehat x_{15}\wedge\widehat y_{15}-\widehat x_{16}\wedge\widehat y_{16}+\widehat x_{17}\wedge\widehat y_{17}-\widehat x_{18}\wedge\widehat y_{18}-\widehat x_{19}\wedge\widehat y_{19}\\
 \hphantom{c_0=}{}
 +\widehat x_{20}\wedge\widehat y_{20}+\widehat x_{21}\wedge\widehat y_{21}.
\end{gather*}
\end{enumerate}
\end{Lemma}

\begin{Lemma}\label{RigRk5} Let $p\geq 3$. For $\fg=\fosp(1|10)$, $\fosp(2|8)$, $\fosp(6|4)$, $\fosp(3|8)$,
$\fosp(4|6)$, $\fosp(5|6)$, $\fosp(7|4)$, $\fosp(8|2)$,
$\fosp(9|2)$, $\fgl(4|2)$, and $\fgl(1|5)$, we have $H^1(\fg;\fg)=0$ and ${H^2(\fg)=0}$.
 \end{Lemma}

\begin{Lemma}\label{RigRk53}Let $p=3$. For $\fg=\fg(8,3)$ and
$\mathfrak{el}(5;3)$, we have $H^1(\fg;\fg)=0$ and ${H^2(\fg)=0}$.
 \end{Lemma}

\begin{Lemma}\label{RigRk55}
Let $p=5$. For $\fg=\mathfrak{el}(5;5)$, we have $H^1(\fg;\fg)=0$ and $H^2(\fg)=0$.
\end{Lemma}

\subsection[Rank 6: $\fsl(7)$ and $\fgl(7)$ for $p=7$;
$\fg(4,6)$, $\fg(6,6)$, and $\fe^{(1)}(6)/\fc$ for $p=3$; $\fe(6)$ and
$\fe(6,i)$ for $p=2$]{Rank 6: $\boldsymbol{\fsl(7)}$ and $\boldsymbol{\fgl(7)}$ for $\boldsymbol{p=7}$;
$\boldsymbol{\fg(4,6)}$, $\boldsymbol{\fg(6,6)}$,\\ and $\boldsymbol{\fe^{(1)}(6)/\fc}$ for $\boldsymbol{p=3}$; $\boldsymbol{\fe(6)}$ and
$\boldsymbol{\fe(6,i)}$ for $\boldsymbol{p=2}$}

\begin{Lemma}\label{2.6.1a} \underline{Let $p=3$}, $\fg=\fe(6)$ and the simple relative. For the Chevalley basis
\begin{gather*}
x_1, \ x_2, \ x_3, \ x_4, \ x_5, \ x_6,\
x_7=[x_1, x_2], \ x_8=[x_2, x_3], \ x_9= [x_3, x_4], \ x_{10}=[x_3, x_6], \\
 x_{11}=[x_4, x_5], \ x_{12}= [x_3, [x_1, x_2]], \
x_{13}=[x_4, [x_2, x_3]], \ x_{14}= [x_5, [x_3, x_4]],\\
 x_{15}= [x_6, [x_2, x_3]], \ x_{16}= [x_6, [x_3, x_4]], \
x_{17}=[x_6, [x_4, [x_2, x_3]]], \ x_{18}= [[x_1, x_2], [x_3, x_4]], \\
 x_{19}= [[x_1, x_2], [x_3, x_6]], \ x_{20}=[[x_2, x_3], [x_4, x_5]], \ x_{21}=[[x_3, x_6], [x_4, x_5]],\\
x_{22}= [[x_1, x_2], [x_6, [x_3, x_4]]], \ x_{23}= [[x_3, x_6], [x_4, [x_2, x_3]]], \ x_{24}=[[x_4, x_5], [x_3, [x_1, x_2]]],\\
 x_{25}= [[x_4, x_5], [x_6, [x_2, x_3]]],\
 x_{26}=[[x_4, x_5], [[x_1, x_2], [x_3, x_6]]], \\
 x_{27}= [[x_3, [x_1, x_2]], [x_6, [x_3, x_4]]], \ x_{28}= [[x_5, [x_3, x_4]], [x_6, [x_2, x_3]]], \\
 x_{29}=[[x_5, [x_3, x_4]], [[x_1, x_2], [x_3, x_6]]], \ x_{30}=[[x_6, [x_2, x_3]], [[x_1, x_2], [x_3, x_4]]],\\
 x_{31}=[[x_6, [x_3, x_4]], [[x_2, x_3], [x_4, x_5]]], \
 x_{32}= [[[x_1, x_2], [x_3, x_4]], [[x_3, x_6], [x_4, x_5]]], \\
 x_{33}= [[[x_1, x_2], [x_3, x_6]], [[x_2, x_3], [x_4, x_5]]],\
 x_{34}= [[[x_2, x_3], [x_4, x_5]], [[x_1, x_2], [x_6, [x_3, x_4]]]], \\ x_{35}=[[[x_3, x_6], [x_4, [x_2, x_3]]], [[x_4, x_5], [x_3, [x_1, x_2]]]], \\
 x_{36}=[[[x_4, x_5], [x_6, [x_2, x_3]]], [[x_3, [x_1,x_2]], [x_6, [x_3, x_4]]]]
 \end{gather*}
we have
\begin{enumerate}\itemsep=0pt
\item[$(a1)$] For $\fg=\fe(6)$, we
have $H^1(\fg;\fg)$ spanned by
\[
2h_1\otimes \widehat D_7+h_2\otimes \widehat D_7+2h_4\otimes \widehat D_7+h_5\otimes \widehat D_7,
\]
where $D_7$ is the outer derivation, see Remark $\ref{RemD}$.

\item[$(a2)$] We have $H^2(\fg)=0$.

\item[$(b1)$] For $\fg=\fe^{(1)}(6)/\fc$, we have $\dim\big(H^1(\fg;\fg)\big)=1$, and
$\fder\ \fg=\fe(6)/\fc$.

\item[$(b2)$] The space $H^2(\fg)$ is spanned by the following cocycle:
\begin{gather*}
c_0=\widehat x_5\wedge\widehat y_5-\widehat x_1\wedge\widehat y_1+\widehat x_7\wedge\widehat y_7-\widehat x_{11}\wedge\widehat y_{11}-\widehat x_{12}\wedge\widehat y_{12}+\widehat x_{14}\wedge\widehat y_{14}+\widehat x_{18}\wedge\widehat y_{18}\\
\hphantom{c_0=}{}
+\widehat x_{19}\wedge\widehat y_{19}-
\widehat x_{20}\wedge\widehat y_{20}-\widehat x_{21}\wedge\widehat y_{21}-\widehat x_{22}\wedge\widehat y_{22}+\widehat x_{25}\wedge\widehat y_{25}+\widehat x_{27}\wedge\widehat y_{27}\\
\hphantom{c_0=}{}
-\widehat x_{28}\wedge\widehat y_{28}-\widehat x_{30}\wedge\widehat y_{30}+\widehat x_{31}\wedge\widehat y_{31}.
\end{gather*}

\item[$(C1)$] For $p=2$ and $\fg=\fe (6;1)$, $\fe (6,6)$, we
have $H^1(\fg;\fg)=0$.

\item[$(C2)$] For $\fe(6,1)$, for a~basis of $H^2(\fg)$ we can take the following cocycles $($recall convention~\eqref{conv}$)$
\begin{alignat*}{5}
& \deg=-22\colon \quad && \widehat x_{36}\wedge\widehat x_{36},\qquad && \deg=-10\colon \quad &&\widehat x_{22}\wedge\widehat x_{22},\ \widehat x_{24}\wedge\widehat x_{24},&\\
& \deg=-20\colon \quad &&\widehat x_{35}\wedge\widehat x_{35},\qquad && \deg=-8\colon \quad &&\widehat x_{18}\wedge\widehat x_{18},\ \widehat x_{19}\wedge\widehat x_{19},&\\
& \deg=-18\colon \quad &&\widehat x_{34}\wedge\widehat x_{34},\qquad && \deg=-6\colon \quad &&\widehat x_{12}\wedge\widehat x_{12},&\\
& \deg=-16\colon \quad &&\widehat x_{33}\wedge\widehat x_{33}, \ \widehat x_{32}\wedge\widehat x_{32},\qquad && \deg=-4\colon \quad &&\widehat x_{7}\wedge\widehat x_{7},&\\
& \deg=-14\colon \quad &&\widehat x_{30}\wedge\widehat x_{30},\ \widehat x_{29}\wedge\widehat x_{29},\qquad &&\deg=-2\colon \quad &&\widehat x_{1}\wedge\widehat x_{1}.&\\
& \deg=-12\colon \quad &&\widehat x_{27}\wedge\widehat x_{27},\ \widehat x_{26}\wedge\widehat x_{26},\qquad && && &
\end{alignat*}

\item[$(C3)$] For $\fe(6,6)$, for a~basis of $H^2(\fg)$ we can take the following cocycles $($recall convention~\eqref{conv}$)$
\begin{alignat*}{3}
&\deg=-20\colon \quad&&\widehat x_{35}\wedge\widehat x_{35},&\\
&\deg=-18\colon \quad&&\widehat x_{34}\wedge\widehat x_{34},&\\
&\deg=-16\colon \quad&&\widehat x_{33}\wedge\widehat x_{33},\
\widehat x_{32}\wedge\widehat x_{32},&\\
&\deg=-14\colon \quad&&\widehat x_{30}\wedge\widehat x_{30},\
\widehat x_{29}\wedge\widehat x_{29},\
\widehat x_{31}\wedge\widehat x_{31},&\\
&\deg=-12\colon \quad&&\widehat x_{27}\wedge\widehat x_{27},\ \widehat x_{28}\wedge\widehat x_{28}, \
\widehat x_{26}\wedge\widehat x_{26},&\\
&\deg=-10\colon \quad&&\widehat x_{23}\wedge\widehat x_{23},\ \widehat x_{22}\wedge\widehat x_{22}, \
\widehat x_{25}\wedge\widehat x_{25},&\\
&\deg=-8\colon \quad&&\widehat x_{19}\wedge\widehat x_{19},\ \widehat x_{17}\wedge\widehat x_{17},\
\widehat x_{21}\wedge\widehat x_{21},&\\
&\deg=-6\colon \quad&&\widehat x_{15}\wedge\widehat x_{15},\ \widehat x_{16}\wedge\widehat x_{16},\\
&\deg=-4\colon \quad&&\widehat x_{10}\wedge\widehat x_{10},&\\
&\deg=-2\colon \quad&&\widehat x_{6}\wedge\widehat x_{6}.&
\end{alignat*}
\end{enumerate}
\end{Lemma}

\begin{Lemma}\label{2.6.2a}
For $p=3$ and $\fg=\mathfrak{g}(4,8)$,
$\mathfrak{g}(6,6)$, we have $H^1(\fg;\fg)=0$ and $H^2(\fg)=0$.
\end{Lemma}

\subsection[Rank 7: $\fsl(8)$ and $\fgl(8)$ for $p=2$; $\fe^{(1)}(7)/\fc$ and
$\fe^{(1)}(7,1)/\fc$, $\fe^{(1)}(7,6)/\fc$, $\fe^{(1)}(7,7)/\fc$ for $p=2$;
$\fg(8,6)$ for $p=3$]{Rank 7: $\boldsymbol{\fsl(8)}$ and $\boldsymbol{\fgl(8)}$ for $\boldsymbol{p=2}$; $\boldsymbol{\fe^{(1)}(7)/\fc}$\\ and
$\boldsymbol{\fe^{(1)}(7,1)/\fc}$, $\boldsymbol{\fe^{(1)}(7,6)/\fc}$, $\boldsymbol{\fe^{(1)}(7,7)/\fc}$ for $\boldsymbol{p=2}$;
$\boldsymbol{\fg(8,6)}$ for $\boldsymbol{p=3}$}

\begin{Lemma}\label{2.7.1a} For $p=2$, let the Chevalley basis
of $\fg=\fe(7,i)$, where $i=1, 6$ or $7$, differ from that of $\fg=\fe(7)$
in that $p(x_i)=p(y_i)=\od$.

\emph{(A)} For $\fg=\fe(7,i)$, where $i=1, 6$ or $7$; we have \emph{(where $D_8$ is the outer derivation, see Remark $\ref{RemD}$)}
\[
H^1(\fg;\fg)=\Kee\big[h_1\otimes \widehat D_8 + h_3\otimes \widehat D_8 + h_7\otimes \widehat D_8\big].
\]

In the rest of Lemma we consider the Chevalley basis
\begin{gather*}
x_1, \ x_2,\ x_3,\ x_4, \ x_5,\ x_6,\ x_7,\
x_8=[x_1, x_2], \ x_9=[x_2, x_3],\ x_{10}=[x_3, x_4], \ x_{11}=[x_4, x_5], \\
 x_{12}= [x_4, x_7],\
 x_{13}= [x_5, x_6], \ x_{14}=[x_3, [x_1, x_2]],\
x_{15}=[x_4, [x_2, x_3]],\ x_{16}=[x_5, [x_3, x_4]], \\
 x_{17}=[x_6, [x_4, x_5]], \
 x_{18}=[x_7, [x_3, x_4]], \ x_{19}=[x_7,[x_4, x_5]], \
 x_{20}=[x_7, [x_5, [x_3, x_4]]], \\
 x_{21}= [[x_1, x_2], [x_3, x_4]], \
 x_{22}=[[x_2, x_3], [x_4, x_5]],\ x_{23}=[[x_2, x_3], [x_4, x_7]], \\
 x_{24}= [[x_3, x_4], [x_5, x_6]],\
 x_{25}=[[x_4, x_7], [x_5, x_6]], \
x_{26}=[[x_2, x_3], [x_7, [x_4, x_5]]],\\
 x_{27}=[[x_4, x_5], [x_3, [x_1, x_2]]], \
 x_{28}=[[x_4, x_7], [x_3, [x_1, x_2]]], \ x_{29}=[[x_4, x_7], [x_5, [x_3, x_4]]],\\
 x_{30}=[[x_5, x_6], [x_4, [x_2, x_3]]], \ x_{31}=[[x_5, x_6], [x_7, [x_3, x_4]]], \\
 x_{32}=[[x_5, x_6], [[x_2, x_3], [x_4, x_7]]], \ x_{33}=[[x_3, [x_1, x_2]], [x_6, [x_4, x_5]]],\\
 x_{34}=[[x_3, [x_1, x_2]], [x_7, [x_4, x_5]]], \
 x_{35}=[[x_4, [x_2, x_3]], [x_7, [x_4, x_5]]], \\
 x_{36}=[[x_6, [x_4, x_5]], [x_7, [x_3, x_4]]], \ x_{37}=[[x_3, [x_1, x_2]], [[x_4, x_7], [x_5, x_6]]], \\
 x_{38}=[[x_6, [x_4, x_5]], [[x_2, x_3], [x_4, x_7]]], \ x_{39}=[[x_7, [x_3, x_4]], [[x_2, x_3], [x_4, x_5]]],\\
 x_{40}=[[x_7, [x_4, x_5]], [[x_1, x_2], [x_3, x_4]]], \ x_{41}= [[x_7, [x_4, x_5]], [[x_3, x_4], [x_5, x_6]]],\\
 x_{42}=[[x_7, [x_5,[x_3, x_4]]], [[x_1, x_2], [x_3\, x_4]]], \ x_{43}= [[[x_1, x_2] [x_3, x_4]], [[x_4, x_7], [x_5, x_6]]],\\
 x_{44}= [[[x_2, x_3], [x_4\, x_5]], [[x_4, x_7], [x_5,x_6]]],\ x_{45}=[[[x_2, x_3], [x_4, x_7]], [[x_3, x_4], [x_5, x_6]]], \\
 x_{46}=[[[x_2, x_3], [x_4, x_7]], [[x_4, x_5], [x_3, [x_1, x_2]]]],\\
 x_{47}= [[[x_3, x_4], [x_5, x_6]], [[x_2, x_3], [x_7, [x_4, x_5]]]], \\
 x_{48}=[[[x_3, x_4], [x_5, x_6]], [[x_4, x_7], [x_3, [x_1, x_2]]]],\\
 x_{49}=[[[x_4, x_7], [x_5, x_6]],[[x_4, x_5],[x_3, [x_1, x_2]]]], \\
 x_{50}=[[[x_4, x_5], [x_3, [x_1, x_2]]], [[x_5, x_6], [x_7, [x_3, x_4]]]],\\
 x_{51}= [[[x_4, x_7], x_3\, [x_1, x_2]]], [[x_5, x_6], [x_4, [x_2, x_3]]]], \\
 x_{52}=[[[x_4, x_7], [x_5, [x_3, x_4]]], [[x_5, x_6], [x_4, [x_2, x_3]]]],\\
 x_{53}=[[[x_4, x_7], [x_5, [x_3, x_4]]], [[x_3, [x_1, x_2]], [x_6, [x_4, x_5]]]], \\
 x_{54}=[[[x_5, x_6], [x_4, [x_2, x_3]]], [[x_3, [x_1, x_2]], [x_7, [x_4, x_5]]]],\\
 x_{55}= [[[x_5, x_6], [x_7, [x_3, x_4]]], [[x_4, [x_2,x_3]], [x_7, [x_4, x_5]]]], \\
 x_{56}=[[[x_3, [x_1, x_2]], [x_6, [x_4, x_5]]], [[x_4, [x_2, x_3]], [x_7, [x_4, x_5]]]],\\
 x_{57}=[[[x_3, [x_1, x_2]], [x_7, [x_4, x_5]]], [[x_6, [x_4, x_5]], [x_7, [x_3, x_4]]]], \\
 x_{58}=[[[x_3, [x_1, x_2]], [x_6, [x_4, x_5]]], [[x_7, [x_3, x_4]], [[x_2, x_3], [x_4, x_5]]]],\\
 x_{59}=[[[x_4, [x_2, x_3]], [x_7, [x_4, x_5]]], [[x_3, [x_1, x_2]], [[x_4, x_7], [x_5, x_6]]]], \\
 x_{60}= [[[x_3, [x_1, x_2]], [[x_4, x_7], [x_5, x_6]]], [[x_7, [x_3, x_4]], [[x_2, x_3], [x_4\, x_5]]]],\\
 x_{61}=[[[x_7, [x_4, x_5]], [[x_1, x_2], [x_3, x_4]]], [[[x_2, x_3], [x_4, x_7]], [[x_3, x_4], [x_5,x_6]]]], \\
 x_{62}=[[[x_7\, [x_5, [x_3, x_4]]], [[x_1, x_2], [x_3, x_4]]], [[[x_2, x_3], [x_4, x_5]], [[x_4, x_7], [x_5, x_6]]]], \\
 x_{63}=[[[[x_2, x_3], [x_4, x_7]], [[x_3, x_4], [x_5, x_6]]], [[[x_4, x_7],[x_5, x_6]], [[x_4, x_5], [x_3, [x_1, x_2]]]]].
 \end{gather*}
\begin{enumerate}\itemsep=0pt
\item[$(A1)$] For $\fg=\fe(7,1)$, for a~basis of $H^2(\fg)$ we can take the following cocycles $($recall convention~\eqref{conv}$)$
\begin{alignat*}{5}
&\deg=-34\colon \quad &&\widehat x_{63}\wedge\widehat x_{63},\quad && \deg=-16\colon \quad &&\widehat x_{42}\wedge\widehat x_{42},\ \widehat x_{43}\wedge\widehat x_{43},&\\
&\deg=-32\colon \quad &&\widehat x_{62}\wedge\widehat x_{62},\quad &&\deg=-14\colon \quad && \widehat x_{40}\wedge\widehat x_{40},\ \widehat x_{37}\wedge\widehat x_{37},&\\
&\deg=-30\colon \quad &&\widehat x_{61}\wedge\widehat x_{61},\quad &&\deg=-12\colon \quad && \widehat x_{34}\wedge\widehat x_{34},\ \widehat x_{33}\wedge\widehat x_{33},&\\
&\deg=-28\colon \quad &&\widehat x_{60}\wedge\widehat x_{60},\qquad &&\deg=-10\colon \quad && \widehat x_{27}\wedge\widehat x_{27},\ \widehat x_{28}\wedge\widehat x_{28},&\\
&\deg=-26\colon \quad &&\widehat x_{58}\wedge\widehat x_{58},\ \widehat x_{59}\wedge\widehat x_{59},\quad &&\deg=-8\colon \quad && \widehat x_{21}\wedge\widehat x_{21},&\\
&\deg=-24\colon \quad && \widehat x_{56}\wedge\widehat x_{56},\ \widehat x_{57}\wedge\widehat x_{57},\quad &&\deg=-6\colon \quad && \widehat x_{14}\wedge\widehat x_{14},&\\
&\deg=-22\colon \quad && \widehat x_{54}\wedge\widehat x_{54},\ \widehat x_{53}\wedge\widehat x_{53},\quad &&\deg=-4\colon \quad && \widehat x_{8}\wedge\widehat x_{8},&\\
&\deg=-20\colon \quad && \widehat x_{51}\wedge\widehat x_{51},\ \widehat x_{50}\wedge\widehat x_{50},\quad &&\deg=-2\colon \quad && \widehat x_{1}\wedge\widehat x_{1}.&\\
&\deg=-18\colon \quad && \widehat x_{46}\wedge\widehat x_{46},\ \widehat x_{48}\wedge\widehat x_{48},\ \widehat x_{49}\wedge\widehat x_{49},\quad && && &
\end{alignat*}

\item[$(A6)$] For $\fg=\fe(7,6)$, for a~basis of $H^2(\fg)$ we can take the following cocycles $($recall convention~\eqref{conv}$)$
\begin{alignat*}{3}
& \deg=-24\colon \quad &&\widehat x_{57}\wedge\widehat x_{57},& \\
& \deg=-22\colon \quad &&\widehat x_{53}\wedge\widehat x_{53},\ \widehat x_{55}\wedge\widehat x_{55},& \\
& \deg=-20\colon \quad &&\widehat x_{50}\wedge\widehat x_{50},\ \widehat x_{52}\wedge\widehat x_{52},&\\
& \deg=-18\colon \quad &&\widehat x_{48}\wedge\widehat x_{48},\ \widehat x_{47}\wedge\widehat x_{47},\ \widehat x_{49}\wedge\widehat x_{49},&\\
& \deg=-16\colon \quad &&\widehat x_{42}\wedge\widehat x_{42},\ \widehat x_{45}\wedge\widehat x_{45},\ \widehat x_{43}\wedge\widehat x_{43},\ \widehat x_{44}\wedge\widehat x_{44},&\\
& \deg=-14\colon \quad &&\widehat x_{39}\wedge\widehat x_{39},\ \widehat x_{40}\wedge\widehat x_{40},\ \widehat x_{37}\wedge\widehat x_{37},\ \widehat x_{38}\wedge\widehat x_{38},&\\
& \deg=-12\colon \quad &&\widehat x_{34}\wedge\widehat x_{34},\ \widehat x_{33}\wedge\widehat x_{33},\ \widehat x_{35}\wedge\widehat x_{35}, \ \widehat x_{32}\wedge\widehat x_{32},& \\
&\deg=-10\colon \quad &&\widehat x_{27}\wedge\widehat x_{27},\ \widehat x_{28}\wedge\widehat x_{28},\ \widehat x_{26}\wedge\widehat x_{26},\ \widehat x_{30}\wedge\widehat x_{30},& \\
& \deg=-8\colon \quad &&\widehat x_{21}\wedge\widehat x_{21},\ \widehat x_{22}\wedge\widehat x_{22},\ \widehat x_{23}\wedge\widehat x_{23},& \\
& \deg=-6\colon \quad &&\widehat x_{14}\wedge\widehat x_{14},& \\
& \deg=-4\colon \quad &&\widehat x_{9}\wedge\widehat x_{9},\ \widehat x_{8}\wedge\widehat x_{8},&\\
& \deg=-2\colon \quad &&\widehat x_{2}\wedge\widehat x_{2}.&
\end{alignat*}

\item[$(A7)$] For $\fg=\fe(7,7)$, for a~basis of $H^2(\fg)$ we can take the following cocycles $($recall convention~\eqref{conv}$)$
\begin{alignat*}{3}
&\deg=-34\colon \quad &&\widehat x_{63}\wedge\widehat x_{63},&\\
&\deg=-32\colon \quad &&\widehat x_{62}\wedge\widehat x_{62},&\\
&\deg=-30\colon \quad &&\widehat x_{61}\wedge\widehat x_{61},&\\
&\deg=-28\colon \quad &&\widehat x_{60}\wedge\widehat x_{60},&\\
&\deg=-26\colon \quad &&\widehat x_{58}\wedge\widehat x_{58},&\\
&\deg=-18\colon \quad &&\widehat x_{49}\wedge\widehat x_{49},&\\
&\deg=-16\colon \quad &&\widehat x_{44}\wedge\widehat x_{44},\ \widehat x_{43}\wedge\widehat x_{43},&\\
& \deg=-14\colon \quad &&\widehat x_{40}\wedge\widehat x_{40},\ \widehat x_{37}\wedge\widehat x_{37}, \ \widehat x_{38}\wedge\widehat x_{38},\ \widehat x_{41}\wedge\widehat x_{41},&\\
&\deg=-12\colon \quad &&\widehat x_{34}\wedge\widehat x_{34},\ \widehat x_{33}\wedge\widehat x_{33},\ \widehat x_{35}\wedge\widehat x_{35},\ \widehat x_{32}\wedge\widehat x_{32},\ \widehat x_{36}\wedge\widehat x_{36},&\\
& \deg=-10\colon \quad && \widehat x_{27}\wedge\widehat x_{27},\ \widehat x_{28}\wedge\widehat x_{28},\ \widehat x_{26}\wedge\widehat x_{26},\ \widehat x_{30}\wedge\widehat x_{30},\ \widehat x_{29}\wedge\widehat x_{29},\ \widehat x_{31}\wedge\widehat x_{31},&\\
& \deg=-8\colon \quad && \widehat x_{21}\wedge\widehat x_{21},\ \widehat x_{22}\wedge\widehat x_{22},\ \widehat x_{23}\wedge\widehat x_{23},\ \widehat x_{20}\wedge\widehat x_{20},\ \widehat x_{24}\wedge\widehat x_{24},&\\
& \deg=-6\colon \quad &&\widehat x_{14}\wedge\widehat x_{14},\ \widehat x_{15}\wedge\widehat x_{15},\ \widehat x_{16}\wedge\widehat x_{16},\ \widehat x_{18}\wedge\widehat x_{18},&\\
&\deg=-4\colon \quad &&\widehat x_{9}\wedge\widehat x_{9},\ \widehat x_{10}\wedge\widehat x_{10},&\\
&\deg=-2\colon \quad && \widehat x_{3}\wedge\widehat x_{3}.&
\end{alignat*}

\item[$(b1)$] For $\fg=\fe^{(1)}(7,i)/\fc $, where $i=1, 6$ or $7$, we
have $\dim H^1(\fg;\fg)=1$ and $\fder\ \fg= \fe(7,i)/\fc$.

\item[$(b2)$] For $\fg=\fe^{(1)}(7,i)/\fc $, where $i=1, 6$ or $7$, we have: $H^2(\fg)$ is
spanned by the same cocycles as for $\fe(7,i)$ and the cocycle
\begin{gather*}
c_0= \widehat x_1 \wedge \widehat y_1+ \widehat x_3 \wedge
 \widehat y_3+ \widehat x_4 \wedge \widehat y_4+
 \widehat x_7 \wedge \widehat y_7+ \widehat x_8\wedge
 \widehat y_8+ \widehat x_9\wedge \widehat y_9+
 \widehat x_{11} \wedge \widehat y_{11} \\
\hphantom{c_0=}{}
+ \widehat x_{17}\wedge
 \widehat y_{17} +\widehat x_{18} \wedge \widehat y_{18} +
 \widehat x_{20} \wedge \widehat y_{20}+ \widehat x_{21}\wedge
 \widehat y_{21}+ \widehat x_{23} \wedge \widehat y_{23}+
 \widehat x_{26} \wedge \widehat y_{26}\\
\hphantom{c_0=}{}
 + \widehat x_{27}\wedge
 \widehat y_{27}+ \widehat x_{31}\wedge \widehat y_{31}+
 \widehat x_{32}\wedge \widehat y_{32}+ \widehat x_{33}\wedge
 \widehat y_{33}+ \widehat x_{39} \wedge \widehat y_{39}+
 \widehat x_{40} \wedge \widehat y_{40}\\
\hphantom{c_0=}{}
 +\widehat x_{43}\wedge
 \widehat y_{43}+ \widehat x_{45}\wedge \widehat y_{45}+
 \widehat x_{47} \wedge \widehat y_{47}+
 \widehat x_{49} \wedge
 \widehat y_{49}+ \widehat x_{53} \wedge \widehat y_{53}+
 \widehat x_{55} \wedge \widehat y_{55}\\
\hphantom{c_0=}{}
 + \widehat x_{56} \wedge
 \widehat y_{56}+ \widehat x_{60} \wedge\widehat y_{60}.
\end{gather*}
\end{enumerate}
\end{Lemma}

\begin{Lemma}\label{2.7.2a} For $p=3$ and $\fg=\mathfrak{g}(8,8)$, we have
$H^1(\fg;\fg)=0$ and $H^2(\fg)=0$.
\end{Lemma}

\subsection[Rank 8: $\fe(8)$ and $\fe(8,1)$, $\fe(8,8)$ for $p=2$]{Rank 8: $\boldsymbol{\fe(8)}$ and $\boldsymbol{\fe(8,1)}$, $\boldsymbol{\fe(8,8)}$ for $\boldsymbol{p=2}$}

\begin{Lemma}\label{2.8.1a} For $p=2$: for $\fg=\fe(8)$, we have $H^1(\fg;\fg)=0$ and $H^2(\fg)=0$.
For $\fg=\fe(8;1)$ and $\fg=\fe(8,8)$, we have $H^1(\fg;\fg)=0$ whereas
\[
H^2(\fg)=\operatorname{Span}(\widehat x_i\wedge \widehat x_i, \ \widehat y_i\wedge \widehat y_i\text{~for all $i$ such
that $x_i$ and $y_i$ are odd}).
\]
\end{Lemma}

\section[Orthogonal series $\fo_\Pi(2n)$ and its super versions in characteristic 2]{Orthogonal series $\boldsymbol{\fo_\Pi(2n)}$ and its super versions\\ in characteristic 2}\label{Sortho}

Here we consider Lie (super)algebras with Cartan matrix, and their relatives.

\begin{Lemma}\label{2.7.2x} Let $H_s$ be the space of cocycles of degree $s$
representing classes in $H^1(\fg;\fg)$. For $\fg=\fo(i,2k)/\fc$, where $k>4$,
 we have
\begin{gather}
\dim H_0=3\qquad \text{for $k$ odd},\nonumber\\
\dim H_0=4\qquad \text{for $k$ even},\nonumber\\
\dim H_{2s}=1\qquad \text{for $-2k+2\leq 2s<0$ and $0<2s\leq 2k-2$,}\nonumber\\
H_{2s+1}=0.\label{dimsH_s}
\end{gather}

Let now $\widetilde H_s$ denote the space of cocycles of degree $s$
representing classes in $H^2(\fg)$. We have the same distribution of dimensions of $\widetilde H_s$ as that of $\dim H_s$ in \eqref{dimsH_s}.
\end{Lemma}

\subsection{Double-checking our results} Permyakov's $\overline{L}$, see \cite{Per}, for the root system $D_n$ is precisely the simple relative of $\fo_\Pi(2n)$, at least for $n>4$. If $p=0$, the elements of $\fo_\Pi(2n)$
can be represented by matrices of the form
\[
 \begin{pmatrix}
 A & B \\
 C & -A^{\rm T}
 \end{pmatrix},
\]
where $A,B,C\in\fgl(n)$ with $B$ and $C$ antisymmetric.

The root elements of the Chevalley basis are
\[
E^{i,j} - E^{n+j, n+i},\qquad E^{i, n+j} - E^{j, n+i},\qquad E^{n+i, j} - E^{n+j, i},
\]
where $1\leq i, j \leq n$ and $i\neq j$. For a~basis of the Cartan subalgebra take
\begin{gather*}
h_i = E^{i,i} - E^{i+1, i+1} - E^{n+i, n+i} + E^{n+i+1, n+i+1}, \qquad \text{where $1\leq i\leq n-1$,}\\
h_n = E^{n-1, n-1} + E^{n, n} - E^{2n-1, 2n-1} - E^{2n, 2n}.
\end{gather*}

These same elements for $p=2$ span $\fo^{(2)}_\Pi(2n)$, but they are not linearly independent: $h_{n-1}=h_n$. But $h_{n-1} + h_n$ is a~central element if $p=2$, so $\fo^{(2)}_\Pi(2n)$ is the quotient of the corresponding Permyakov's $L$ modulo the ideal spanned by this central element.

If $n$ is odd, this ideal is the whole center of $L$, so $\overline{L} = \fo^{(2)}_\Pi(2n)$.

If $n$ is even, then the dimension of the center of $L$ is equal to 2; it is spanned by
$h_{n-1} +h_n$ and $h_1 + h_3 +\dots + h_{n-1}$ (to which the unit matrix corresponds in $\fo^{(2)}_\Pi(2n)$), and hence $\overline{L}$ is the quotient of $\fo^{(2)}_\Pi(2n)$ modulo the 1-dimensional center, i.e., anyway, the simple relative of
$\fo_\Pi(2n)$.

Therefore, describing the algebras of derivations of simple relatives of $\fo_\Pi(2n)$ for
$n>4$, we can refer to~\cite{Per}, same as for $\fsl$ and $\fpsl$.

It is easy to see that the algebra denoted by
$\tilde{\fo}_\Pi$ in~\cite{KrLe} does indeed consist of derivations of simple relatives of
$\fo_\Pi$, and the dimension of its outer derivations is exactly as Permyakov computed, so it contains all of them.

Permyakov's answer $\dim(\fout~\fo^{(2)}_\Pi(8)/\fc)=26$ obtained analytically is also confirmed:
for $\fg=\fo^{(2)}_\Pi(8)/\fc$, \textsc{SuperLie} tells us that $H^1(\fg;\fg)$ is spanned by
\begin{gather*}
c_{-6,1} = y_1\otimes\widehat x_{12}+y_5\otimes\widehat x_{11}+y_8\otimes\widehat x_9+y_9\otimes\widehat x_8+y_{11}\otimes\widehat x_5+y_{12}\otimes\widehat x_1,\\
c_{-6,2} = y_3\otimes\widehat x_{12}+y_6\otimes\widehat x_{11}+y_8\otimes\widehat x_{10}+y_{10}\otimes\widehat x_8+y_{11}\otimes\widehat x_6+y_{12}\otimes\widehat x_3,\\
c_{-6,3} = y_4\otimes\widehat x_{12}+y_7\otimes\widehat x_{11}+y_9\otimes\widehat x_{10}+y_{10}\otimes\widehat x_9+y_{11}\otimes\widehat x_7+y_{12}\otimes\widehat x_4,\\
c_{-4,1} = y_2\otimes\widehat x_{10}+y_6\otimes\widehat x_7+y_7\otimes\widehat x_6+y_{10}\otimes\widehat x_2+x_1\otimes\widehat x_{12}+y_{12}\otimes\widehat y_1,\\
c_{-4,2} = y_2\otimes\widehat x_9+y_5\otimes\widehat x_7+y_7\otimes\widehat x_5+y_9\otimes\widehat x_2+x_3\otimes\widehat x_{12}+y_{12}\otimes\widehat y_3,\\
c_{-4,3} = y_2\otimes\widehat x_8+y_5\otimes\widehat x_6+y_6\otimes\widehat x_5+y_8\otimes\widehat x_2+x_4\otimes\widehat x_{12}+y_{12}\otimes\widehat y_4,\\
c_{-2,1} = y_3\otimes\widehat x_4+y_4\otimes\widehat x_3+x_2\otimes\widehat x_{10}+x_5\otimes\widehat x_{11}+y_{10}\otimes\widehat y_2+y_{11}\otimes\widehat y_5,\\
c_{-2,2} = y_1\otimes\widehat x_4+y_4\otimes\widehat x_1+x_2\otimes\widehat x_9+x_6\otimes\widehat x_{11}+y_9\otimes\widehat y_2+y_{11}\otimes\widehat y_6,\\
c_{-2,3} = y_1\otimes\widehat x_3+y_3\otimes\widehat x_1+x_2\otimes\widehat x_8+x_7\otimes\widehat x_{11}+y_8\otimes\widehat y_2+y_{11}\otimes\widehat y_7,\\
c_{0,1} = x_4\otimes\widehat x_3+x_7\otimes\widehat x_6+x_9\otimes\widehat x_8+y_3\otimes\widehat y_4+y_6\otimes\widehat y_7+y_8\otimes\widehat y_9,\\
c_{0,2} = x_4\otimes\widehat x_1+x_7\otimes\widehat x_5+x_{10}\otimes\widehat x_8+y_1\otimes\widehat y_4+y_5\otimes\widehat y_7+y_8\otimes\widehat y_{10},\\
c_{0,3} = x_3\otimes\widehat x_4+x_6\otimes\widehat x_7+x_8\otimes\widehat x_9+y_4\otimes\widehat y_3+y_7\otimes\widehat y_6+y_9\otimes\widehat y_8,\\
c_{0,4} = x_3\otimes\widehat x_1+x_6\otimes\widehat x_5+x_{10}\otimes\widehat x_9+y_1\otimes\widehat y_3+y_5\otimes\widehat y_6+y_9\otimes\widehat y_{10},\\
c_{0,5} = x_1\otimes\widehat x_4+x_5\otimes\widehat x_7+x_8\otimes\widehat x_{10}+y_4\otimes\widehat y_1+y_7\otimes\widehat y_5+y_{10}\otimes\widehat y_8,\\
c_{0,6} = x_1\otimes\widehat x_3+x_5\otimes\widehat x_6+x_9\otimes\widehat x_{10}+y_3\otimes\widehat y_1+y_6\otimes\widehat y_5+y_{10}\otimes\widehat y_9,\\
c_{0,7} = x_1\otimes\widehat x_1+x_3\otimes\widehat x_3+x_5\otimes\widehat x_5+x_6\otimes\widehat x_6+x_9\otimes\widehat x_9+x_{10}\otimes\widehat x_{10}+y_1\otimes\widehat y_1+y_3\otimes\widehat y_3\\
\hphantom{c_{0,7} =}{}
+y_5\otimes\widehat y_5 +y_6\otimes\widehat y_6+y_9\otimes\widehat y_9+y_{10}\otimes\widehat y_{10},\\
c_{0,8} = x_1\otimes\widehat x_1+x_4\otimes\widehat x_4+x_5\otimes\widehat x_5+x_7\otimes\widehat x_7+x_8\otimes\widehat x_8+x_{10}\otimes\widehat x_{10}+y_1\otimes\widehat y_1+y_4\otimes\widehat y_4\\
\hphantom{c_{0,8} =}{}
+y_5\otimes\widehat y_5 +y_7\otimes\widehat y_7+y_8\otimes\widehat y_8+y_{10}\otimes\widehat y_{10}.
\end{gather*}

Since $\fo^{(2)}_\Pi(6)\simeq\fpsl(4)$, our description of $\fder\, \fo^{(2)}_\Pi(6)$ matches that by Permyakov.
The isomorphism is given by the following correspondence:
\begin{alignat*}{3}
& E^{2,3}\longmapsto E^{1,2} + E^{5,4},\qquad && E^{1,3}\longmapsto E^{1,3} + E^{6,4},& \\
& E^{1,2}\longmapsto E^{2,3} + E^{6,5},\qquad && E^{2,4}\longmapsto E^{1,6} + E^{3,4},&\\
& E^{3,4}\longmapsto E^{2,6} + E^{3,5},\qquad && E^{1,4}\longmapsto E^{1,5} + E^{2,4},&
\end{alignat*}
and the basis elements, obtained from the above ones on the left, go to the transposed images on the right, i.e., $E^{3,2}\longmapsto E^{2,1} + E^{4,5}$ and so on; besides,
\[
E^{1,1} + E^{2,2}\bigl(=E^{3,3} + E^{4,4}\text{~in $\fpsl$}\bigr)\longmapsto E^{2,2} + E^{3,3} + E^{5,5} + E^{6,6}
\]
whereas
\[
E^{2,2} + E^{3,3}\bigl(=E^{1,1} + E^{4,4}\text{~in $\fpsl$}\bigr)\longmapsto E^{1,1} + E^{2,2} + E^{4,4} + E^{5,5}.
\]

\subsection{Rank 2}

\begin{Lemma}\label{2.2.1bO} Let $\fg=\fo\fo_{I \Pi}^{(1)}(1|4)$. Then, $H^1(\fg;\fg)=0$ and $H^2(\fg)=0$.
\end{Lemma}

\begin{Lemma}\label{2.2.1aO} Let $\fg=\fo \fo_{I \Pi}^{(1)}(3|2)$ with Cartan matrix
\[
\begin{pmatrix}
 0& 1 \\
 1 & 1 \end{pmatrix}
\]
and a~basis
 \begin{equation*}
x_1,\ x_2,\ x_3=x_2^2,\ x_4=[x_1, x_2],\ x_5=[x_2, x_2^2],
\end{equation*}
\begin{enumerate}\itemsep=0pt
\item[$(a1)$] For a~basis of $H^1(\fg;\fg)$ we can take the following derivations $($recall convention~\eqref{conv}$)$
\[
c_{-4}=y_2\otimes\widehat x_5+y_4\otimes\widehat x_4 +y_5\otimes\widehat x_2.
\]

\item[$(a2)$] For a~basis of $H^2(\fg)$ we can take (the classes of) the following cocycles $($recall convention~\eqref{conv}$)$
\[
\deg=-6\colon \quad \widehat x_5\wedge\widehat x_5,\qquad \deg=-2\colon \quad \widehat x_2\wedge\widehat x_2.
\]
\end{enumerate}
\end{Lemma}

\subsection{Rank 3}

\subsubsection[The analog of the Lie superalgebra $\fosp_a(4|2)$ for $p=2$]{The analog of the Lie superalgebra $\boldsymbol{\fosp_a(4|2)}$ for $\boldsymbol{p=2}$}\label{sssUnclear}
As we know (e.g., from \cite{BGL1}), over $\Cee$, one of the Cartan
matrices of $\fosp_a(4|2)$ is
\[
\begin{pmatrix} \hphantom{-}2 & -1 & \hphantom{-}0 \\ -1 & \hphantom{-}0 & -a \\ \hphantom{-}0 & -1 & \hphantom{-}2
\end{pmatrix},
\]
 so $\fosp_1(4|2)=\fosp(4|2)$, and $\fosp_a(4|2)$ is simple
except for $a=0$ and $-1$, where we have exact sequences
\begin{gather}
0\tto \fpsl(2|2)\tto \fosp_{-1}(4|2)\tto \fsl(2)\tto 0,\nonumber\\
0\tto \fsl(2)\tto \fosp_0(4|2)\tto \fpsl(2|2)\tto 0.\label{exSeqOsp4,2,a}
\end{gather}
In other words, $\fder\ \fpsl(2|2)\simeq \fosp_{-1}(4|2)$ and on the
space of three nontrivial central extensions of $\fpsl(2|2)$, as well on the space of outer derivations of $\fpsl(2|2)$, a~natural structure of the Lie algebra $\fsl(2)$ can be defined.

\begin{Remark}
We wished to conjecture that over $\Kee$ the central extensions of $\fpsl(4)$---i.e., of the desuperization of $\fpsl(2|2)$---constitute the Lie algebra $\fo(3)\oplus \Kee$ and something
similar is true for $\fpsl(2|2)$:
\[
0\tto \fo(3)\oplus \Kee\tto \text{?} \tto \fpsl(4)\tto 0
\]
but the actual answers \eqref{centextrpsl4} and \eqref{exSeqOsp4,2,a1} are
quite different. 
For $p=2$, there are considerably more of both outer derivations and
central extensions of $\fpsl(2|2)$ than there are outer derivations and central extensions
of $\fpsl(2|2)$ over $\Cee$; the same applies to its desuperization,
$\fpsl(4)$. Thanks to the fact that $\fpsl(2|2)\simeq\fh_\Pi^{(1)}(0|4)$ and
$\fpsl(4)=\fh_\Pi^{(1)}(4;\One)\simeq \textbf{F}(\fh_\Pi^{(1)}(0|4))$, where $\textbf{F}$ is the \textit{desuperization} functor, and recent results \cite{BeBou1}, we were able to understand the meaning of the answers \eqref{centextrpsl4} and \eqref{exSeqOsp4,2,a1}, see
Lemma \ref{gl4_22}.
\end{Remark}

\begin{Lemma}\label{ospForP2}
Over $\Kee$ for $p=2$, the role of $\fosp_a(4|2)$ is played by
$\fbgl(3;a)$, where $a\neq 0, 1$, and its desuperization,
$\fwk(3;a)$. For these algebras we have the following analogs of exact sequences~\eqref{exSeqOsp4,2,a}:
\begin{gather*}
0\tto \Kee^3\tto \fwk(3;1)\tto \fpsl(4)\tto 0,\\
0\tto \Kee^7\tto \fwk(3;0)\tto \Kee^2\oplus \fsl(3)\tto 0,
\end{gather*}
where $\Kee^7$ is a~commutative ideal which, as the $\fwk(3;0)$-module, is the direct sum of three irreducible modules of dimensions $1\oplus 3\oplus 3$,
\begin{gather}
0\tto \Kee^3\tto \fbgl(3;1)\tto \fpsl(2|2)\tto 0,\nonumber\\
0\tto \Kee^{5|2}\tto \fbgl(3;0)\tto \Pi(\Kee^2)\oplus \fsl(1|2)\tto 0,\label{exSeqOsp4,2,a1}
\end{gather}
where $\Kee^{5|2}$ is a~commutative ideal which, as the $\fbgl(3;0)$-module, is the direct sum of three irreducible modules of dimensions $1\oplus 2|1\oplus
2|1$.
\end{Lemma}

\subsection{Rank 4}

\begin{Lemma}\label{2.4.1aO} \quad
\begin{enumerate}\itemsep=0pt
\item[$(a)$] Let $\fg=\fo\fc(1;8)\ltimes \Kee I_0$,
$\fo\fo\fc(1;4|4)\ltimes \Kee I_0$, and $\fo\fo\fc(1;6|2)/\fc$ with
Cartan matrix
\[
\begin{pmatrix}
 \ev & 0 & 0 & 1 \\
 0 & \ev & 0 & 1 \\
 0 & 0 & \ev& 1 \\
 1 & 1 & 1 & \ev \end{pmatrix}
 \]
 and basis
\begin{gather*}
x_1,\ x_2,\ x_3,\ x_4,\ x_5=[x_1, x_4],\ x_6=[x_2, x_4],\ x_7=[x_3, x_4],\
x_8=[x_2, [x_1, x_4]],\nonumber\\
 x_9=[x_3, [x_1, x_4]],\
 x_{10}=[x_3, [x_2, x_4]], \ x_{11}=[x_3, [x_2, [x_1, x_4]]],\nonumber\\
x_{12}=[[x_3, x_4], [x_2, [x_1, x_4]]].
\end{gather*}
 For a~basis of $H^1(\fg;\fg)$ we can take the following derivations:
\begin{gather*}
c_{0}^1=h_1\otimes\widehat h_5 + h_2\otimes\widehat h_5, \qquad
c_{0}^2=h_1\otimes\widehat h_6 + h_2\otimes\widehat h_6, \\
c_{0}^3=h_1\otimes\widehat h_5 + h_3\otimes\widehat h_5, \qquad
c_{0}^4=h_1\otimes\widehat h_6 + h_3\otimes\widehat h_6.
 \end{gather*}

For $\fo\fo\fc(1;8)\ltimes \Kee I_0$, we have $H^2(\fg)=\operatorname{Span}(\widehat h_5\wedge\widehat h_6)$.

Let $\fg=\fo\fo\fc(1;4|4)\ltimes \Kee I_0$, where the parity is distributed as follows: all simple roots are even except $x_4$ and $y_4$. For a~basis of $H^1(\fg;\fg)$ we can take
\begin{alignat*}{3}
&\deg=-8\colon \quad && \widehat x_{11}\wedge\widehat x_{11},&\\
&\deg=-6\colon \quad && c_{-6}^1=\widehat x_{8}\wedge\widehat x_{8}, \ c_{-6}^2=\widehat x_{9}\wedge\widehat x_{9}, \ c_{-6}^3=\widehat x_{10}\wedge\widehat x_{10},&\\
&\deg=-4\colon \quad && c_{-4}^1=\widehat x_{5}\wedge\widehat x_{5}, \ c_{-4}^2=\widehat x_{6}\wedge\widehat x_{6}, \ c_{-4}^3=\widehat x_{7}\wedge\widehat x_{7},&\\
&\deg=-2\colon \quad &&\widehat x_{4}\wedge\widehat x_{4},&\\
&\deg=0\colon \quad &&\widehat h_{5}\wedge\widehat h_{6}.&
\end{alignat*}

Let $\fg=\fo\fo\fc(1;6|2)\ltimes \Kee I_0$, where the parity is distributed as follows: all simple roots are even except $x_1$ and $y_1$. For a~basis of $H^1(\fg;\fg)$ we can take
\begin{alignat*}{5}
& \deg=-12\colon \quad &&\widehat x_{12}\wedge\widehat x_{12},\qquad && \deg=-4\colon \quad &&\widehat x_{5}\wedge\widehat x_{5},& \\
& \deg=-8\colon \quad &&\widehat x_{11}\wedge\widehat x_{11},\qquad && \deg=-2\colon \quad &&\widehat x_{1}\wedge\widehat x_{1},& \\
& \deg=-6\colon \quad && c_{-6}^1=\widehat x_{8}\wedge\widehat x_{8}, \ c_{-6}^2=\widehat x_{9}\wedge\widehat x_{9},\qquad && \deg=0\colon \quad &&\widehat h_{5}\wedge\widehat h_{6}.&
\end{alignat*}

\item[$(b)$] Let $\fg=\fo\fc(1;8)/\fc$, $\fo\fo\fc(1;4|4)/\fc$ and
$\fo\fo\fc(1;6|2)/\fc$. For a~basis of $H^1(\fg;\fg)$ we can take
the following derivations $($recall convention~\eqref{conv}$)$
\begin{gather*}
c_{-6}^1= y_1\otimes \widehat x_{12} +y_5\otimes \widehat x_{11} +y_8\otimes \widehat x_9
+y_9\otimes \widehat x_8 +y_{11}\otimes \widehat x_5
+y_{12}\otimes \widehat x_1,\\
c_{-6}^2= y_2\otimes \widehat x_{12} +y_6\otimes \widehat x_{11} +y_8\otimes \widehat
x_{10} +y_{10}\otimes \widehat x_8 +y_{11}\otimes \widehat x_6
+y_{12}\otimes \widehat x_2,\\
c_{-6}^3= y_3\otimes \widehat x_{12} +y_7\otimes \widehat x_{11} +y_9\otimes\widehat
x_{10} +y_{10}\otimes \widehat x_9 +y_{11}\otimes \widehat x_7
+y_{12}\otimes \widehat x_3,\\
c_{-4}^1= x_1\otimes \widehat x_{12} +y_4\otimes \widehat x_{10} +y_6\otimes \widehat x_7
+y_7\otimes \widehat x_6 +y_{10}\otimes \widehat x_4
+y_{12}\otimes \widehat y_1,\\
c_{-4}^2= x_2\otimes \widehat x_{12} +y_4\otimes \widehat x_9 +y_5\otimes \widehat x_7
+y_7\otimes \widehat x_5 +y_9\otimes \widehat x_4
+y_{12}\otimes \widehat y_2,\\
c_{-4}^3= x_3\otimes \widehat x_{12} +y_4\otimes \widehat x_8 +y_5\otimes \widehat x_6
+y_6\otimes \widehat x_5 +y_8\otimes \widehat x_4
+y_{12}\otimes \widehat y_3,\\
c_{-2}^1= x_4\otimes \widehat x_{10} +x_5\otimes \widehat x_{11} +y_2\otimes \widehat x_3
+y_3\otimes \widehat x_2 +y_{10}\otimes \widehat y_4
+y_{11}\otimes \widehat y_5,\\
c_{-2}^2= x_4\otimes \widehat x_9 +x_6\otimes \widehat x_{11} +y_1\otimes \widehat x_3
+y_3\otimes \widehat x_1 +y_9\otimes \widehat y_4
+y_{11}\otimes \widehat y_6,\\
c_{-2}^3= x_4\otimes \widehat x_8 +x_7\otimes \widehat x_{11} +y_1\otimes \widehat x_2
+y_2\otimes \widehat x_1 +y_8\otimes \widehat y_4
+y_{11}\otimes \widehat y_7,\\
c_{0}^1= x_3\otimes \widehat x_2 +x_7\otimes \widehat x_6 +x_9\otimes \widehat x_8
+y_2\otimes \widehat y_3
+y_6\otimes \widehat y_7 +y_8\otimes \widehat y_9,\\
c_{0}^2= x_3\otimes \widehat x_1 +x_7\otimes \widehat x_5 +x_{10}\otimes \widehat x_8
+y_1\otimes \widehat y_3
+y_5\otimes \widehat y_7 +y_8\otimes \widehat y_{10},\\
c_{0}^3= x_2\otimes \widehat x_3 +x_6\otimes \widehat x_7 +x_8\otimes \widehat x_9
+y_3\otimes \widehat y_2
+y_7\otimes \widehat y_6 +y_9\otimes \text{dy}_{8},\\
c_{0}^4= x_2\otimes \widehat x_1 +x_6\otimes \widehat x_5 +x_{10}\otimes \widehat x_9
+y_1\otimes \widehat y_2
+y_5\otimes \widehat y_6 +y_9\otimes \widehat y_{10},\\
c_{0}^5= x_1\otimes \widehat x_3 +x_5\otimes \widehat x_7 +x_8\otimes \widehat x_{10}
+y_3\otimes \widehat y_1
+y_7\otimes \widehat y_5 +y_{10}\otimes \widehat y_8,\\
c_{0}^6= x_1\otimes \widehat x_2 +x_5\otimes \widehat x_6 +x_9\otimes \widehat x_{10}
+y_2\otimes \widehat y_1
+y_6\otimes \widehat y_5 +y_{10}\otimes \widehat y_9,\\
c_{0}^7= x_1\otimes \widehat x_1 +x_2\otimes \widehat x_2 +x_5\otimes \widehat x_5
+x_6\otimes \widehat x_6 +
x_9\otimes \widehat x_9 +x_{10}\otimes \widehat x_{10}+ y_1\otimes \widehat y_1 \\
\hphantom{c_{0}^7=}{}
+y_2\otimes \widehat y_2 +y_5\otimes \widehat y_5
+y_6\otimes \widehat y_6 +y_9\otimes \widehat y_9 +y_{10}\otimes \widehat y_{10},\\
c_{0}^8= x_1\otimes \widehat x_1+x_3\otimes \widehat x_3+x_5\otimes \widehat
x_5+x_7\otimes \widehat x_7+ x_8\otimes \widehat x_8+x_{10}\otimes \widehat
x_{10}+
 y_1\otimes \widehat y_1\\
\hphantom{c_{0}^8=}{}
 +y_3\otimes \widehat y_3+ y_5\otimes \widehat y_5+y_7\otimes \widehat
y_7+y_8\otimes \widehat y_8+y_{10}\otimes \widehat y_{10}.
 \end{gather*}

\item[$(c)$] For a~basis of $H^2(\fg)$ we can take the following cocycles $($recall convention~\eqref{conv}$)$.

\item[$(c1)$] For $\fg=\fo\fc(1;8)/\fc$:
\begin{alignat*}{3}
& \deg=-6\colon \quad && c_{-6}^1=\widehat x_1 \wedge \widehat x_{12}+ \widehat x_5 \wedge \widehat x_{11}+ \widehat x_8 \wedge \widehat x_9, &\\
&&& c_{-6}^2=x_2 \wedge \widehat x_{12}+ \widehat x_6 \wedge \widehat x_{11}+ \widehat x_8 \wedge \widehat x_{10}, &\\
&&& c_{-6}^3=
 \widehat x_3 \wedge \widehat x_{12}+ \widehat x_7 \wedge \widehat x_{11} + \widehat x_9 \wedge \widehat x_{10},&\\
& \deg=-4\colon \quad && c_{-4}^1=\widehat x_4 \wedge \widehat x_{10} + \widehat x_6 \wedge \widehat x_7+ \widehat x_{12} \wedge \widehat y_1, &\\
&&& c_{-4}^2=\widehat x_4 \wedge \widehat x_9 +\widehat x_5 \wedge \widehat x_7 +\widehat x_{12} \wedge \widehat y_2, &\\
&&& c_{-4}^3= \widehat x_4 \wedge \widehat x_8 +\widehat x_5 \wedge \widehat x_6 +\widehat x_{12} \wedge \widehat y_3, &\\
& \deg=-2\colon \quad && c_{-2}^1= \widehat x_2 \wedge \widehat x_3+\widehat x_{10} \wedge \widehat y_4+\widehat x_{11} \wedge \widehat y_5, &\\
&&& c_{-2}^2=\widehat x_1\wedge \widehat x_3+\widehat x_9 \wedge \widehat y_4 +x_{11}\wedge \widehat y_6,&\\
&&&c_{-2}^3=\widehat x_1\wedge \widehat x_2+\widehat x_8 \wedge \widehat y_4+\widehat x_{11} \wedge \widehat y_7,&\\
& \deg=0\colon \quad &&c_{0}^1= \widehat x_2 \wedge \widehat y_3 +
\widehat x_6 \wedge \widehat y_7+
\widehat x_8 \wedge \widehat y_9, &\\
&&&c_{0}^2= \widehat x_1\wedge \widehat y_3+\widehat x_5 \wedge \widehat y_7 +\widehat x_8 \wedge \widehat y_{10}, &\\
&&&c_{0}^3= \widehat x_3 \wedge \widehat y_2 +\widehat x_7 \wedge \widehat y_6 +\widehat x_9 \wedge\widehat y_8, &\\
&&&c_{0}^4= \widehat x_1 \wedge \widehat y_2+\widehat x_5 \wedge \widehat y_6 +\widehat x_9 \wedge \widehat y_{10},&\\
&&&c_{0}^5= \widehat x_3 \wedge \widehat y_1+
 \widehat x_7 \wedge \widehat y_5+x_{10}\wedge \widehat y_8, &\\
&&&c_{0}^6= \widehat x_2\wedge \widehat y_1+\widehat x_6 \wedge \widehat y_5+\widehat x_{10} \wedge \widehat y_9,&\\
&&&c_{0}^7= \widehat x_2 \wedge \widehat y_2+\widehat x_3 \wedge \widehat y_3+\widehat x_6 \wedge \widehat y_6 +\widehat x_7 \wedge \widehat y_7+x_8 \wedge \widehat y_8+\widehat x_9 \wedge
 \widehat y_9,&\\
 &&&c_{0}^8= \widehat x_1\wedge \widehat y_1+
 \widehat x_3\wedge \widehat y_3+x_5\wedge \widehat y_5+x_7\wedge \widehat y_7+x_8 \wedge \widehat y_8 +\widehat x_{10} \wedge \widehat y_{10}.&
\end{alignat*}

\item[$(c2)$] For $\fg=\fo\fc(1;4|4)/\fc$:
\begin{alignat*}{3}
&\deg=-8\colon \quad &&c_{-8}=\widehat x_{11}\wedge\widehat x_{11},& \\
&\deg=-6\colon \quad &&c_{-6}^1=\widehat x_{8}\wedge\widehat x_{8}, \ \widehat x_{9}\wedge\widehat x_{9}, \ \widehat x_{10}\wedge\widehat x_{10},&\\
&&&c_{-6}^2=\widehat x_1\wedge\widehat x_{12}+\widehat x_5\wedge\widehat x_{11}+\widehat x_8\wedge\widehat x_9,&\\
&&&c_{-6}^3=\widehat x_2\wedge\widehat x_{12}+\widehat x_6\wedge\widehat x_{11}+\widehat x_8\wedge\widehat x_{10},&\\
&&&c_{-6}^4=\widehat x_3\wedge\widehat x_{12}+\widehat x_7\wedge\widehat x_{11}+\widehat x_9\wedge\widehat x_{10},&\\
&\deg=-4\colon \quad &&c_{-4}^1=\widehat x_{5}\wedge\widehat x_{5},\widehat x_{6}\wedge\widehat x_{6},\widehat x_{7}\wedge\widehat x_{7},&\\
&&&c_{-4}^2=\widehat x_4\wedge\widehat x_{10}+\widehat x_6\wedge\widehat x_{7}+\widehat x_{12}\wedge\widehat y_1,&\\
&&&c_{-4}^3=\widehat x_4\wedge\widehat x_{9}+\widehat x_5\wedge\widehat x_{7}+\widehat x_{12}\wedge\widehat y_2,&\\
&&&c_{-4}^4=\widehat x_4\wedge\widehat x_{8}+\widehat x_5\wedge\widehat x_{6}+\widehat x_{12}\wedge\widehat y_3,&\\
&\deg=-2\colon \quad &&c_{-2}^1=\widehat x_{4}\wedge\widehat x_{4},&\\
&&&c_{-2}^2=\widehat x_2\wedge\widehat x_{3}+\widehat x_{10}\wedge\widehat y_{4}+\widehat x_{11}\wedge\widehat y_5,&\\
&&&c_{-2}^3=\widehat x_1\wedge\widehat x_{3}+\widehat x_{9}\wedge\widehat y_{4}+\widehat x_{11}\wedge\widehat y_6,&\\
&&&c_{-2}^4=\widehat x_1\wedge\widehat x_{2}+\widehat x_{8}\wedge\widehat y_{4}+\widehat x_{11}\wedge\widehat y_7,	&\\
& \deg=0\colon \quad &&c_{0}^1=\widehat x_2\wedge\widehat y_{3}+\widehat x_{5}\wedge\widehat y_{7}+\widehat x_{8}\wedge\widehat y_9,&\\
&&&c_{0}^2=\widehat x_1\wedge\widehat y_{3}+\widehat x_{5}\wedge\widehat y_{7}+\widehat x_{8}\wedge\widehat y_{10},&\\
&&&c_{0}^3=\widehat x_3\wedge\widehat y_{2}+\widehat x_{7}\wedge\widehat y_{6}+\widehat x_{9}\wedge\widehat y_{8},&\\
&&&c_{0}^4=\widehat x_1\wedge\widehat y_{2}+\widehat x_{5}\wedge\widehat y_{6}+\widehat x_{9}\wedge\widehat y_{10},&\\
&&&c_{0}^5=\widehat x_3\wedge\widehat y_{1}+\widehat x_{7}\wedge\widehat y_{5}+\widehat x_{10}\wedge\widehat y_{8},&\\
&&&c_{0}^6=\widehat x_2\wedge\widehat y_{1}+\widehat x_{6}\wedge\widehat y_{5}+\widehat x_{10}\wedge\widehat y_{9},&\\
&&&c_{0}^7=\widehat x_2\wedge\widehat y_{2}+\widehat x_{3}\wedge\widehat y_{3}+\widehat x_{6}\wedge\widehat y_{6} +\widehat x_7\wedge\widehat y_{7} +\widehat x_8\wedge\widehat y_{8}+\widehat x_9\wedge\widehat y_{9} ,&\\
&&&c_{0}^8=\widehat x_1\wedge\widehat y_{1}+\widehat x_{3}\wedge\widehat y_{3}+\widehat x_{5}\wedge\widehat y_{5} +\widehat x_7\wedge\widehat y_{7} +\widehat x_8\wedge\widehat y_{8}+\widehat x_{10}\wedge\widehat y_{10}.&
\end{alignat*}

\item[$(c3)$] For $\fg=\fo\fc(1;6|2)/\fc$:
\begin{alignat*}{3}
&\deg=-10\colon \quad &&c_{-10}=\widehat x_{12}\wedge\widehat x_{12},&\\
&\deg=-8\colon \quad &&c_{-8}=\widehat x_{11}\wedge\widehat x_{11},&\\
&\deg=-6\colon \quad &&c_{-6}^1=\widehat x_{8}\wedge\widehat x_{8},&\\
&&&c_{-6}^2=\widehat x_{1}\wedge\widehat x_{12}+\widehat x_{5}\wedge\widehat x_{11}+\widehat x_{8}\wedge\widehat x_{9},&\\
&&&c_{-6}^3=\widehat x_{2}\wedge\widehat x_{12}+\widehat x_{6}\wedge\widehat x_{11}+\widehat x_{8}\wedge\widehat x_{10},&\\
&&&c_{-6}^4=\widehat x_{9}\wedge\widehat x_{9},&\\
&&&c_{-6}^5=\widehat x_{3}\wedge\widehat x_{12}+\widehat x_{7}\wedge\widehat x_{11}+\widehat x_{9}\wedge\widehat x_{10},&\\
&\deg=-4\colon \quad &&c_{-4}^1=\widehat x_{5}\wedge\widehat x_{5},&\\
&&&c_{-4}^2=\widehat x_{4}\wedge\widehat x_{10}+\widehat x_{6}\wedge\widehat x_{7}+\widehat x_{12}\wedge\widehat y_{1},&\\
&&&c_{-4}^3=\widehat x_{4}\wedge\widehat x_{9}+\widehat x_{5}\wedge\widehat x_{7}+\widehat x_{12}\wedge\widehat y_{2},&\\
&&&c_{-4}^4=\widehat x_{4}\wedge\widehat x_{8}+\widehat x_{5}\wedge\widehat x_{6}+\widehat x_{12}\wedge\widehat y_{3},&\\
&\deg=-2\colon \quad &&c_{-2}^1=\widehat x_{1}\wedge\widehat x_{1},&\\
&&&c_{-2}^2=\widehat x_{2}\wedge\widehat x_{3}+\widehat x_{10}\wedge\widehat y_{4}+\widehat x_{11}\wedge\widehat y_{5},&\\
&&&c_{-2}^3=\widehat x_{1}\wedge\widehat x_{3}+\widehat x_{9}\wedge\widehat y_{4}+\widehat x_{11}\wedge\widehat y_{6},&\\
&&&c_{-2}^4=\widehat x_{1}\wedge\widehat x_{2}+\widehat x_{8}\wedge\widehat y_{4}+\widehat x_{11}\wedge\widehat y_{7},&\\
&\deg=0\colon \quad &&c_{0}^1=\widehat x_{2}\wedge\widehat y_{3}+\widehat x_{6}\wedge\widehat y_{7}+\widehat x_{8}\wedge\widehat y_{9},&\\
&&&c_{0}^2=\widehat x_{1}\wedge\widehat y_{3}+\widehat x_{5}\wedge\widehat y_{7}+\widehat x_{8}\wedge\widehat y_{10},&\\
&&&c_{0}^3=\widehat x_{3}\wedge\widehat y_{2}+\widehat x_{7}\wedge\widehat y_{6}+\widehat x_{9}\wedge\widehat y_{8},&\\
&&&c_{0}^4=\widehat x_{1}\wedge\widehat y_{2}+\widehat x_{5}\wedge\widehat y_{6}+\widehat x_{9}\wedge\widehat y_{10},&\\
&&&c_{0}^5=\widehat x_{3}\wedge\widehat y_{1}+\widehat x_{7}\wedge\widehat y_{5}+\widehat x_{10}\wedge\widehat y_{8},&\\
&&&c_{0}^6=\widehat x_{2}\wedge\widehat y_{1}+\widehat x_{6}\wedge\widehat y_{5}+\widehat x_{10}\wedge\widehat y_{9},&\\
&&&c_{0}^7=\widehat x_{2}\wedge\widehat y_{2}+\widehat x_{3}\wedge\widehat y_{3}+\widehat x_{6}\wedge\widehat y_{6}+\widehat x_{7}\wedge\widehat y_{7}+\widehat x_{8}\wedge\widehat y_{8}+\widehat x_{9}\wedge\widehat y_{9},&\\
&&&c_{0}^8=\widehat x_{1}\wedge\widehat y_{1}+\widehat x_{3}\wedge\widehat y_{3}+\widehat x_{5}\wedge\widehat y_{5}+\widehat x_{7}\wedge\widehat y_{7}+\widehat x_{8}\wedge\widehat y_{8}+\widehat x_{10}\wedge\widehat y_{10}.&
\end{alignat*}
\end{enumerate}
\end{Lemma}

\begin{Lemma}\label{2.4.1bO} \quad
\begin{enumerate}\itemsep=0pt
\item[$(a1)$] Let $\fg=\fo\fo_{I\Pi}^{(1)}(3|6)$ with Cartan matrix
\[
\begin{pmatrix}
 \ev & 1 & 0 & 0 \\
 1 & \ev & 1 & 0 \\
 0 & 1 & 0 & 1 \\
 0 & 0 & 1 & \od \end{pmatrix}
 \]
 and basis
\begin{gather*}
x_1,\ x_2,\ x_3,\ x_4,\
x_5=[x_1, x_2],\ x_6=[x_2, x_3],\ x_7=[x_3, x_4],\
x_8=[x_3, [x_1, x_2]],\nonumber\\
 x_9=[x_4, [x_2, x_3]],\ x_{10}=[x_4, [x_3,
x_4]],\
x_{11}=[x_4, [x_4, [x_2, x_3]]],\nonumber\\
 x_{12}=[[x_1, x_2], [x_3, x_4]],\
x_{13}=[x_3, x_4]^2, \ x_{14}=[[x_1, x_2], [x_4, [x_3,
x_4]]],\nonumber\\
x_{15}=[[x_3, x_4], [x_4, [x_2, x_3]]],\
x_{16}=[[x_3, [x_1, x_2]], [x_4, [x_3, x_4]]],\
x_{17}=[x_4, [x_2, x_3]]^2,\nonumber\\
 x_{18}=[[x_4, [x_2, x_3]],
[[x_1, x_2], [x_3, x_4]]],\
x_{19}=[[x_1, x_2], [x_3, x_4]]^2.
\end{gather*}
For a~basis of $H^1(\fg;\fg)$ we can take the following derivations $($recall convention~\eqref{conv}$)$
\begin{gather*}
c_{-2}= x_3\otimes\widehat x_{10} +x_6\otimes\widehat x_{11} + x_8\otimes\widehat x_{14}
+y_4\otimes\widehat x_4 + y_{10}\otimes\widehat y_3 +y_{11}\otimes\widehat y_6
+y_{14}\otimes\widehat y_8.
\end{gather*}

\item[$(a2)$] For a~basis of $H^2(\fg)$ we can take the following cocycles \emph{(Recall convention~\eqref{conv})}
\begin{alignat*}{3}
&\deg=-10 \colon \quad && c_{-10}=\widehat x_{14}\wedge\widehat x_{14}, &\\
&\deg=-8 \colon \quad && c_{-8}=\widehat x_{11}\wedge\widehat x_{11}, &\\
&\deg=-6\colon \quad && c_{-6}^1=\widehat x_{8}\wedge\widehat x_{8}, \
c_{-6}^2=\widehat x_{10}\wedge\widehat x_{10}, &\\
& \deg=-4\colon \quad && c_{-4}=\widehat x_{6}\wedge\widehat x_{6}, &\\
& \deg=-2\colon \quad && c_{-2}=\widehat x_{3}\wedge\widehat x_{3}.&
\end{alignat*}

\item[$(b1)$] Let $\fg=\fo\fo_{I\Pi}^{(1)}(5|4)$ with Cartan matrix
\[
\begin{pmatrix}
 \ev & 1 & 0 & 0 \\
 1 & 0 & 1 & 0 \\
 0 & 1 & \ev& 1 \\
 0 & 0 & 1 & \od \end{pmatrix}
 \]
 and basis
 \begin{gather*}
x_1,\ x_2,\ x_3,\ x_4,\
x_5=[x_1, x_2],\ x_6=[x_2, x_3],\ x_7=[x_3, x_4], \
x_8=[x_3, [x_1, x_2]],\nonumber\\
 x_9=[x_4, [x_2, x_3]],\ x_{10}=[x_4, [x_3, x_4]],\
x_{11}=[x_4, [x_4, [x_2, x_3]]],\nonumber\\
 x_{12}=[[x_1, x_2], [x_3, x_4]],\
x_{13}=[[x_1, x_2], [x_4, [x_3, x_4]]],\ x_{14}=[[x_3, x_4], [x_4,
[x_2, x_3]]],\nonumber\\
x_{15}=[[x_3, [x_1, x_2]], [x_4, [x_3, x_4]]],\
x_{16}=[x_4, [x_2, x_3]]^2,\nonumber\\
 x_{17}=[[x_4, [x_2, x_3]],
[[x_1, x_2], [x_3, x_4]]],\
x_{18}=[[x_1, x_2], [x_3, x_4]]^2.
\end{gather*}
For a~basis of $H^1(\fg;\fg)$ we can take the following derivations $($recall convention~\eqref{conv}$)$
\begin{gather*}
c_{-4}= x_2\otimes\widehat x_{14} +x_5\otimes\widehat x_{15} +y_3\otimes\widehat x_{10}
+ y_7\otimes\widehat x_7 +y_{10}\otimes\widehat x_3 +y_{14}\otimes\widehat y_2 +
y_{15}\otimes\widehat y_5, \\
c_{-2}= x_3\otimes\widehat x_{10} +x_6\otimes\widehat x_{11} +x_8\otimes\widehat x_{13}
+y_4\otimes\widehat x_4 +y_{10}\otimes\widehat y_3 +y_{11}\otimes\widehat y_6 +
y_{13}\otimes\widehat y_8.
\end{gather*}

\item[$(b2)$] For a~basis of $H^2(\fg)$ we can take the following cocycles $($recall convention~\eqref{conv}$)$
\begin{gather*}
c_{-12}=\widehat x_{15}\wedge\widehat x_{15}, \qquad
c_{-10}^1=\widehat x_{14}\wedge\widehat x_{14},\qquad c_{-10}^2=\widehat x_{13}\wedge\widehat x_{13}, \qquad
c_{-8}=\widehat x_{11}\wedge\widehat x_{11}, \\
c_{-6}=\widehat x_{8}\wedge\widehat x_{8}, \qquad
c_{-4}^1=\widehat x_{6}\wedge\widehat x_{6}, \qquad c_{-4}^2=\widehat x_{5}\wedge\widehat x_{5}, \qquad
c_{-2}=\widehat x_{2}\wedge\widehat x_{2}.
 \end{gather*}

\item[$(c1)$] Let $\fg=\fo\fo_{I\Pi}^{(1)}(7|2)$ with Cartan matrix
\[
\begin{pmatrix}
 0 & 1 & 0 & 0 \\
 1 & \ev & 1 & 0 \\
 0 & 1 & \ev & 1 \\
 0 & 0 & 1 & 1\end{pmatrix}
 \]
 and basis
\begin{gather*}
x_1,\ x_2,\ x_3,\ x_4,\
x_5=[x_1, x_2],\ x_6=[x_2, x_3],\ x_7=[x_3, x_4], \ x_8=[x_3, [x_1, x_2]], \nonumber\\
x_9=[x_4, [x_2, x_3]], \ x_{10}=[x_4, [x_3,
x_4]],\
x_{11}=[x_4, [x_4, [x_2, x_3]]],\nonumber\\
 x_{12}=[[x_1, x_2], [x_3, x_4]], \
x_{13}=[[x_1, x_2], [x_4, [x_3, x_4]]],\ x_{14}=[[x_3, x_4], [x_4,
[x_2, x_3]]], \nonumber\\
x_{15}=[[x_3, [x_1, x_2]], [x_4, [x_3, x_4]]],\
x_{16}=[[x_4, [x_2, x_3]], [[x_1, x_2], [x_3, x_4]]],\nonumber\\
x_{17}=[[x_1, x_2], [x_3, x_4]]^2.
\end{gather*}
For a~basis of $H^1(\fg;\fg)$ we can take
the following derivations $($recall convention~\eqref{conv}$)$
\begin{gather*}
c_{-6}= x_1\otimes\widehat x_{16} +y_2\otimes\widehat x_{14} +y_6\otimes
\widehat x_{11} +y_9\otimes\widehat x_9 +y_{11}\otimes\widehat x_6 +y_{14}\otimes\widehat x_2 +y_{16}\otimes\widehat y_1, \\
c_{-4}= x_2\otimes\widehat x_{14} +x_5\otimes\widehat x_{15} +y_3\otimes\widehat x_{10}
+y_7\otimes\widehat x_7 +y_{10}\otimes\widehat x_3 +y_{14}\otimes\widehat y_2 +y_{15}\otimes\widehat y_5, \\
c_{-2}= x_3\otimes\widehat x_{10} +x_6\otimes\widehat x_{11} +x_8\otimes\widehat x_{13}
+y_4\otimes\widehat x_4 +y_{10}\otimes\widehat y_3 +y_{11}\otimes\widehat y_6 +y_{13}\otimes\widehat y_8.
\end{gather*}

\item[$(c2)$] For a~basis of $H^2(\fg)$ we can take $($the classes of$)$ the following cocycles $($recall convention~\eqref{conv}$)$
\begin{gather*}
c_{-6}= \widehat x_{8}\wedge\widehat x_{8}, \qquad
c_{-4}=\widehat x_{5}\wedge\widehat x_{5}, \qquad
c_{-2}=\widehat x_{1}\wedge\widehat x_{1}.
\end{gather*}

\item[$(d1)$] For $\fg=\fo\fo^{(1)}_{I\Pi}(1|8)$, we have $H^1(\fg;\fg)=0$ and $H^2(\fg)=0$.
\end{enumerate}
\end{Lemma}

\subsection{Rank 5}

\begin{Lemma}\label{2.5.1cO)} \quad
\begin{enumerate}\itemsep=0pt
\item[$(a1)$] Let $\fg=\fo\fo^{(1)}_{I\Pi}(7|4)$ with the following Cartan matrix
\[
\begin{pmatrix}
 \ev & 1 & 0 & 0 & 0 \\
 1 & 0 & 1 & 0 & 0 \\
 0 & 1 & \ev & 1 & 0 \\
 0 & 0 & 1 & \ev & 1 \\
 0 & 0 & 0 & 1 &\od \end{pmatrix}
 \]
 and
basis
\begin{gather*}
x_1,\ x_2,\ x_3,\ x_4,\ x_5,\
x_6=[x_1, x_2],\ x_7=[x_2, x_3],\ x_8=[x_3, x_4],\ x_9=[x_4,
x_5],\nonumber\\
x_{10}=[x_3, [x_1, x_2]],\ x_{11}=[x_4, [x_2, x_3]],\ x_{12}=[x_5,
[x_3,
x_4]],\
 x_{13}=[x_5, [x_4, x_5]],\nonumber\\
x_{14}=[x_5, [x_5, [x_3, x_4]]],\ x_{15}=[[x_1, x_2], [x_3,
x_4]],\ x_{16}=[[x_2, x_3], [x_4, x_5]],\nonumber\\
x_{17}=[[x_2, x_3], [x_5, [x_4, x_5]]],\ x_{18}=[[x_4, x_5], [x_3,
[x_1, x_2]]],\nonumber\\
 x_{19}=[[x_4, x_5], [x_5, [x_3, x_4]]],\
x_{20}=[[x_3, [x_1, x_2]], [x_5, [x_4, x_5]]],\nonumber\\
x_{21}=[[x_4, [x_2,
x_3]], [x_5, [x_4, x_5]]],\
x_{22}=[[x_5, [x_3, x_4]], [[x_2, x_3], [x_4, x_5]]],\nonumber\\
 x_{23}=[[x_5,
[x_4, x_5]], [[x_1, x_2], [x_3, x_4]]],\
x_{24}=[[x_5, [x_5, [x_3, x_4]]], [[x_1, x_2], [x_3, x_4]]]\nonumber\\
x_{25}=[[x_2, x_3], [x_4,
x_5]]^2,\
x_{26}=[[[x_2, x_3], [x_4, x_5]], [[x_4, x_5], [x_3, [x_1,
x_2]]]],\nonumber\\
x_{27}=[[x_4, x_5], [x_3, [x_1, x_2]]]^2.
\end{gather*}
For a~basis of $H^1(\fg;\fg)$ we can take the following derivations $($recall convention~\eqref{conv}$)$
\begin{gather*}
c_{-6}= x_2\otimes\widehat x_{22} +x_6\otimes\widehat x_{24} + y_3\otimes\widehat x_{19}
+y_8\otimes\widehat x_{14} + y_{12}\otimes\widehat x_{12} +y_{14}\otimes \widehat
x_8 +y_{19}\otimes\widehat x_3 \\
\hphantom{c_{-6}=}{}
+
y_{22}\otimes\widehat y_2 +y_{24}\otimes\widehat y_6 ,\\
c_{-4}= x_3\otimes\widehat x_{19} +x_7\otimes\widehat x_{21} + x_{10}\otimes
\widehat x_{23} +y_4\otimes\widehat x_{13} + y_9\otimes\widehat x_9 +y_{13}\otimes\widehat x_4 + y_{19}\otimes\widehat y_3 \\
\hphantom{c_{-4}=}{}
+y_{21}\otimes\widehat y_7 +
y_{23}\otimes\widehat y_{10},\\
c_{-2}= x_4\otimes\widehat x_{13} +x_8\otimes\widehat x_{14} + x_{11}\otimes
\widehat x_{17} +x_{15}\otimes\widehat x_{20} + y_5\otimes\widehat x_5 +y_{13}\otimes \widehat
y_4 +y_{14}\otimes\widehat y_8 \\
\hphantom{c_{-2}=}{}
+ y_{17}\otimes\widehat y_{11} +y_{20}\otimes\widehat y_{15}.
\end{gather*}

\item[$(a2)$] For a~basis of $H^2(\fg)$ we can take the following cocycles $($recall convention~\eqref{conv}$)$
\begin{alignat*}{5}
&\deg=-16\colon \quad && \widehat x_{24}\wedge\widehat x_{24},\qquad && \deg=-8\colon \quad &&\widehat x_{15}\wedge\widehat x_{15},&\\
&\deg=-14\colon \quad &&\widehat x_{22}\wedge\widehat x_{22},\
\widehat x_{23}\wedge\widehat x_{23},\qquad && \deg=-6\colon \quad &&\widehat x_{10}\wedge\widehat x_{10},\
\widehat x_{11}\wedge\widehat x_{11},& \\
& \deg=-12\colon \quad &&\widehat x_{20}\wedge\widehat x_{20},\
\widehat x_{21}\wedge\widehat x_{21},\qquad && \deg=-4\colon \quad &&\widehat x_{7}\wedge\widehat x_{7},\
\widehat x_{6}\wedge\widehat x_{6},& \\
&\deg=-10\colon \quad &&\widehat x_{17}\wedge\widehat x_{17},\qquad && \deg=-2\colon \quad && \widehat x_{2}\wedge\widehat x_{2}.&
\end{alignat*}

\item[$(b1)$] Let $\fg=\fo\fo^{(1)}_{I\Pi}(9|2)$ with the following Cartan matrix
\[
\begin{pmatrix}
 0 & 1 & 0 & 0 & 0 \\
 1 & \ev & 1 & 0 & 0 \\
 0 & 1 & \ev & 1 & 0 \\
 0 & 0 & 1 & \ev & 1 \\
 0 & 0 & 0 & 1 &\od \end{pmatrix}
\]
 and basis
\begin{gather*}
x_1,\ x_2,\ x_3,\ x_4,\ x_5,\
x_6=[x_1, x_2],\ x_7=[x_2, x_3],\ x_8=[x_3, x_4],\ x_9=[x_4,
x_5],\nonumber\\
x_{10}=[x_3, [x_1, x_2]],\ x_{11}=[x_4, [x_2, x_3]],\ x_{12}=[x_5,
[x_3,
x_4]],\
 x_{13}=[x_5, [x_4, x_5]],\nonumber\\
x_{14}=[x_5, [x_5, [x_3, x_4]]],\ x_{15}=[[x_1, x_2], [x_3,
x_4]],\ x_{16}=[[x_2, x_3], [x_4, x_5]],\nonumber\\
x_{17}=[[x_2, x_3], [x_5, [x_4, x_5]]],\ x_{18}=[[x_4, x_5], [x_3,
[x_1, x_2]]],\nonumber\\
 x_{19}=[[x_4, x_5], [x_5, [x_3, x_4]]],\
x_{20}=[[x_3, [x_1, x_2]], [x_5, [x_4, x_5]]],\nonumber\\
x_{21}=[[x_4, [x_2,
x_3]], [x_5, [x_4, x_5]]],\
x_{22}=[[x_5, [x_3, x_4]], [[x_2, x_3], [x_4, x_5]]],\nonumber\\
 x_{23}=[[x_5,
[x_4, x_5]], [[x_1, x_2], [x_3, x_4]]],\
x_{24}=[[x_5, [x_5, [x_3, x_4]]], [[x_1, x_2], [x_3, x_4]]],\nonumber\\
x_{25}=[[x_2, x_3], [x_4,
x_5]],[[x_4,x_5],[x_3,[x_1,x_2]]]],\
x_{26}=[[x_4,x_5],[x_3,[x_1,x_2]]]^2.
\end{gather*}
For a~basis of $H^1(\fg;\fg)$ we can take the following derivations $($recall convention~\eqref{conv}$)$
\begin{gather*}
c_{-8}= y_2 \otimes\widehat x_{22} +y_7 \otimes\widehat x_{21}+y_{11} \otimes\widehat x_{17}+
 y_{16} \otimes\widehat x_{16} +y_{17} \otimes\widehat x_{11} +y_{21} \otimes\widehat x_7 +y_{22} \otimes\widehat x_2 \\
\hphantom{c_{-8}=}{}
 +x_1\otimes \widehat x_{25}+y_{25} \otimes\widehat y_1,
 \\
c_{-6}= x_2\otimes\widehat x_{22} +x_6\otimes\widehat x_{24} + y_3\otimes\widehat x_{19}
+y_8\otimes\widehat x_{14} + y_{12}\otimes\widehat x_{12} +y_{14}\otimes d
x_8+y_{19}\otimes\widehat x_3 \\
\hphantom{c_{-6}=}{}
+
y_{22}\otimes\widehat y_2 +y_{24}\otimes\widehat y_6 ,\\
c_{-4}= x_3\otimes\widehat x_{19} +x_7\otimes\widehat x_{21} + x_{10}\otimes
\widehat x_{23} +y_4\otimes\widehat x_{13} + y_9\otimes\widehat x_9 +y_{13}\otimes\widehat x_4+ y_{19}\otimes\widehat y_3\\
\hphantom{c_{-4}=}{}
 +y_{21}\otimes\widehat y_7 +
y_{23}\otimes\widehat y_{10} ,\\
c_{-2}= x_4\otimes\widehat x_{13} +x_8\otimes\widehat x_{14} + x_{11}\otimes
\widehat x_{17} +x_{15}\otimes\widehat x_{20} + y_5\otimes\widehat x_5 +y_{13}\otimes d
y_4+y_{14}\otimes\widehat y_8 \\
\hphantom{c_{-2}=}{}
+ y_{17}\otimes\widehat y_{11} +y_{20}\otimes\widehat y_{15}.
\end{gather*}

\item[$(b2)$] For a~basis of $H^2(\fg)$ we can take the following cocycles $($recall convention~\eqref{conv}$)$
\begin{alignat*}{5}
& \deg=-18\colon \quad && \widehat x_{25}\wedge\widehat x_{25},\qquad && \deg=-8\colon \quad &&\widehat x_{15}\wedge\widehat x_{15},&\\
&\deg=-16\colon \quad &&\widehat x_{24}\wedge\widehat x_{24},\qquad && \deg=-6\colon \quad &&\widehat x_{10}\wedge\widehat x_{10},&\\
&\deg=-14\colon \quad &&\widehat x_{23}\wedge\widehat x_{23},\qquad && \deg=-4\colon \quad &&\widehat x_{6}\wedge\widehat x_{6},&\\
&\deg=-12\colon \quad && \widehat x_{20}\wedge\widehat x_{20},\qquad && \deg=-2\colon \quad &&\widehat x_{1}\wedge\widehat x_{1}.&
\end{alignat*}

\item[$(c1)$] Let $\fg=\fo\fo^{(1)}_{I\Pi}(5|6)$ with the following Cartan matrix
\[
\begin{pmatrix}
 \ev & 1 & 0 & 0 & 0 \\
 1 & \ev & 1 & 0 & 0 \\
 0 & 1 & 0 & 1 & 0 \\
 0 & 0 & 1 & \ev & 1 \\
 0 & 0 & 0 & 1 &\od \end{pmatrix}
 \]
 and
basis
\begin{gather*}
x_1,\ x_2,\ x_3,\ x_4,\ x_5,\
x_6=[x_1, x_2],\ x_7=[x_2, x_3],\ x_8=[x_3, x_4],\ x_9=[x_4,
x_5],\nonumber\\
x_{10}=[x_3, [x_1, x_2]],\ x_{11}=[x_4, [x_2, x_3]],\ x_{12}=[x_5,[x_3,
x_4]],\
 x_{13}=[x_5, [x_4, x_5]],\nonumber\\
x_{14}=[x_5, [x_5, [x_3, x_4]]],\ x_{15}=[[x_1, x_2], [x_3,
x_4]],\ x_{16}=[[x_2, x_3], [x_4, x_5]],\nonumber\\
x_{17}=[[x_2, x_3], [x_5, [x_4, x_5]]],\ x_{18}=[[x_4, x_5], [x_3,
[x_1, x_2]]],\nonumber\\
 x_{19}=[[x_4, x_5], [x_5, [x_3, x_4]]],\
x_{20}=[[x_3, [x_1, x_2]], [x_5, [x_4, x_5]]],\nonumber\\
x_{21}=[[x_4, [x_2,
x_3]], [x_5, [x_4, x_5]]],\
x_{22}=[[x_5, [x_3, x_4]]^2,\nonumber\\
 x_{23}=[[x_5,
[x_3, x_4]], [[x_2, x_3], [x_4, x_5]]],\
x_{24}=[ [x_5, [x_4, x_5]], [[x_1, x_2], [x_3, x_4]]],\nonumber\\
x_{25}=[[x_5,[x_5,[x_3,x_4]]], [[x_1,x_2], [x_3,x_4]]],\
x_{26}=[[x_2,x_3],[x_4,x_5]]^2,\nonumber\\
x_{27}=[[[x_2,x_3], [x_4,x_5]], [[x_4,x_5], [x_3,[x_1,x_2]]]],\
x_{28}=[[x_4,x_5], [x_3, [x_1,x_2]]]^2.
\end{gather*}
For a~basis of $H^1(\fg;\fg)$ we can take the following derivations $($recall convention~\eqref{conv}$)$
\begin{gather*}
c_{-4}= x_3\otimes\widehat x_{19} +x_7\otimes\widehat x_{21} + x_{10}\otimes
\widehat x_{24} +y_4\otimes\widehat x_{13} + y_9\otimes\widehat x_9 +y_{13}\otimes\widehat x_4+ y_{19}\otimes\widehat y_3\\
\hphantom{c_{-4}=}{}
 +y_{21}\otimes\widehat y_7 +
y_{24}\otimes\widehat y_{10} ;\\
c_{-2}= x_4\otimes\widehat x_{13} +x_8\otimes\widehat x_{14} + x_{11}\otimes
\widehat x_{17} +x_{15}\otimes\widehat x_{20} + y_5\otimes\widehat x_5 +y_{13}\otimes \widehat
y_4+y_{14}\otimes\widehat y_8 \\
\hphantom{c_{-2}=}{}
+ y_{17}\otimes\widehat y_{11} +y_{20}\otimes\widehat y_{15}.
\end{gather*}

\item[$(c2)$] For a~basis of $H^2(\fg)$ we can take the following cocycles $($recall convention~\eqref{conv}$)$
\begin{alignat*}{5}
&\deg=-14\colon \quad && \widehat x_{24}\wedge\widehat x_{24},\qquad && \deg=-6\colon \quad && \widehat x_{10}\wedge\widehat x_{10},\ \widehat x_{11}\wedge\widehat x_{11},&\\
& \deg=-12\colon \quad && \widehat x_{20}\wedge\widehat x_{20},\ \widehat x_{21}\wedge\widehat x_{21},\qquad &&\deg=-4\colon \quad && \widehat x_{7}\wedge\widehat x_{7},\
\widehat x_{8}\wedge\widehat x_{8},& \\
&\deg=-10\colon \quad && \widehat x_{17}\wedge\widehat x_{17},\ \widehat x_{19}\wedge\widehat x_{19},\qquad && \deg=-2\colon \quad && \widehat x_{3}\wedge\widehat x_{3}.&\\
&\deg=-8\colon \quad && \widehat x_{15}\wedge\widehat x_{15},\ \widehat x_{14}\wedge\widehat x_{14},\qquad && && &
\end{alignat*}

\item[$(d1)$] Let $\fg=\fo\fo^{(1)}_{I\Pi}(3|8)$ with the following Cartan matrix
\[
\begin{pmatrix}
 \ev & 1 & 0 & 0 & 0 \\
 1 & \ev & 1 & 0 & 0 \\
 0 & 1 & \ev & 1 & 0 \\
 0 & 0 & 1 & 0 & 1 \\
 0 & 0 & 0 & 1 &\od \end{pmatrix}
 \]
 and
basis
\begin{gather*}
x_1,\ x_2,\ x_3,\ x_4,\ x_5,\
x_6=[x_1, x_2],\ x_7=[x_2, x_3],\ x_8=[x_3, x_4],\ x_9=[x_4,
x_5],\nonumber\\
x_{10}=[x_3, [x_1, x_2]],\ x_{11}=[x_4, [x_2, x_3]],\ x_{12}=[x_5,
[x_3,
x_4]],\
 x_{13}=[x_5, [x_4, x_5]],\nonumber\\
x_{14}=[x_5, [x_5, [x_3, x_4]]],\ x_{15}=[[x_1, x_2], [x_3,
x_4]],\ x_{16}=[[x_2, x_3], [x_4, x_5]],\nonumber\\
x_{17}= [x_4, x_5]^2,\ x_{18}=[[x_2, x_3], [x_5,
[x_4, x_5]]],\
 x_{19}=[[x_4, x_5], [x_3, [x_1, x_2]]],\nonumber\\
x_{20}=[[[x_4,x_5], [x_5,[x_3,x_4]]], \ x_{21}=[[x_3, [x_1,
x_2]], [x_5, [x_4, x_5]]],\nonumber\\
x_{22}=[[x_4, [x_2, x_3]] , [x_5, [x_4,x_5]]],\
 x_{23}=[[x_5,
[x_3, x_4]]^2,\nonumber\\
x_{24}=[ [x_5, [x_3, x_4]], [[x_2, x_3], [x_4, x_5]]],\
x_{25}=[ [x_5, [x_4, x_5]], [[x_1, x_2], [x_3, x_4]]],\nonumber\\
x_{26}=[[x_5,[x_5,[x_3,x_4]]], [[x_1,x_2], [x_3,x_4]]],\
x_{27}=[[x_2,x_3], [x_4,x_5]]^2,\nonumber\\
x_{28}=[[[x_2,x_3], [x_4,x_5]], [[x_4,x_5], [x_3,[x_1,x_2]]]].
\end{gather*}
For a~basis of $H^1(\fg;\fg)$ we can take the following derivations $($recall convention~\eqref{conv}$)$
\begin{gather*}
c_{-2}= x_4\otimes\widehat x_{13} +x_8\otimes\widehat x_{14} + x_{11}\otimes
\widehat x_{18} +x_{15}\otimes\widehat x_{21} + y_5\otimes\widehat x_5 +y_{13}\otimes \widehat
y_4\\
\hphantom{c_{-2}=}{}
+y_{14}\otimes\widehat y_8 + y_{17}\otimes\widehat y_{11} +y_{21}\otimes\widehat y_{15}.
\end{gather*}

\item[$(d2)$] For a~basis of $H^2(\fg)$ we can take the following cocycles $($recall convention~\eqref{conv}$)$
\begin{alignat*}{5}
& \deg=-12\colon \quad && \widehat x_{21}\wedge\widehat x_{21},\qquad && \deg=-6\colon \quad && \widehat x_{11}\wedge\widehat x_{11},\ \widehat x_{13}\wedge\widehat x_{13},&\\
& \deg=-10\colon \quad && \widehat x_{18}\wedge\widehat x_{18},\qquad && \deg=-4\colon \quad && \widehat x_{8}\wedge\widehat x_{8}, & \\
& \deg=-8\colon \quad && \widehat x_{15}\wedge\widehat x_{15},\ \widehat x_{14}\wedge\widehat x_{14},\qquad && \deg=-2\colon \quad && \widehat x_{4}\wedge\widehat x_{4}.&
\end{alignat*}

\item[$(e)$] For $\fg=\fo\fo_{I\Pi}^{(1)}(1|10)$, we have $H^2(\fg)=0$ and $H^1(\fg;\fg)=0$.
\end{enumerate}
\end{Lemma}

\begin{Lemma}\label{2.5.1aO} \quad
\begin{enumerate}\itemsep=0pt
\item[$(a)$] For $\fg=\fo\fc(2;10)\ltimes
\Kee I_0$,
with the following Cartan matrix
\[
\begin{pmatrix}
 \ev & 0 & 1 & 0 & 0 \\
 0 & \ev & 1 & 0 & 0 \\
 1 & 1 & \ev & 1 & 0 \\
 0 & 0 & 1 & \ev & 1 \\
 0 & 0 & 0 & 1 & \ev\end{pmatrix}
 \]
 and basis
\begin{gather*}
x_1,\ x_2,\ x_3,\ x_4,\ x_5,\
x_6=[x_1, x_3], \ x_7=[x_2, x_3],\ x_8=[x_3, x_4],\ x_9=[x_4,
x_5],\nonumber\\
x_{10}=[x_2, [x_1, x_3]],\ x_{11}=[x_4, [x_1, x_3]],\ x_{12}=[x_4, [x_2,
x_3]],\ x_{13}=[x_5, [x_3, x_4]],\nonumber\\
x_{14}=[x_4, [x_2, [x_1, x_3]]],\ x_{15}=[[x_1, x_3], [x_4,
x_5]],\ x_{16}=[[x_2, x_3], [x_4, x_5]],\nonumber\\
x_{17}=[[x_3, x_4], [x_2, [x_1, x_3]]],\ x_{18}=[[x_4, x_5], [x_2,
[x_1, x_3]]],\nonumber\\
x_{19}=[[x_2, [x_1, x_3]], [x_5, [x_3, x_4]]],\
x_{20}=[[x_5, [x_3,x_4]], [x_4,[x_2, [x_1, x_3]]]].
\end{gather*} as well as for $\fo
\fo \fc(2;8|2)\ltimes \Kee I_0$, and $\fo\fo \fc(2;6|4)\ltimes
\Kee I_0$, and $\fg=\mathfrak{pec}(2;5)\ltimes \Kee I_0$, with the same Cartan matrix, but $0$ instead of $\ev$
in the respective places on the main diagonal, for a~basis of $H^1(\fg;\fg)$ we can take
the following derivation:
\[
\begin{array}{ll}
c_0=h_1\otimes \widehat D_6 + h_2\otimes \widehat D_6,
\end{array}
\]
$($see Remark $\ref{RemD})$.
We have $H^2(\fg)=0$.

\item[$(b)$] For $\fg=\fo \fc(2; 10)/\fc$, $\fo \fo \fc (2; 8|2)/\fc$, and $\fo\fo \fc
(2; 6|4)/\fc$, as well as $\fg=\mathfrak{pec}(2; 5)/\fc$, for a~basis of
$H^1(\fg;\fg)$ we can take the following derivations $($recall convention~\eqref{conv}$)$
\begin{gather*}
c_{-8}= y_5\otimes \widehat x_{20} +y_9\otimes \widehat x_{19} +y_{13}\otimes \widehat
x_{18} +y_{15}\otimes \widehat x_{16} +y_{16}\otimes \widehat x_{15}
+y_{18}\otimes \widehat x_{13}\\
\hphantom{c_{-8}=}{}
+y_{19}\otimes \widehat x_9 +y_{20}\otimes \widehat x_5, \\
c_{-6}= x_5\otimes \widehat x_{20} +y_4\otimes \widehat x_{17} +y_8\otimes \widehat
x_{14} +y_{11}\otimes \widehat x_{12} +y_{12}\otimes \widehat x_{11}
+y_{14}\otimes \widehat x_8\\
\hphantom{c_{-6}=}{}
+y_{17}\otimes \widehat x_4 +y_{20}\otimes \widehat y_5, \\
c_{-4}= x_4\otimes \widehat x_{17} +x_9\otimes \widehat x_{19} +y_3\otimes \widehat
x_{10} +y_6\otimes \widehat x_7 +y_7\otimes \widehat x_6 +y_{10}\otimes \widehat x_3\\
\hphantom{c_{-4}=}{}
+y_{17}\otimes \widehat y_4 +y_{19}\otimes \widehat y_9, \\
c_{-2}= x_3\otimes \widehat x_{10} +x_8\otimes \widehat x_{14} +x_{13}\otimes \widehat
x_{18} +y_1\otimes \widehat x_2 +y_2\otimes \widehat x_1 +y_{10}\otimes \widehat y_3\\
\hphantom{c_{-2}=}{}
+y_{14}\otimes \widehat y_8 +y_{18}\otimes \widehat y_{13}, \\
c_{0}^1= x_2\otimes \widehat x_1 +x_7\otimes \widehat x_6 +x_{12}\otimes \widehat x_{11}
+x_{16}\otimes \widehat x_{15} +y_1\otimes \widehat y_2 +y_6\otimes \widehat y_7\\
\hphantom{c_{0}^1=}{}
+y_{11}\otimes \widehat y_{12} +y_{15}\otimes \widehat y_{16}, \\
c_{0}^2= x_1\otimes \widehat x_2 +x_6\otimes \widehat x_7 +x_{11}\otimes \widehat x_{12}
+x_{15}\otimes \widehat x_{16}
+y_2\otimes \widehat y_1 +y_7\otimes \widehat y_6 \\
\hphantom{c_{0}^2=}{}
+y_{12}\otimes \widehat y_{11} +y_{16}\otimes \widehat y_{15}, \\
c_{0}^3= x_2\otimes \widehat x_2 +x_4\otimes \widehat x_4 +x_7\otimes \widehat x_7
+x_8\otimes \widehat x_8 +x_9\otimes \widehat x_9 +x_{10}\otimes \widehat x_{10}
+x_{11}\otimes \widehat x_{11} \\
\hphantom{c_{0}^3=}{}
+x_{13}\otimes \widehat x_{13} +
 x_{15}\otimes \widehat x_{15} +x_{20}\otimes \widehat x_{20} +y_2\otimes \widehat y_2
+y_4\otimes \widehat y_4 +y_7\otimes \widehat y_7 +y_8\otimes \widehat y_8\\
\hphantom{c_{0}^3=}{}
+y_9\otimes \widehat
y_9 +y_{10}\otimes \widehat y_{10}
 +y_{11}\otimes \widehat y_{11} +y_{13}\otimes \widehat y_{13} +y_{15}\otimes \widehat
y_{15} +y_{20}\otimes \widehat y_{20}.
\end{gather*}

\item[$(c)$] For $\fg=\fo \fo\fc (2; 6|4)/\fc$ and parities of the Chevalley
generators being $(0,0,1,0,0)$, for a~basis of $H^2(\fg)$ we can take the following cocycles $($recall convention~\eqref{conv}$)$
\begin{gather*}
c_{-10}= \widehat x_{18} \wedge\widehat x_{18}, \\
c_{-8}^1= \widehat x_{14} \wedge\widehat x_{14}, \
c_{-8}^2= \widehat x_{15} \wedge\widehat x_{15}, \
c_{-8}^3= \widehat x_{16}\wedge\widehat x_{16},\\
c_{-8}^4=
\widehat x_5 \wedge \widehat x_{20}+ \widehat x_9 \wedge
 \widehat x_{19}+\widehat x_{13}\wedge \widehat x_{18}+ \widehat x_{15}\wedge \widehat x_{16}, \\
c_{-6}^1= \widehat x_{11} \wedge\widehat x_{11}, \
c_{-6}^2= \widehat x_{12} \wedge\widehat x_{12}, \
c_{-6}^3= \widehat x_{10} \wedge\widehat x_{10}, \
c_{-6}^4= \widehat x_{13} \wedge\widehat x_{13}, \\
c_{-6}^5= \widehat x_4 \wedge \widehat x_{17}+ \widehat x_8 \wedge \widehat x_{14}+ \widehat x_{11}\wedge
 \widehat x_{12}+
 \widehat x_{20} \wedge \widehat y_5, \\
c_{-4}^1= \widehat x_{6} \wedge\widehat x_{6}, \
c_{-4}^2= \widehat x_{7} \wedge\widehat x_{7}, \
c_{-4}^3= \widehat x_{3} \wedge\widehat x_{10}+\widehat x_{6} \wedge\widehat x_{7}+\widehat x_{17} \wedge\widehat y_{4}+\widehat x_{19} \wedge\widehat y_{9}, \\
c_{-2}^1= \widehat x_{3} \wedge\widehat x_{3}, \
c_{-2}^2= \widehat x_{1} \wedge\widehat x_{2}+\widehat x_{10} \wedge\widehat y_{3}+\widehat x_{14} \wedge\widehat y_{8}+\widehat x_{18} \wedge\widehat y_{13}, \\
c_{0}^1= \widehat x_{1} \wedge\widehat y_{2}+\widehat x_{6} \wedge\widehat y_{7}+\widehat x_{11} \wedge\widehat y_{12}+\widehat x_{15} \wedge\widehat y_{16},
\\
c_{0}^2= \widehat x_{2} \wedge\widehat y_{1}+\widehat x_{7} \wedge\widehat y_{6}+\widehat x_{12} \wedge\widehat y_{11}+\widehat x_{16} \wedge\widehat y_{15}, \\
c_{0}^3= \widehat x_{2} \wedge\widehat y_{2}+\widehat x_{5} \wedge\widehat y_{5}+\widehat x_{7} \wedge\widehat y_{7}+\widehat x_{9} \wedge\widehat y_{9}+\widehat x_{10} \wedge\widehat y_{10}+\widehat x_{12} \wedge\widehat y_{12}+\widehat x_{13} \wedge\widehat y_{13}\\
\hphantom{c_{0}^3=}{}
+\widehat x_{14} \wedge\widehat y_{14}+\widehat x_{15} \wedge\widehat y_{15}+\widehat x_{17} \wedge\widehat y_{17}.
\end{gather*}

\item[$(d)$] For $\fg=\fo \fo\fc (2; 8|2)/\fc$ and parities of the Chevalley
generators being $(0,0,0,0,1)$, for a~basis of $H^2(\fg)$ we can take the following cocycles $($recall convention~\eqref{conv}$)$
\begin{gather*}
c_{-14}= \widehat x_{20} \wedge\widehat x_{20}, \\
c_{-12}= \widehat x_{19} \wedge\widehat x_{19}, \\
c_{-10}= \widehat x_{18} \wedge\widehat x_{18}, \\
c_{-8}^1= \widehat x_{15} \wedge\widehat x_{15}, \
c_{-8}^2= \widehat x_{16}\wedge\widehat x_{16},\
c_{-8}^3= \widehat x_5 \wedge \widehat x_{20}+ \widehat x_9 \wedge
\widehat x_{19}+
\widehat x_{13}\wedge \widehat x_{18}+ \widehat x_{15}\wedge\widehat x_{16}, \\
c_{-6}^1= \widehat x_4 \wedge \widehat x_{17}+ \widehat x_8 \wedge \widehat x_{14}+ \widehat x_{11}\wedge
 \widehat x_{12}+
 \widehat x_{20} \wedge \widehat y_5, \
c_{-6}^2= \widehat x_{13} \wedge\widehat x_{13}, \\
c_{-4}^1= \widehat x_{6} \wedge\widehat x_{6}, \
c_{-4}^2= \widehat x_{9} \wedge\widehat x_{9}, \
c_{-4}^3= \widehat x_{3} \wedge\widehat x_{10}+\widehat x_{6} \wedge\widehat x_{7}+\widehat x_{17} \wedge\widehat y_{4}+\widehat x_{19} \wedge\widehat y_{9}, \\
c_{-2}^1= \widehat x_{1} \wedge\widehat x_{2}+\widehat x_{10} \wedge\widehat y_{3}+\widehat x_{14} \wedge\widehat y_{8}+\widehat x_{18} \wedge\widehat y_{13}, \
c_{-2}^2= \widehat x_{5} \wedge\widehat x_{5}, \\
c_{0}^1= \widehat x_{1} \wedge\widehat y_{2}+\widehat x_{6} \wedge\widehat y_{7}+\widehat x_{11} \wedge\widehat y_{12}+\widehat x_{15} \wedge\widehat y_{16}, \
c_{0}^2= \widehat x_{2} \wedge\widehat y_{1}+\widehat x_{7} \wedge\widehat y_{6}+\widehat x_{12} \wedge\widehat y_{11}\\
\hphantom{c_{0}^2=}{}
+\widehat x_{16} \wedge\widehat y_{15}, \\
c_{0}^3= \widehat x_{2} \wedge\widehat y_{2}+\widehat x_{5} \wedge\widehat y_{5}+\widehat x_{7} \wedge\widehat y_{7}+\widehat x_{9} \wedge\widehat y_{9}+\widehat x_{10} \wedge\widehat y_{10}+\widehat x_{12} \wedge\widehat y_{12}+\widehat x_{13} \wedge\widehat y_{13}\\
\hphantom{c_{0}^3=}{}
+\widehat x_{14} \wedge\widehat y_{14}+\widehat x_{15} \wedge\widehat y_{15}+\widehat x_{17} \wedge\widehat y_{17}.
\end{gather*}

\item[$(e)$] For $\fg=\mathfrak{pec} (2; 5)/\fc$ and parities of the Chevalley
generators being $(1,0,0,0,0)$, for a~basis of $H^2(\fg)$ we can take the following cocycles $($recall convention~\eqref{conv}$)$
\begin{gather*}
c_{-14}= \widehat x_{20} \wedge\widehat x_{20}, \
c_{-12}= \widehat x_{19} \wedge\widehat x_{19}, \
c_{-10}= \widehat x_{18} \wedge\widehat x_{17}, \\
c_{-8}^1= \widehat x_{14} \wedge\widehat x_{14}, \
c_{-8}^2= \widehat x_{15} \wedge\widehat x_{15}, \
c_{-8}^3= \widehat x_5 \wedge \widehat x_{20}+ \widehat x_9 \wedge
\widehat x_{19}+
\widehat x_{13}\wedge \widehat x_{18}+ \widehat x_{15}\wedge \widehat x_{16}, \\
c_{-6}^1= \widehat x_{11} \wedge\widehat x_{11}, \ c_{-6}^2= \widehat x_{10} \wedge\widehat x_{10}, \
c_{-6}^3= \widehat x_4 \wedge \widehat x_{17}+ \widehat x_8 \wedge \widehat x_{14}+ \widehat x_{11}\wedge
 \widehat x_{12}+
\widehat x_{20} \wedge \widehat y_5, \\
c_{-4}^1= \widehat x_{6} \wedge\widehat x_{6}, \
c_{-4}^2= \widehat x_{3} \wedge\widehat x_{10}+\widehat x_{6} \wedge\widehat x_{7}+\widehat x_{17} \wedge\widehat y_{4}+\widehat x_{19} \wedge\widehat y_{9}, \\
c_{-2}^1= \widehat x_{1} \wedge\widehat x_{2}+\widehat x_{10} \wedge\widehat y_{3}+\widehat x_{14} \wedge\widehat y_{8}+\widehat x_{18} \wedge\widehat y_{13}, \
c_{-2}^2= \widehat x_{1} \wedge\widehat x_{1}, \\
c_{0}^1= \widehat x_{1} \wedge\widehat y_{2}+\widehat x_{6} \wedge\widehat y_{7}+\widehat x_{11} \wedge\widehat y_{12}+\widehat x_{15} \wedge\widehat y_{16}, \\
c_{0}^2= \widehat x_{2} \wedge\widehat y_{1}+\widehat x_{7} \wedge\widehat y_{6}+\widehat x_{12} \wedge\widehat y_{11}+\widehat x_{16} \wedge\widehat y_{15}, \\
c_{0}^3= \widehat x_{2} \wedge\widehat y_{2}+\widehat x_{5} \wedge\widehat y_{5}+\widehat x_{7} \wedge\widehat y_{7}+
\widehat x_{9} \wedge\widehat y_{9}+\widehat x_{10} \wedge\widehat y_{10}+\widehat x_{12} \wedge\widehat y_{12}+\widehat x_{13} \wedge\widehat y_{13}\\
\hphantom{c_{0}^3=}{} +
\widehat x_{14} \wedge\widehat y_{14}+\widehat x_{15} \wedge\widehat y_{15}+\widehat x_{17} \wedge\widehat y_{17}.
\end{gather*}

\item[$(f)$] For $\fg=\fo \fc(2; 10)/\fc$, for a~basis of $H^2(\fg)$ we can take the following cocycles $($recall convention~\eqref{conv}$)$
\begin{gather*}
c_{-8}^1= \widehat x_5 \wedge\widehat x_{20}+ \widehat x_9 \wedge
 \widehat x_{19}+
 \widehat x_{13}\wedge \widehat x_{18}+ \widehat x_{15}\wedge
 \widehat x_{16}, \\
c_{-6}^1= \widehat x_4 \wedge \widehat x_{17}+ \widehat x_8 \wedge
 \widehat x_{14}+
 \widehat x_{11}\wedge \widehat x_{12}+ \widehat x_{20}\wedge
 \widehat y_5, \\
c_{-4}^1= \widehat x_3 \wedge \widehat x_{10}+ \widehat x_6 \wedge \widehat x_7+ \widehat x_{17}\wedge
 \widehat y_4+
 \widehat x_{19} \wedge \widehat y_9, \\
c_{-2}^1= \widehat x_1 \wedge \widehat x_2+
 \widehat x_{10} \wedge \widehat y_3+ \widehat x_{14} \wedge
 \widehat y_8 +
 \widehat x_{18} \wedge \widehat y_{13}, \\
c_{0}^1= \widehat x_1 \wedge \widehat y_2+
 \widehat x_6 \wedge \widehat y_7+
 \widehat x_{11} \wedge \widehat y_{12}+ \widehat x_{15}\wedge
 \widehat y_{16}, \\
c_{0}^2= \widehat x_2 \wedge \widehat y_1+
 \widehat x_7 \wedge \widehat y_6 +
 \widehat x_{12} \wedge \widehat y_{11}+ \widehat x_{16} \wedge
 \widehat y_{15}, \\
c_{0}^3= \widehat x_2 \wedge \widehat y_2 +
 \widehat x_5 \wedge \widehat y_5 +
 \widehat x_7 \wedge \widehat y_7 +
 \widehat x_9 \wedge \widehat y_9 +
 \widehat x_{10} \wedge \widehat y_{10}+ \widehat x_{12} \wedge
 \widehat y_{12} +
 \widehat x_{13} \wedge \widehat y_{13} \\
 \hphantom{c_{0}^3=}{}
 + \widehat x_{14} \wedge
 \widehat y_{14}+
 \widehat x_{15} \wedge \widehat y_{15}+ \widehat x_{17} \wedge
 \widehat y_{17}.
\end{gather*}
\end{enumerate}
\end{Lemma}

\subsection{Any rank}
For the following orthogonal Lie algebras and ortho-orthogonal Lie superalgebras, the pattern is clear for any rank, and the answer is given in the following lemmas.

\begin{Lemma}\label{o_pi(2n)}\quad
\begin{enumerate}\itemsep=0pt
\item[$(a1)$] For $\fg=\fo_\Pi(2)$, for a~basis of $H^1(\fg;\fg)$ we can take the following derivations:
\begin{gather*}
1_2 \otimes \widehat 1_2 +E_{1,2}\otimes \widehat E_{1,2},\qquad
E_{1,2}\otimes \widehat E_{1,2}+ E_{2,1} \otimes \widehat E_{2,1} .
\end{gather*}

\item[$(a2)$] For $\fg=\fo_\Pi(2n)$ and $n>2$, for a~basis of $H^1(\fg;\fg)$ we can take the following derivations:
\begin{gather*}
1_{2n} \otimes \widehat E_{i,n+i} \qquad \text{for $i=1,\dots,n$},\nonumber\\
1_{2n} \otimes \widehat E_{n+i,i} \qquad \text{for $i=1,\dots,n$};\nonumber\\
\sum\limits_{1\leq i\leq n} (E_{i,n}+E_{2n,n+i}) \otimes \big(\widehat{E_{i,n}+E_{2n,n+i}}\big) + \sum\limits_{1\leq i\leq n} (E_{n,i}+E_{n+i,2n}) \otimes \big(\widehat{E_{n,i}+E_{n+i,2n}}\big)\nonumber\\
\qquad{} +\sum\limits_{1\leq i\leq n} E_{i,n+i} \otimes \widehat E_{i,n+i} + \sum\limits_{1\leq i\leq n} E_{n+i,i} \otimes \widehat E_{n+i,i} .
\end{gather*}

\item[$(b1)$] For $\fg=\fo_\Pi(2)$, for a~basis of $H^2(\fg)$ we can take the following cocycles:
\[
1_{2}\wedge \widehat E^{1,2}, \qquad 1_{2}\wedge \widehat E^{2,1}.
\]

\item[$(b2)$] For $\fg=\fo_\Pi(2n)$, where $n>1$, for a~basis of $H^2(\fg)$ we can take the following cocycles:
\begin{gather*}
\widehat E^{i,i+n} \wedge \widehat E^{j,j+n}, \qquad 1\leq i <j \leq n,\\
\big(\widehat E^{i,i} + \widehat E^{i+n,i+n}\big)\wedge \widehat E^{i, i+n}\\
\qquad{} +\sum_{1\leq j\leq n, i\not = j} \big(\widehat E^{i,j}+\widehat E^{j+n,i+n}\big) \wedge \big(\widehat E^{j,i+n}+ \widehat E^{i,j+n}\big), \quad i=1,\ldots, n,\\
\widehat E^{i,i+n} \wedge \widehat E^{j+n,j}, \qquad 1\leq i, j \leq n, \quad \text{except for $(i,j)=(n,n)$},\\
\big(\widehat E^{i,i} + \widehat E^{i+n,i+n}\big)\wedge \widehat E^{i+n,i}\\
\qquad{} +\sum_{1\leq j\leq n, i\not = j} \big(\widehat E^{j,i}+\widehat E^{i+n,j+n}\big) \wedge \big(\widehat E^{i+n,j}+ \widehat E^{j+n,i}\big), \qquad i=1,\ldots, n,\\
\widehat E^{i+n,i} \wedge \widehat E^{j+n,j}, \qquad 1\leq i <j \leq n.
\end{gather*}
\end{enumerate}
\end{Lemma}

\section[Orthogonal series $\fo_I(2n)$ and its super versions in characteristic 2]{Orthogonal series $\boldsymbol{\fo_I(2n)}$ and its super versions\\ in characteristic 2}\label{Sortho_I}
In this section, we consider Lie superalgebras none of whose relatives has Cartan matrix.
Let $X_{i,j}:=E_{i,j}+E_{j,i}$.

\begin{Lemma}\label{2.11.1aO} \quad
\begin{enumerate}\itemsep=0pt
\item[$(a)$]
Observe that the Lie algebra $\fo^{(1)}_I(4)$ is not simple; we consider it for completeness of the picture. For $\fg=\fo^{(1)}_I(4)$, for a~basis of $H^1(\fg;\fg)$ we can take
the following derivations of degree $0$:
\begin{itemize}\itemsep=0pt
\item regular cases as for any $n$
\begin{gather*}
X_{1,2}\otimes \widehat X_{1,2} + X_{1,3}\otimes \widehat X_{1,3}
+ X_{1,4}\otimes \widehat X_{1,4}, \\
X_{1,2}\otimes \widehat X_{1,2} + X_{2,3}\otimes \widehat X_{2,3}
+ X_{2,4}\otimes \widehat X_{2,4}, \\
X_{1,3}\otimes \widehat X_{1,3} + X_{2,3}\otimes \widehat X_{2,3} +
X_{3,4}\otimes \widehat X_{3,4};
\end{gather*}
\item exceptional cases
\begin{gather*}X_{1,2}\otimes \widehat X_{3,4} +
X_{1,3}\otimes \widehat X_{2,4}+ X_{2,3}\otimes \widehat X_{1,4}, \\
X_{1,2}\otimes \widehat X_{3,4} + X_{1,4}\otimes \widehat X_{2,3}
+ X_{2,4}\otimes \widehat X_{1,3}, \\
X_{1,3}\otimes \widehat X_{2,4} + X_{1,4}\otimes \widehat X_{2,3}
+ X_{3,4}\otimes \widehat X_{1,2}.
\end{gather*}
\end{itemize}
The space $H^2(\fg)$ is spanned by two cocycles:
\[
\widehat X_{1,2}\wedge \widehat X_{3,4}+ \widehat X_{1,3}\wedge \widehat X_{2,4},\qquad
\widehat X_{1,2}\wedge \widehat X_{3,4}+ \widehat X_{1,4}\wedge \widehat X_{2,3}.
\]

\item[$(b)$] For $\fg=\fo^{(1)}_I(n)$, where $n>4$
or $n=3$, for a~basis of $H^1(\fg;\fg)$ we can take the following
$n-1$ derivations of degree $0$ and the following form:
\begin{gather*}
X_{1,2}\otimes \widehat X_{1,2}
+ X_{1,3}\otimes \widehat X_{1,3} +\dots + X_{1,n}\otimes \widehat X_{1,n}, \\
X_{1,2}\otimes \widehat X_{1,2} + X_{2,3}\otimes \widehat X_{2,3}+\dots +
X_{2,n}\otimes \widehat X_{2,n},\\
X_{1,3}\otimes \widehat X_{1,3} + X_{2,3}\otimes \widehat X_{2,3}+\dots +
X_{3,n}\otimes \widehat X_{3,n},\\
\cdots\cdots\cdots\cdots\cdots\cdots\cdots\cdots\cdots\cdots\cdots\cdots\cdots\cdots\cdots\cdots\cdots\cdots\\
X_{1,n-1}\otimes \widehat X_{1,n-1} + X_{2,n-1}\otimes \widehat X_{2,n-1}+\dots +
X_{n-1,n}\otimes \widehat X_{n-1,n}.
\end{gather*}
We have $H^2(\fg)=0$.
\end{enumerate}
\end{Lemma}

\begin{Statement}[\cite{KrLe}] For $n=3$ or $n\geq 5$, the algebra $\fder
\fo^{(1)}_I(n)$ can be identified with $\fo_I(n)/\fc$
in the sense that for any $D\in\fder\ \fo^{(1)}_I(n)$, there is
$A_D\in\fo_I(n)$ such that $D$ coincides with the restriction of
$\operatorname{ad}_{A_D}$ to $\fo^{(1)}_I(n)$; for a~given $D$, the element $A_D$ is uniquely defined up to adding a~multiple of $1_n$.
\end{Statement}

\begin{proof} In this proof, $i$, $j$, $k$, $l$, $m$ are always indices from $1$ through $n$.

The algebra $\fo^{(1)}_I(n)$ consists of zero-diagonal symmetric $n\times n$ matrices, which means that the elements $\bo{i}{j} := E^{i,j}+E^{j,i}$, where $\{i,j\}$ are all two-element subsets of $\{1,\ldots,n\}$, form a~basis of $\fo^{(1)}_I(n)$. Their commutation relations are (we assume that $i\neq j$ and $k\neq l$):
\[
\big[\bo{i}{j}, \bo{k}{l}\big] = \begin{cases}0&\text{if~} \{k,l\} = \{i,j\} \text{~or~} \{k,l\} \cap \{i,j\} = \varnothing,\\
 \bo{j}{k}&\text{for~} k\neq i, j.
 \end{cases}
\]
Alternatively, we can say that for an arbitrary matrix $M\in \fo^{(1)}_I(n)$ and $i\neq j$, we have
 \begin{alignat*}{3}
&\big[M, \bo{i}{j}\big]_{kl} = 0 \qquad &&\text{if~} \{k,l\} = \{i,j\} \text{~or~} \{k,l\} \cap \{i,j\} = \varnothing,&\\
&\big[M, \bo{i}{j}\big]_{ik} = \big[M, \bo{i}{j}\big]_{ki} = M_{jk}\qquad &&\text{for~} k\neq i,j,& \\
&\big[M, \bo{i}{j}\big]_{kk} = 0\qquad &&\text{for an arbitrary~}k.&
\end{alignat*}
Let $D$ be a~derivation of $\fo^{(1)}_I(n)$. Let us prove that for arbitrary three pairwise distinct indices~$i$,~$j$,~$k$,
\begin{gather}\label{oI-eq1}
\big(D\bo{i}{j}\big)_{ij} + \big(D\bo{i}{k}\big)_{ik} + \big(D\bo{j}{k}\big)_{jk} = 0,
\end{gather}
and for arbitrary four pairwise distinct indices $i$, $j$, $k$, $l$,
\begin{gather}\label{oI-eq2}
\big(D\bo{i}{j}\big)_{kl} = 0,\\
\label{oI-eq3}
\big(D\bo{i}{k}\big)_{il} = \big(D\bo{j}{k}\big)_{jl} = \big(D\bo{i}{l}\big)_{ik}.
\end{gather}

Since $\bo{i}{k} = \big[\bo{i}{j}, \bo{j}{k}\big]$, we have
\begin{align*}
\big(D\bo{i}{k}\big)_{ik} & = \big(D[\bo{i}{j}, \bo{j}{k}]\big)_{ik} = \big[D\bo{i}{j}, \bo{j}{k}\big]_{ik} + \big[\bo{i}{j}, D\bo{j}{k}\big]_{ik} \\
& = \big(D\bo{i}{j}\big)_{ij} + \big(D\bo{k}{j}\big)_{jk},
\end{align*}
which proves (\ref{oI-eq1}).

Note that (\ref{oI-eq2}) and (\ref{oI-eq3}) are vacuously true for $n=3$, since in this case, there exist no four pairwise distinct indices. In case of $n\geq 5$, let $m$ be an index different from all of $i$, $j$, $k$, $l$. Then, $\big[\bo{i}{j}, \bo{l}{m}\big] = 0$, which means that
\[
0 = \big(D[\bo{i}{j}, \bo{l}{m}]\big)_{km} = \big[D\bo{i}{j}, \bo{l}{m}\big]_{km} + \big[\bo{i}{j}, D\bo{l}{m}\big]_{km} = \big(D\bo{i}{j}\big)_{kl} + 0,
\]
which proves (\ref{oI-eq1}). Since $\bo{i}{k} = \big[\bo{i}{j}, \bo{j}{k}\big]$, we have
\[
\big(D\bo{i}{k}\big)_{il} = \big(D[\bo{i}{j}, \bo{j}{k}]\big)_{il} = \big[D\bo{i}{j}, \bo{j}{k}\big]_{il} + \big[\bo{i}{j}, D\bo{j}{k}\big]_{il} = 0 + \big(D\bo{j}{k}\big)_{jl},
\]
which proves the first part of (\ref{oI-eq3}). And since $\bo{i}{k} = [\bo{i}{l}, \bo{k}{l}]$, we have
\[
\big(D\bo{i}{k}\big)_{il} = \big(D[\bo{i}{l}, \bo{k}{l}]\big)_{il} = \big[D\bo{i}{l}, \bo{k}{l}\big]_{il} + \big[\bo{i}{l}, D\bo{k}{l}\big]_{il} = \big(D\bo{i}{l}\big)_{ik} + 0,
\]
which proves the second part of (\ref{oI-eq3}).

Now consider matrix $A_D$ whose entries are as follows:
 \begin{gather*}
 (A_D)_{11} = 0,\\
 (A_D)_{ii} = \big(D\bo{1}{i}\big)_{1i} \qquad \text{for~} i\neq 1,\\
 (A_D)_{ij} = \big(D\bo{k}{i}\big)_{kj} \qquad \text{for~} i\neq j, \text{~where~} k\neq i,j.
 \end{gather*}
Note that (\ref{oI-eq3}) shows that $(A_D)_{ij}$ does not depend on the choice of $k$ and that the matrix is symmetric, i.e., $(A_D)_{ij} = (A_D)_{ji}$.

It follows from (\ref{oI-eq1}), (\ref{oI-eq2}) and (\ref{oI-eq3}) that $D\bo{i}{j} = \big[A_D, \bo{i}{j}\big]$ for an arbitrary $\bo{i}{j}$. (More specifically, using~(\ref{oI-eq1}), (\ref{oI-eq2}) and (\ref{oI-eq3}), one can show that $\big(D\bo{i}{j}\big)_{kl} = \big([A_D, \bo{i}{j}]\big)_{kl}$ for any indices~$k$ and~$l$.)

To have a~clear conscience, one has to consider also several cases depending on whether some of the indices $i$, $j$, $k$, $l$ coincide. We did consider these cases, but skip here these boring details of the proof.

Since the $\bo{i}{j}$ form a~basis of $\fo^{(1)}_I(n)$, it means that $D$ coincides with the restriction of $\operatorname{ad}_{A_D}$ to $\fo^{(1)}_I(n)$. On the other hand, it is easy to check that for any matrix $A\in\fo_I(n)$, the restriction of $\operatorname{ad}_{A}$ to $\fo^{(1)}_I(n)$ is a~derivation of $\fo^{(1)}_I(n)$, and two matrices $A, A'\in\fo_I(n)$ determine the same derivation of $\fo^{(1)}_I(n)$ if and only if $A-A' = c1_n$ for some $c\in \Kee$. This proves our claim.
\end{proof}

\begin{Lemma}\label{2.1.1aO}\quad
\begin{enumerate}\itemsep=0pt
\item[$(a)$] For $\fg=\fo_{\Pi}^{(1)}(2n+1)$,
for a~basis of $H^1(\fg;\fg)$ we can take $2n$ derivations of weight
$\pm 2, \pm 4,\ldots,\allowbreak \pm2n$, and hence $\fder \fg=\fo_{\Pi}(2n+1)$. We have $H^2(\fg)=0$ $($checked for $n\leq6)$.

\item[$(b)$] For $\fg=\fo_{I}(2n)$,
for a~basis of $H^1(\fg;\fg)$ we can take $n$ derivations
$c_i= 1_n \otimes \widehat{E_{i,i}}$, where $i=1, \dots, n$.
The space $H^2(\fg)$ is spanned by $c_{i,j}=\widehat{E_{i,i}}\wedge\widehat{E_{j,j}}$, where $i,j=1,\dots,n$ and $i<j$ $($checked for $n\leq6)$.

\item[$(c)$] For $\fg=\fo\fo_{I\Pi}^{(1)}(1|2k)$, for a~basis of
$H^1(\fg;\fg)$ we can take $2k$ derivations whose weights are $\pm 2,\pm
\ldots, \pm2k$, and hence $\fder\ \fg=\fo\fo_{I\Pi}(1|2k)$.
We have $H^2(\fg)=0$ $($checked for $n\leq6)$.

\item[$(d)$] For $\fg=\fo\fo_{I\Pi}^{(1)}(2n+1|2k)$ and $n\not=0$, for a
basis of $H^1(\fg;\fg)$ we can take $2n$ derivations of weight
$\pm2,\pm4,\ldots,\pm 2n$, and hence $\fder\ \fg=\fo\fo_{I\Pi}(2n+1|2k)$.
The space $H^2(\fg)$ is spanned by $\widehat x_i\wedge \widehat x_i$ such that $x_i$ is odd and $x_i^2=0$, as well as $\widehat y_i\wedge \widehat y_i$
$($checked for $(n,k)$ such that $n\leq3$ and $k\leq4)$.
\end{enumerate}
\end{Lemma}

\section[Simple symmetric Lie (super)algebras without Cartan matrix for $p=2$]{Simple symmetric Lie (super)algebras without Cartan matrix\\ for $\boldsymbol{p=2}$}\label{SymmWithoutCM}

\subsection[Shen's analog $\mathfrak{gs}(2)$ of $\fg(2)$ and the Melikyan algebra for $p=2$]{Shen's analog $\boldsymbol{\mathfrak{gs}(2)}$ of $\boldsymbol{\fg(2)}$ and the Melikyan algebra for $\boldsymbol{p=2}$}\label{2.9} For a
description of the Lie algebra, which we call $\mathfrak{gs}(2)$, in terms
of a~multiplication table, see~\cite{Shen1}; as a~CTS-prolong it is
described in \cite{BGLLS}.

In \cite{Bro}, Brown described analogs of
the Melikyan algebras in characteristic 2 as follows. Recall that
for any $\mathfrak{vect}(n; \un)$-${\mathcal O}(n; \un)$-bimodule $M$ with the
$\mathfrak{vect}(n; \un)$-action $\rho$, we denote by $M_{a \operatorname{div}}$ a~copy of
$M$ with the affine $\mathfrak{vect}(n\un)$-action given by
\begin{gather*}
\rho_{a \operatorname{div}}(D)(\mu)=\rho(D)(\mu)+a \operatorname{div}(D)(\mu),
\end{gather*}
for any $D\in \mathfrak{vect}(n; \un)$, $\mu\in M$ and $a\in \Kee$.
As spaces, and $\Zee/3$-graded Lie algebras, let
\begin{equation}\label{defme2}
L(\un):=\fg_\ev \oplus\fg_{\bar 1}\oplus\fg_{\bar 2} \simeq
\mathfrak{vect}(2; \un)\oplus {\mathcal O}(2; \un)_{\operatorname{div}}\oplus {\mathcal O}(2; \un)\ ,
\end{equation}
 where ${\mathcal O}(2; \un)=\Kee[u_1,
u_2; \un]$ is the space of functions; ${\mathcal O}(2; \un)_{\operatorname{div}}$ is the
space of volume forms with volume element $v:=\operatorname{vol}(u)$ as the
generator of rank 1 module over the algebra ${\mathcal O}(2; \un)$ of
functions. The $\fg_\ev$-action on the $\fg_{\bar i}$ is natural (adjoint, on the space of
volume forms and functions, respectively).

The multiplication in $L(\un)$ is given, for any
$f,g\in{\mathcal O}(2; \un)$, by the following formulas:
\begin{equation*}
 [fv, gv]=0,\qquad [fv, g]=fH_g,\qquad [f, g]:=H_f(g)v,
\end{equation*}
where
\[
H_f=\frac{\partial f}{\partial u_1}\partial_{u_2}+\frac{\partial f}{\partial u_2}\partial_{u_1}.
\]

Define a~$\Zee$-grading of $L(\un)$ by setting
\begin{equation*}
\deg u^{\underline{r}}\partial_i =3|\underline{r}|-3,\qquad \deg
u^{\underline{r}}v =3|\underline{r}|-2,\qquad
\deg u^{\underline{r}} =3|\underline{r}|-4.
\end{equation*}
Now, set $\fme(5;\un):=L(\un)/L(\un)_{-4}$. This algebra is not
simple, because ${\mathcal O}(2; \un)_{\operatorname{div}}$ has a~submodule of codimension
1; but $\fme^{(1)}(5;\un)$ is simple. This is an analog of the \textit{Melikyan algebra}, hence the name.

In \cite{Ei}, the Lie algebra $\fme^{(1)}(5;\One)$, where
$\One=(1,\dots, 1)$, is denoted by $\text{Bro}_2(1,1)$. This
algebra was discovered by Shen Guangyu, see \cite{Shen1}, and should
be denoted somehow to commemorate this, so notation $\mathfrak{gs}(2):=\fme^{(1)}(5;\One)$ for
this analog of $\fg(2)$ (and Guangyu Shen in Western order: first
name, last name) seems most appropriate.

The nonpositive components of $\fme^{(1)}(5;\un)$ and $\mathfrak{gs}(2)$
are the same; they are expressed in terms of vector fields, where
$\partial_i:=\partial_{x_i}$ to distinguish from $\partial_{u_i}$; we use both
representations in terms of $x$ and $u$, whichever is more
convenient:
\begin{equation*}
\renewcommand{\arraystretch}{1.4}
\begin{tabular}{|l|l|} \hline
$\fg_{i}$&the generators\\ \hline\hline $\fg_{-3}$ &
$\partial_{u_1}\longleftrightarrow\partial_{1}$, $\partial_{u_2}\longleftrightarrow\partial_{2}$\\
\hline
$\fg_{-2}$ & $v\longleftrightarrow\partial_{3}$\\
\hline
$\fg_{-1}$ & $u_1\longleftrightarrow (x_3+x_4 x_5)\partial_{2}+\partial_4$, $u_2\longleftrightarrow x_3\partial_{1}+x_4\partial_3+ \partial_5$
\\
\hline
$\fg_{0}$ & $u_1\partial_{1}\longleftrightarrow
x_1\partial_{1}+x_3\partial_3+x_4\partial_4$,\\
& $X_1^+:=u_1\partial_{u_2}:=x_5^{(3)}\partial_{1}+\big(x_1+
x_4x_5^{(2)}\big)\partial_{2}+x_5^{(2)}\partial_3+x_5\partial_4$
\\
&
$X_1^-:=u_2\partial_{u_1}:=\big(x_2+x_4^{(2)} x_5\big)\partial_{1}+x_4^{(3)}\partial_{2}+
+x_4^{(2)}\partial_3+ x_4\partial_5$
\\
& $u_2\partial_{u_2}\longleftrightarrow x_2\partial_{2}+x_3\partial_3+
x_5\partial_5$\\
\hline
\end{tabular}
\end{equation*}
 The highest weight vector in $\fg_{-1}$ is
$X_2^-:=u_1$. Consider the positive part of $\fg=\mathfrak{gs}(2)$. The
lowest weight vector in $\fg_1$ is given by the vector field
\[X_2^+:=x_4^{(3)} x_5\partial_{2}+\big(x_2+ x_4^{(2)}x_5\big)\partial_3+x_4 x_5\partial_5\quad (=u_2v).
\]
So far, the generators and the dimensions of the components look as
their namesakes of $\fg(2)$ for $p>3$; however, the relations are
different: To facilitate comparison with presentations in terms of
Chevalley generators, set $H_i:=[X_i^+, X_i^-]$, i.e.,
\begin{gather*}
H_1=x_1\partial_{1}+x_2\partial_{2}+x_4\partial_4+x_5\partial_5\quad (=u_1\partial_1+u_2\partial_2),\nonumber\\
H_2=x_2\partial_{2}+x_3\partial_3+
x_5\partial_5\quad (= u_2\partial_2).
\end{gather*}
Clearly, $H_1$ is the central element of $\fg_0$. It takes less
space to designate certain invariants of the algebra $\mathfrak{gs}(2)$ in
terms of the indeterminates $u$, as in \eqref{defme2}; e.g., the
cocycles in Lemma~\ref{LemmaShen}; for the grading element of
$\fg_0$ we can take $u_2\partial_{u_2}$, see~\cite{BGL1}.

\begin{Lemma}\label{LemmaShen} For $\fg= \mathfrak{gs}(2)$:
\begin{enumerate}\itemsep=0pt
\item[$1)$] The multiplication table
is as follows:
\begin{alignat*}{4}
&[H_1,X_1^+]= 0,\qquad &&[H_2,X_1^+]= X_1^+,\qquad &&[H_1,H_2]= 0,& \nonumber\\
& [H_1,X_2^+]= X_2^+,\qquad &&[H_2,X_2^+]= 0,\qquad && [X_1^-, X_2^-]= u_2,& \nonumber\\
&[H_1,X_1^-]= 0,\qquad &&[H_2,X_1^-]= X^-_1,\qquad && [X_1^+, X_2^+]= u_1v,& \nonumber\\
&[H_1,X_2^-]= X^-_2,\qquad &&[H_2,X_2^-]= 0,\qquad &&[X_1^\pm,X_2^\mp]= 0.&
\end{alignat*}
The defining relations between the positive and negative root
vectors corresponding to simple roots are as follows:
\begin{alignat*}{3}
& \operatorname{ad}_{X_1^+}^2(X_2^+)=0,\qquad && \operatorname{ad}_{X_2^+}^2(X_1^+)=0,&\nonumber\\
&\operatorname{ad}_{X_2^-}^4(X_1^-)=0,\qquad && \operatorname{ad}_{X_1^-}^2\big(\operatorname{ad}_{X_2^-}^3(X_1^-)\big)=0,& \nonumber\\
& \operatorname{ad}_{X_2^-}\operatorname{ad}_{X_1^-}\big(\operatorname{ad}_{X_2^-}^3(X_1^-)\big)=0,\qquad && \operatorname{ad}_{X_1^-}^2(X_2^-)=0.& 
\end{alignat*}

\item[$2)$] Outer derivations. For a~basis of the space $H^1(\fg;\fg)$
we take the following derivations:
\begin{gather*}
c_{-4}=\partial_{u_1}\otimes \widehat{u_2v}+\partial_{u_2}\otimes \widehat{u_1v}+v\otimes
\widehat{u_1u_2}+u_1\otimes \widehat{u_1u_2\partial_{u_2}}+ u_2\otimes
\widehat{u_1u_2\partial_{u_1}},
\nonumber\\
c_{-2}^1=\partial_{u_2}\otimes \widehat{u_2} +v\otimes \widehat{u_2\partial_{u_1}}+
u_1v\otimes \widehat{u_1u_2\partial_{u_1}}+u_1\partial_{u_2}\otimes \widehat{u_1u_2},
\nonumber\\
c_{-2}^2=\partial_{u_1}\otimes \widehat{u_1} +v\otimes \widehat{u_1\partial_{u_2}}
+u_2v\otimes \widehat{u_1u_2\partial_{u_2}} +u_2\partial_{u_1}\otimes \widehat{u_1u_2},
\nonumber\\
c_{2}^1=u_2\otimes \big(\widehat{\partial_{u_2}}\big) +u_2v\otimes (\widehat{u_1})+
u_2\partial_{u_1}\otimes (\widehat{v})+u_1u_2\otimes \big(\widehat{u_1\partial_{u_2}}\big)
+u_1u_2\partial_{u_1}\otimes
(\widehat{u_1v}),\nonumber\\
c_{2}^2=u_1\otimes \big(\widehat{\partial_{u_1}}\big) +u_1v\otimes (\widehat{u_2})+
u_1\partial_{u_2}\otimes (\widehat{v})+u_1u_2\otimes \big(\widehat{u_2\partial_{u_1}}\big)
+u_1u_2\partial_{u_2}\otimes
(\widehat{u_2v}),\nonumber\\
c_{4}=u_2v\otimes \big(\widehat{\partial_{u_1}}\big)+u_1v\otimes
\big(\widehat{\partial_{u_2}}\big)+u_1u_2\partial_{u_2}\otimes (\widehat{u_1})+u_1u_2\partial_{u_1}\otimes
(\widehat{u_2}).
\end{gather*}

\item[$3)$] Central extensions. For a~basis of the space $H^2(\fg)$ we
take $($the classes of$)$ the following cocycles:
\begin{gather*}
c_{-2}^1= (\widehat{u_2} )\wedge \big(\widehat{u_1u_2\partial_{u_1}}\big) +\big(\widehat{u_2\partial_{u_1}}\big)\wedge (\widehat{u_1u_2} ),\nonumber\\
c_{-2}^2= (\widehat{u_1} )\wedge \big(\widehat{u_1u_2\partial_{u_2}}\big) +\big(\widehat{u_1\partial_{u_2}}\big)\wedge (\widehat{u_1u_2} ),\nonumber\\
c_{4}=\big(\widehat{\partial_{u_1}}\big)\wedge \big(\widehat{u_1}\big)+\big(\widehat{\partial_{u_2}}\big)\wedge (\widehat{u_2} ).
\end{gather*}
\end{enumerate}
\end{Lemma}

\section{Queerifications} Let us recall necessary definitions from
\cite{BLLSq}.

\subsection[Queerification for $p\neq 2$]{Queerification for $\boldsymbol{p\neq 2}$}\label{111}

Let ${\mathcal A}$ be an
associative algebra; let ${\mathcal A}_L$ be the Lie algebra with the same space as ${\mathcal A}$ and
the multiplication being defined by the commutator instead of the dot product (usually denoted by juxtaposition).
The space of the Lie superalgebra $\fq({\mathcal A})$, which we call the
\textit{queerification} of ${\mathcal A}$, is ${\mathcal A}_L\oplus\Pi({\mathcal A})$, where
$\Pi$ is the change of parity functor,
so $\fq({\mathcal A})_\ev={\mathcal A}_L$ and $\fq({\mathcal A})_\od=\Pi({\mathcal A})$, with the multiplication given
by the following expressions for any $x,y\in {\mathcal A}$:
\begin{gather*}
[x,y]=xy-yx\in\fq({\mathcal A})_\ev,\qquad [x,\Pi(y)]=\Pi(xy-yx)\in\fq({\mathcal A})_\od,\nonumber\\
[\Pi(x),\Pi(y)]=
xy+yx\in\fq({\mathcal A})_\ev .
\end{gather*}

The term ``queer'', now universally accepted, comes from the name of the Lie
superalgebra $\fq(n):=\fq({\rm Mat}(n))$, a~``queer'' analog (for reasons given in~\cite{Lsos}) of $\fgl(n)$, where ${\rm Mat}(n)$ is the
associative algebra of $n\times n$ matrices. We express the elements
of the Lie superalgebra $\fg=\fq(n)$ by means of a~pair of matrices
\begin{equation*}
 (A,B)\longleftrightarrow
\begin{pmatrix}A&B\\B&A
\end{pmatrix}
\in \fq(n),
\end{equation*}
where $A,B \in \fgl(n)$.

We will similarly denote the elements of $\fq({\mathcal A})$ by pairs $(A,B)$, where $A,B \in{\mathcal A}$. The brackets between these basis elements are as follows:
\begin{gather}
[(A_1,0),(A_2,0)]=([A_1,A_2],0),\qquad
[(A,0),(0,B)]=(0,[A,B]),\nonumber\\
[(0,B_1), (0,B_2)]=(B_1B_2+B_2B_1,0).\label{.2}
\end{gather}

\subsubsection[The simple part of $\fq(n)$ for $p\ne2$]{The simple part of $\boldsymbol{\fq(n)}$ for $\boldsymbol{p\ne2}$}
Let $\operatorname{qtr} \colon (A,B) \mapsto \operatorname{tr} B$ be called the \textit{queer
trace}; let $\mathfrak{sq}(n):=\fq(n)^{(1)}$ denote the
subsuperalgebra of queer-traceless matrices. The Lie superalgebras $\fq(n)$ and
$\mathfrak{sq}(n)$ are specifically ``super'' analogs of the general Lie
algebra $\fgl(n)$ and its special (traceless) subalgebra $\fsl(n)$; we define their \textit{projectivisations} to be $\fp\fq(n):=\fq(n)/\Kee 1_{2n}$ and
$\mathfrak{psq}(n):=\mathfrak{sq}(n)/\Kee 1_{2n}$.

One might think that if $n=pk$, then $\fp\mathfrak{sq}(n)$ is not simple: its even part is $\mathfrak{pgl}(n)$, its odd part is $\Pi(\fsl(n))$.
However, the two highest weight vectors, one in each homogeneous component, together with the two lowest weight vectors, one in each homogeneous component, generate a~simple Lie superalgebra.

\subsection[Queerification for $p=2$]{Queerification for $\boldsymbol{p=2}$}\label{112} For $\fg=\fq({\mathcal A})$, where ${\mathcal A}$ is an associative algebra, the
multiplication is defined by the expressions~\eqref{.2}, the bracket
of odd elements is the polarization of squaring of the odd elements:
\begin{equation*}
(0,B)^2=\big(B^2, 0\big).
\end{equation*}

If $p=2$, it is possible to queerify, not only associative algebras, but also \textit{any restricted Lie algebra} $\fg$. Namely, set
$\fq(\fg)_\ev=\fg$ and $\fq(\fg)_\od=\Pi(\fg)$; define the
multiplication involving the odd elements as follows, where $x\longmapsto x^{[2]}$ is the $2$-structure:
\begin{equation*}
 [x,\Pi(y)]=\Pi([x,y]),\qquad (\Pi(x))^2=x^{[2]}\qquad\text{for
any~~}x,y\in\fg.
\end{equation*}
Clearly, if $\fg$ is restricted and $\mathfrak{i}\subset\fq(\fg)$ is an ideal,
then $\mathfrak{i}_\ev$ and $\Pi(\mathfrak{i}_\od)$ are ideals in $\fg$. So, if $\fg$ is
restricted and simple, then $\fq(\fg)$ is a~simple Lie superalgebra.
(Note that $\fg$ has to be simple as a~Lie algebra, not just as a
\textit{restricted} Lie algebra, i.e., $\fg$ is not allowed to have
\textit{any} ideals, not only restricted ones.) A generalization of
the queerification is the following procedure producing a~
simple Lie superalgebras for any simple Lie algebras.

\subsubsection{Generalized queerification}\label{Method1} Let the \textit{$1$-step
restricted closure} $\fg^{\langle 1\rangle }$ of the simple Lie algebra $\fg$ be
the minimal subalgebra of the (classically) restricted closure
$\overline{\fg}$ containing $\fg$ and all the elements $x^{[2]}$,
where $x\in\fg$, see~\cite{BLLSq}. To any simple Lie algebra $\fg$ the \textit{generalized queerification} assigns the simple Lie superalgebra%
\begin{equation*}\label{tildeQ}
\tilde \fq(\fg):=\fg^{\langle 1\rangle }\oplus \Pi(\fg)
\end{equation*}
with squaring given by $(\Pi(x))^2=x^{[2]}$ for any $x\in\fg$.
Obviously, for $\fg$ restricted, the generalized queerification coincides with the
queerification: $\tilde \fq(\fg)=\fq(\fg)$.

\begin{Lemma}\label{2.10.1} Let
$\fg=\widetilde {\fq}(\fwk^{(1)}(3;a)/\fc)$ for $a\not =0,1$ $($for
the excluded values of $a$, the algebra $\fwk^{(1)}(3;a)/\fc$ is not
simple$)$.
\begin{enumerate}\itemsep=0pt
\item[$(a)$] For a~basis of $H^1(\fg;\fg)$
we can take the following derivations \emph{(if $a=1$ the cohomology is the same, but if $a=0$ the
cohomology is $6$-dimensional; since the algebra is not simple, we skip the answer)} for the Chevalley basis of $\fwk(3;a)$ considered in formula~\eqref{wk3a}
\begin{gather*}
\underline{c_{0}^1}= h_1\otimes \Pi(\widehat h_1)+ h_2\otimes \Pi(\widehat h_2)+
x_1\otimes \Pi(\widehat x_1) + x_5\otimes \Pi(\widehat x_5)+ x_6\otimes \Pi(\widehat x_6) \\
\hphantom{\underline{c_{0}^1}=}{}
+y_1\otimes \Pi(\widehat y_1)+ y_5\otimes \Pi(\widehat y_5)+ y_6\otimes \Pi(\widehat y_6)+ \Pi(x_2)\otimes \widehat x_2+ \Pi(x_3) \otimes \widehat x_3\\
\hphantom{\underline{c_{0}^1}=}{}
+ \Pi(x_4)\otimes \widehat x_4 +
\Pi(x_7)\otimes \widehat x_7+ \Pi(y_2)\otimes
\widehat y_2+ \Pi(y_3)\otimes \widehat y_3+ \Pi(y_4)\otimes \widehat y_4\\
\hphantom{\underline{c_{0}^1}=}{}
+ \Pi(y_7)\otimes \widehat y_7,\\
c_{0}^2= x_3\otimes \widehat x_3+ x_5\otimes \widehat x_5+ x_6\otimes \widehat x_6+
x_7\otimes \widehat x_7+ y_3\otimes \widehat y_3+ y_5\otimes \widehat y_5+ y_6\otimes \widehat y_6+y_7\otimes \widehat y_7 \\
\hphantom{c_{0}^2=}{}
+ \Pi(h_1)\otimes \Pi(\widehat h_1)+ \Pi(h_2)\otimes \Pi(\widehat h_2)+
\Pi(x_1)\otimes \Pi(\widehat x_1)+ \Pi(x_2)\otimes \Pi(\widehat x_2)\\
\hphantom{c_{0}^2=}{}
+\Pi(x_4)\otimes \Pi(\widehat x_4)+ \Pi(y_1)\otimes \Pi(\widehat y_1)+ \Pi(y_2)\otimes \Pi(\widehat y_2)
+ \Pi(y_4)\otimes \Pi(\widehat y_4),\\
\underline{c_{0}^3}= h_1\otimes \Pi(\widehat h_1)+ h_2\otimes \Pi(\widehat h_2)+ x_1\otimes
\Pi(\widehat x_1) + x_2\otimes \Pi(\widehat x_2)+ x_3\otimes \Pi(\widehat x_3)\\
\hphantom{\underline{c_{0}^3}=}{}
 +
x_4\otimes \Pi(\widehat x_4)+
 x_5\otimes \Pi(\widehat x_5)+ x_6\otimes \Pi(\widehat x_6)+ x_7\otimes
\Pi(\widehat x_7)+ y_1\otimes \Pi(\widehat y_1) \\
\hphantom{\underline{c_{0}^3}=}{}
+ y_2\otimes \Pi(\widehat y_2)+
y_3\otimes \Pi(\widehat y_3)
+ y_4\otimes \Pi(\widehat y_4)+ y_5\otimes \Pi(\widehat y_5)+ y_6\otimes
\Pi(\widehat y_6)\\
\hphantom{\underline{c_{0}^3}=}{}
 + y_7\otimes \Pi(\widehat y_7).
\end{gather*}

\item[$(b)$] The space $H^2(\fg)$ is spanned by the following cocycles:
\begin{gather*}
c_{-8}= \Pi(\widehat x_7)\wedge \Pi(\widehat x_7), \\
c_{-6}= \Pi(\widehat x_6)\wedge \Pi(\widehat x_6), \\
c_{-4}^1= \Pi(\widehat x_4)\wedge \Pi(\widehat x_4), \qquad
c_{-4}^2= \Pi(\widehat x_5)\wedge \Pi(\widehat x_5), \\
c_{-2}^1= \Pi(\widehat x_1)\wedge \Pi(\widehat x_1), \qquad
c_{-2}^2= \Pi(\widehat x_2)\wedge \Pi(\widehat x_2), \qquad
c_{-2}^3= \Pi(\widehat x_3)\wedge \Pi(\widehat x_3), \\
\underline{c_{0}^1}=\big(\alpha^2+\alpha \big) \widehat x_7\wedge
\Pi(\widehat y_7)+\big(\alpha ^2+\alpha \big)
 \widehat y_7\wedge \Pi(\widehat x_7)+\alpha
 \widehat x_5\wedge \Pi(\widehat y_5)+\alpha
\widehat y_5\wedge \Pi(\widehat x_5)\\
\hphantom{\underline{c_{0}^1}=}{}
+\widehat{x_1}\wedge\Pi(\widehat{y_1})+\widehat{x_2}\wedge\Pi(\widehat{y_2})+\widehat{y_1}\wedge\Pi(\widehat{x_1})+\widehat{y_2}\wedge\Pi(\widehat{x_2}),\\
c_{0}^2=\alpha^2\Pi(\widehat{h_2}) \wedge\Pi(\widehat{h_1}) +\big(\alpha ^2+\alpha\big) \Pi(\widehat{x_7})\wedge\Pi(\widehat{y_7})+\alpha
 \Pi(\widehat{x_6})\wedge\Pi(\widehat{y_6})+\alpha \widehat{x_5}\wedge
 \widehat{y_5}\\
\hphantom{c_{0}^2=}{}
 +\Pi(\widehat{x_2})\wedge
 \Pi(\widehat{y_2})+\widehat{x_1}\wedge
 \widehat{y_1}+\widehat{x_3}\wedge
 \widehat{y_3}+\widehat{x_4}\wedge
 \widehat{y_4},\\
c_{0}^3=\alpha^2\Pi(\widehat{h_2 })\wedge \Pi(\widehat{h_1})+\big(\alpha ^2+\alpha\big)
 \widehat{x_7}\wedge \widehat{y_7}+\alpha
 \Pi(\widehat{x_5})\wedge\Pi(\widehat{
 y_5})+\alpha\Pi(\widehat{x_6})\wedge
 \Pi(\widehat{y_6})\\
\hphantom{c_{0}^3=}{}
+ \Pi(\widehat{x_1})\wedge
 \Pi(\widehat{y_1})+\widehat{x_2}\wedge
 \widehat{y_2}+\widehat{x_3}\wedge
 \widehat{y_3}+\widehat{x_4}\wedge
 \widehat{y_4},\\
c_{0}^4=\Pi(\widehat h_1)\wedge \Pi(\widehat h_1).
\end{gather*}
\end{enumerate}
\end{Lemma}

\begin{Lemma}\label{2.10.2} Let $\fg=\mathfrak{q}(\fwk(4;a))$ for
the following Cartan matrix
\[
\begin{pmatrix}
 0 & a~& 1 & 0 \\
a & 0 & 0 & 0 \\
 1 & 0 & 0 & 1 \\
 0 & 0 & 1 & 0 \end{pmatrix}
\]
 and basis of $\fwk(4;a)$
\begin{gather*}
x_1,\ x_2,\ x_3,\ x_4,\
x_5=[x_1, x_2], \ x_6=[x_1, x_3],\ x_7=[x_3, x_4],\
x_8=[x_3, [x_1, x_2]],\nonumber\\
 x_9=[x_4, [x_1, x_3]],\
x_{10}=[[x_1, x_2], [x_1, x_3]],\ x_{11}=[[x_1, x_2], [x_3,x_4]],\nonumber\\
x_{12}=[[x_1, x_2], [x_4, [x_1, x_3]]], \
x_{13}=[[x_3, [x_1, x_2]], [x_4, [x_1, x_3]]], \nonumber\\
x_{14}=[[x_4, [x_1, x_3]], [[x_1, x_2], [x_1, x_3]]],\
x_{15}=[[[x_1, x_2], [x_1, x_3]], [[x_1, x_2], [x_3,x_4]]].
\end{gather*}
\begin{enumerate}\itemsep=0pt
\item[$(a)$] The space $H^1(\fg;\fg)$ is spanned by $($the classes of$)$ the following derivations:
\begin{gather*}
\underline{c_{0}^1}= h_1\otimes \Pi(\widehat y_1)+ h_2\otimes
\Pi(\widehat h_2)+ h_3\otimes \Pi(\widehat h_3)+ h_4\otimes \Pi(\widehat h_4)+ x_1
\otimes \Pi(\widehat x_1)\\
\hphantom{\underline{c_{0}^1}=}{} +
x_2\otimes \Pi(\widehat x_2)+
 x_3\otimes \Pi(\widehat x_3)+ x_4\otimes \Pi(\widehat x_4) + x_5\otimes
\Pi(\widehat x_5)+ x_6\otimes \Pi(\widehat x_6)\\
\hphantom{\underline{c_{0}^1}=}{}
+ x_7\otimes \Pi(\widehat x_7) + x_8\otimes \Pi(\widehat x_8)+
 x_9\otimes \Pi(\widehat x_9)+ x_{10}\otimes \Pi(\widehat x_{10}) +
x_{11}\otimes \Pi(\widehat x_{11}) \\
\hphantom{\underline{c_{0}^1}=}{}
+ x_{12}\otimes \Pi(\widehat x_{12})+
x_{13}\otimes \Pi(\widehat x_{13}) +
 x_{14}\otimes \Pi(\widehat x_{14})+ x_{15}\otimes \Pi(\widehat x_{15})+
y_{1}\otimes \Pi(\widehat y_{1}) \\
\hphantom{\underline{c_{0}^1}=}{}
+ y_{2}\otimes \Pi(\widehat y_2)+ y_{3}\otimes
\Pi(\widehat y_3)+ y_{4}\otimes \Pi(\widehat y_4)+
 y_{5}\otimes \Pi(\widehat y_5)+ y_{6}\otimes \Pi(\widehat y_6)\\
\hphantom{\underline{c_{0}^1}=}{}
 + y_{7}\otimes
\Pi(\widehat y_7)+ y_{8}\otimes \Pi(\widehat y_8)+ y_{9}\otimes \Pi(\widehat y_9)+
y_{10}\otimes \Pi(\widehat y_{10})+
 y_{11}\otimes \Pi(\widehat y_{11})\\
\hphantom{\underline{c_{0}^1}=}{}
 +y_{12}\otimes \Pi(\widehat y_{12})+
y_{13}\otimes \Pi(\widehat y_{13}) + y_{14}\otimes \Pi(\widehat y_{14})+
y_{15}\otimes
\Pi(\widehat y_{15}),\\
c_{0}^2= x_4\otimes \widehat x_4+ x_7\otimes \widehat x_7 + x_9\otimes \widehat x_9+
x_{11}\otimes \widehat x_{11}+ x_{12}\otimes \widehat x_{12} + x_{13}\otimes \widehat
x_{13}+
x_{14}\otimes \widehat x_{14}\\
\hphantom{c_{0}^2=}{}
 + x_{15}\otimes \widehat x_{15} +
 y_4\otimes \widehat y_4+ y_7\otimes \widehat y_7+ y_9\otimes \widehat y_9 +
y_{11}\otimes \widehat y_{11}+ y_{12}\otimes \widehat y_{12}+ y_{13}\otimes \widehat
y_{13} \\
\hphantom{c_{0}^2=}{}
+
y_{14}\otimes \widehat y_{14}+ y_{15}\otimes \widehat y_{15}+
 \Pi(h_1)\otimes \Pi(\widehat h_1) + \Pi(h_2)\otimes \Pi(\widehat h_2)+
\Pi(h_3)\otimes \Pi(\widehat h_3) \\
\hphantom{c_{0}^2=}{}
+ \Pi(h_4)\otimes \Pi(\widehat h_4) +
\Pi(x_1)\otimes \Pi(\widehat x_1)+
 \Pi(x_2)\otimes \Pi(\widehat x_2)+ \Pi(x_3)\otimes \Pi(\widehat x_3)\\
\hphantom{c_{0}^2=}{}
 +
\Pi(x_5)\otimes \Pi(\widehat x_5)+ \Pi(x_6)\otimes \Pi(\widehat x_6)+
\Pi(x_8)\otimes \Pi(\widehat x_8)
+ \Pi(x_{10})\otimes \Pi(\widehat x_{10})\\
\hphantom{c_{0}^2=}{}
 + \Pi(y_1)\otimes \Pi(\widehat y_{1})
+ \Pi(y_2)\otimes \Pi(\widehat y_2) + \Pi(y_3)\otimes \Pi(\widehat y_3)+
\Pi(y_5)\otimes \Pi(\widehat y_5) \\
\hphantom{c_{0}^2=}{}
+
 \Pi(y_6)\otimes \Pi(\widehat y_6) + \Pi(y_8)\otimes \Pi(\widehat y_8)+
\Pi(y_{10})\otimes \Pi(\widehat y_{10}).
 \end{gather*}

\item[$(b)$] The space $H^2(\fg)$ is spanned by the following cocycles:
\begin{gather*}
c_{-16}= \Pi(\widehat x_{15}) \wedge \Pi(\widehat x_{15}), \\
c_{-14}=\Pi(\widehat x_{14}) \wedge \Pi(\widehat x_{14}), \\
c_{-12}= \Pi(\widehat x_{13}) \wedge \Pi(\widehat x_{13}), \\
c_{-10}= \Pi(\widehat x_{12}) \wedge \Pi(\widehat x_{12}), \\
c_{-8}^1=\Pi(\widehat x_{10}) \wedge \Pi(\widehat x_{10}), \qquad
c_{-8}^2=\Pi(\widehat x_{11}) \wedge \Pi(\widehat x_{11}), \\
c_{-6}^1=\Pi(\widehat x_{8}) \wedge \Pi(\widehat x_{8}), \qquad
c_{-6}^2=\Pi(\widehat x_{9}) \wedge \Pi(\widehat x_{9}), \\
c_{-4}^1=\Pi(\widehat x_{5}) \wedge \Pi(\widehat x_{5}), \qquad
c_{-4}^2=\Pi(\widehat x_{6}) \wedge \Pi(\widehat x_{6}), \qquad
c_{-4}^3=\Pi(\widehat x_{7}) \wedge \Pi(\widehat x_{7}), \\
c_{-4}^2= \Pi(\widehat x_{2}) \wedge \Pi(\widehat x_{2}), \qquad
c_{-4}^3= \Pi(\widehat x_{3}) \wedge \Pi(\widehat x_{3}), \qquad
c_{-4}^4=\Pi(\widehat x_{4}) \wedge \Pi(\widehat x_{4}), \\
c_{-2}^1=\Pi(\widehat x_{1}) \wedge \Pi(\widehat x_{1}), \\
c_{0}^1= \widehat x_2 \wedge \widehat y_2+
 \widehat x_3 \wedge \widehat y_3+
 \alpha \widehat x_5 \wedge \widehat y_5+ \widehat x_6 \wedge \widehat y_6+ \widehat x_7 \wedge \widehat y_7+ \widehat x_9 \wedge \widehat y_9\\
\hphantom{c_{0}^1=}{}
 + \big(\alpha ^2+\alpha\big)
 \widehat x_{13}\wedge \widehat y_{13}+\big(\alpha ^3+\alpha ^2\big) \widehat x_{14} \wedge
 \widehat y_{14}+ \alpha \Pi (\widehat
 h_1)\wedge \Pi (\widehat h_2)+
\Pi (\widehat h_1)\wedge \Pi (\widehat h_3)\\
\hphantom{c_{0}^1=}{}
 + \Pi(\widehat h_3) \wedge \Pi(\widehat
 h_4)+ \Pi(\widehat x_1)\wedge
 \Pi(\widehat y_1)+ \Pi(\widehat
 x_4)\wedge \Pi(\widehat y_4)+
 \alpha \Pi(\widehat x_8) \wedge \Pi(\widehat
 y_8)\\
\hphantom{c_{0}^1=}{}
 + \big(\alpha ^2+\alpha\big)
 \Pi(\widehat x_{10}) \wedge \Pi (\widehat{y_{10}})
 + \alpha \Pi(\widehat x_{11})\wedge \Pi(\widehat
 y_{11})+ \big(\alpha ^2+\alpha\big) \Pi(\widehat x_{12})\wedge \Pi(\widehat
 y_{12})\\
\hphantom{c_{0}^1=}{}
 + \big(\alpha ^4+\alpha
 ^3\big) \Pi(\widehat x_{15})\wedge \Pi(\widehat
 y_{15}), \\
c_{0}^2= \widehat x_2 \wedge \widehat y_2+
 \widehat x_4 \wedge \widehat y_4 +
 \alpha \widehat x_5 \wedge \widehat y_5+\widehat x_7\wedge \widehat y_7+ \alpha \widehat x_8 \wedge
 \widehat y_8+ \widehat x_9\wedge
 \widehat y_9+ \alpha ^2+\alpha
 \widehat x_{10}\wedge \widehat y_{10}\\
\hphantom{c_{0}^2= }{}
 +\big( \alpha ^4+\alpha ^3\big)
 \widehat x_{15} \wedge \widehat y_{15}+
 \alpha \Pi(\widehat h_1)\wedge \Pi(\widehat
 h_2)+ \Pi(\widehat h_1)\wedge
 \Pi(\widehat h_3)+ \Pi(\widehat
 h_3) \wedge \Pi(\widehat h_4)\\
\hphantom{c_{0}^2= }{}
 +
 \Pi(\widehat x_1) \wedge \Pi(\widehat y_1)
+ \Pi(\widehat x_3)\wedge \Pi(\widehat
 y_3)+ \Pi(\widehat x_6)\wedge
 \Pi(\widehat y_6)+ \alpha \Pi(\widehat x_{11}) \wedge \Pi(\widehat y_{11})\\
\hphantom{c_{0}^2= }{}
 + \big(\alpha ^2+\alpha\big) \Pi(\widehat
 x_{12}) \wedge \Pi(\widehat y_{12})
 + \big(\alpha ^2+\alpha\big) \Pi(\widehat x_{13})\wedge
 \Pi(\widehat y_{13})\\
\hphantom{c_{0}^2= }{}
 + \big(\alpha
 ^3+\alpha ^2\big) \Pi(\widehat x_{14})\wedge \Pi(\widehat
 y_{14}),\\
c_{0}^3= \widehat x_1\wedge \widehat y_1+
 \widehat x_3 \wedge \widehat y_3 +
 \alpha \widehat x_5\wedge \widehat y_5+ \widehat x_7 \wedge \widehat y_7+ \big(\alpha ^2+\alpha\big)
 \widehat x_{10}\wedge \widehat y_{10}\\
 \hphantom{c_{0}^3=}{}
 +\big(\alpha ^2+\alpha\big) \widehat x_{12}\wedge
 \widehat y_{12}
 + \big(\alpha ^3+\alpha
 ^2\big) \widehat x_{14} \wedge \widehat y_{14}+ \big(\alpha ^4+\alpha ^3\big)
 \widehat x_{15} \wedge \widehat y_{15}\\
 \hphantom{c_{0}^3=}{}
 +
 \Pi(\widehat x_1) \wedge \Pi(\widehat y_1)+ \Pi(\widehat x_3) \wedge \Pi(\widehat
 y_3)
 + \alpha \Pi(\widehat
 x_5) \wedge \Pi(\widehat y_5) +
 \Pi(\widehat x_7) \wedge \Pi(\widehat y_7)\\
 \hphantom{c_{0}^3=}{}
 + \big(\alpha ^2+\alpha\big) \Pi(\widehat
 x_{10})\wedge \Pi(\widehat y_{10})
 +
 \big(\alpha ^2+\alpha\big) \Pi(\widehat x_{12})\wedge
 \Pi(\widehat y_{12}) \\
 \hphantom{c_{0}^3=}{}
 + \big(\alpha
 ^3+\alpha ^2 \big) \Pi(\widehat x_{14}) \wedge \Pi(\widehat y_{14})+ \big(\alpha ^4+\alpha
 ^3\big) \Pi(\widehat x_{15}) \wedge \Pi(\widehat
 y_{15}),\\
c_{0}^4= \Pi(\widehat h_1)\wedge \Pi(\widehat h_1).
\end{gather*}
\end{enumerate}
\end{Lemma}

\section[Periplectic Lie superalgebras and their desuperizations for $p=2$]{Periplectic Lie superalgebras and their desuperizations\\ for $\boldsymbol{p=2}$}\label{Speripl}

For completeness, we have to consider periplectic Lie superalgebras for $p\neq 2$ although they are ``symmetric'' only for $p=2$.

In this section, we set
\begin{gather}
\operatorname{diag}(A,B):=\begin{pmatrix}A&0\\
0&B\end{pmatrix},\qquad \operatorname{antidiag}(A,B):=\begin{pmatrix}0&A\\
B&0\end{pmatrix}, \nonumber\\
E_{i,j}^{(n)}\text{~~or just $E_{i,j}$ designates the $n\times
n$ matrix unit with a~1 in the $(i,j)$th slot},\nonumber\\
d_i=\operatorname{antidiag}(E_{i,i}^{(n)}, 0_{n}),\qquad
d_{1,1}^{(n)}=\operatorname{diag}(E_{1,1}^{(n)}, E_{1,1}^{(n)}),\qquad
D_n=\operatorname{diag}(0_n, 1_n).\label{notat}
\end{gather}

\subsection[The Lie superalgebras preserving odd symmetric non-degenerate bilinear forms for $p>2$]{The Lie superalgebras preserving odd symmetric non-degenerate\\ bilinear forms for $\boldsymbol{p>2}$}\label{section7.1}

The normal shape of the periplectic (odd
\textit{symmetric} non-degenerate) bilinear form has the Gram matrix
$J_{2n}=\operatorname{antidiag}(1_n,-1_n)$, that of \textit{anti}symmetric form is
$\Pi_{n|n}=\operatorname{antidiag}(1_n,1_n)$ (sic: signs are correct, see~\cite{Lsos}).

With the following types of Lie superalgebras several series of
``odd" analogs of Lie algebras of Hamiltonian and contact vector
fields and an exceptional vectorial Lie superalgebra are associated
for $p\neq 2$.

The case 3) in the following list becomes much more involved for $p=2$, see~\cite{BGLLS}:
\begin{enumerate}
\item[1)] $\mathfrak{pe}(n):=\mathfrak{aut} (J_{2n})\simeq\mathfrak{aut} (\Pi_{n|n})$: if the matrix
of the bilinear form ${\mathcal B}$ is $\Pi_{2n}$ (resp. $\Pi_{n|n}$), then the Lie (super)algebra
$\mathfrak{aut} ({\mathcal B})$ preserving the bilinear form ${\mathcal B}$ consists of the (super)matrices of the form
\begin{gather}
\mathfrak{aut} (J_{2n})=\left\{X=\begin{pmatrix} A & B\\C &
-A^t
\end{pmatrix} \bigg| \, B=B^{\rm t},\ C=-C^{\rm t} \right\}, \nonumber\\
\mathfrak{aut} (\Pi_{n|n})= \left\{X=\begin{pmatrix} A & B\\C &
-A^t
\end{pmatrix} \bigg|\, B=-B^{\rm t},\ C=C^{\rm t} \right\}.\label{matrgen}
\end{gather}

\item[2)] $\mathfrak{spe}(n):=\mathfrak{pe}(n)\cap \mathfrak{sle}(n|n)$, so $\mathfrak{pe}(n)=\mathfrak{spe}(n)\ltimes
\Kee D$, where $D=\operatorname{diag}(1_n,-1_n)$.
\item[3)] $\mathfrak{spe}_{a,b}(n):=\mathfrak{spe}(n)\ltimes \Kee (aD+b1_{2n})$, where
$a,b\in\Kee$.
\item[4)] $\fp(\mathfrak{spe})(pn):=\mathfrak{spe}(pn)/\fc$.
\end{enumerate}

Over $\Cee$, only for $n=4$ does $\mathfrak{spe}(n)$ have a~nontrivial
central extension \eqref{eq3.1.1.0}, discovered by A.~Sergeev; with
its spinor representation a~simple exceptional infinite-dimensional
vectorial Lie superalgebra is associated, see \cite{LSh}.

For $p=2$, there are 8 analogs of $\mathfrak{pe}(n)$ and 8 analogs of $\mathfrak{spe}(n)$, see
Table~\eqref{Table}, and lots of nontrivial central extensions of these analogs are
found in this section.

\begin{Lemma}\label{3.1.1} \quad
\begin{enumerate}\itemsep=0pt
\item[$(a)$] For any $p>2$ and $n>2$, we have
$H^2(\mathfrak{pe}(n))=0$.
For any $p>3$, we have $H^2(\mathfrak{spe}(n))=0$ for $n>4$ and $n=3$ whereas
for a~basis of $H^2(\mathfrak{spe}(3))$ \textit{only if $p=3$} we can take the following three
cocycles
\begin{gather*}
c_0^1=\big(\widehat{E_{2,3}-E_{6,5}}\big)\wedge \big(\widehat{E_{1,3}-E_{6,4}}\big)+\big(\widehat{E_{1,5}-
E_{2,4}}\big)\wedge \big(\widehat{E_{6,3}}\big),\nonumber\\
c_{0}^2= \big(\widehat{E_{2,3}-E_{5,6}}\big)\wedge \big(\widehat{E_{1,2}-E_{5,4}}\big)+
\big(\widehat{E_{1,6}-E_{3,4}}\big)\wedge \big(\widehat{E_{5,2}}\big),\nonumber\\
c_{0}^3= 2 \big(\widehat{E_{2,1}-E_{4,5}}\big)\wedge
\big(\widehat{E_{3,1}-E_{4,6}}\big)+\big(\widehat{E_{2,6}-E_{3,5}}\big)\wedge \big(\widehat{E_{4,1}}\big).
\end{gather*}

For every $p\neq 2$, for a~basis of $H^2(\mathfrak{spe}(4))$ we can take the
cocycle \begin{gather}
c_{-2}=\big(\widehat{E_{1,6}-E_{2,5}}\big)\wedge \big(\widehat{E_{3,8}-E_{4,7}}\big)-
\big(\widehat{E_{1,7}-E_{3,5}}\big)\wedge \big(\widehat{E_{2,8}-E_{4,6}}\big)\nonumber\\
\hphantom{c_{-2}=}{}
+\big(\widehat{E_{1,8}-E_{4,5}}\big)
\wedge \big(\widehat{E_{2,7}-E_{3,8}}\big).\label{eq3.1.1.0}
 \end{gather}

\item[$(b)$] For any $p>2$, we have
$H^1(\mathfrak{pe}(n);\mathfrak{pe}(n))=0$. For a~basis of $H^1(\mathfrak{spe}(n);\mathfrak{spe}(n))$ we
can take the outer derivation $D=\operatorname{diag}(1_n, -1_n)$ for any $n>2$.
\end{enumerate}
\end{Lemma}

\subsection[The Lie superalgebras preserving odd symmetric non-degenerate bilinear forms for $p=2$ and their desuperizations: $\mathfrak{pe}_{\rm gen}(n)$ and $\fo_{\rm gen}(n)$]{The Lie superalgebras preserving odd symmetric non-degenerate bilinear\\ forms for $\boldsymbol{p=2}$ and their desuperizations: $\boldsymbol{\mathfrak{pe}_{\rm gen}(n)}$ and $\boldsymbol{\fo_{\rm gen}(n)}$}

For $p=2$, the \textit{general periplectic Lie superalgebra} accrues additional elements as compared with the cases $p\neq 2$, see \eqref{matrgen}:
\begin{equation*}
\mathfrak{pe}_{\rm gen}(n):=\mathfrak{aut} (\Pi_{n|n})= \left\{X=\begin{pmatrix} A & B\\C &
A^t
\end{pmatrix} \bigg|\, B=B^{\rm t},\ C=C^{\rm t}\right\}.
\end{equation*}
Its desuperization will be called the \textit{general orthogonal Lie algebra}:
\begin{equation*}
\fo_{\rm gen}(2n):=\mathfrak{aut} (\Pi_{2n})= \left\{X=\begin{pmatrix} A & B\\C &
A^t
\end{pmatrix} \bigg| B=B^{\rm t},\ C=C^{\rm t}\right\}.
\end{equation*}

Let
\begin{equation*}
{\rm ZD}:=\text{the space of symmetric matrices with zeros on
their main diagonals.}
\end{equation*} The derived Lie (super)algebra
$\mathfrak{aut}^{(1)}({\mathcal B})$ consists of the (super)matrices of the form
\eqref{matrgen}, where $B,C\in {\rm ZD}$. In other words, these Lie
(super)algebras resemble the orthogonal Lie algebras. On these Lie
(super)algebras $\mathfrak{aut}^{(1)}({\mathcal B})$ the following (super)trace---we call it the \emph{half-trace}---is defined:
\begin{equation*}
\htr \colon \ \begin{pmatrix} A & B\\C & A^{\rm t}
\end{pmatrix}\tto \operatorname{tr} A.
\end{equation*}
The Lie sub(super)algebra of $\mathfrak{aut}^{(1)}({\mathcal B})$ consisting of
half-traceless supermatrices is isomorphic to $\mathfrak{aut}^{(2)}({\mathcal B})$.

Define the Lie superalgebra
\begin{align*}
\widetilde{\mathfrak{aut}} ({\mathcal B}):={}& \{\text{(super)matrices of the form \eqref{matrgen}, where $C\in {\rm ZD}$}\},\nonumber\\
 \simeq {}& \{\text{(super)matrices of the form~\eqref{matrgen}, where $B\in {\rm ZD}$}\}.
\end{align*}

Clearly, $\mathfrak{aut} ^{(1)}({\mathcal B})\subset\widetilde{\mathfrak{aut}}
({\mathcal B})\subset\mathfrak{aut} ({\mathcal B})$. Introduce a~new notation $\fo\fp$:
\begin{equation*}
\widetilde{\mathfrak{aut}} ({\mathcal B}):=\begin{cases}\fo\fp(2n)&\text{if ${\mathcal B}$ is even},\\
\mathfrak{pe}(n)&\text{if ${\mathcal B}$ is odd.}
\end{cases}
\end{equation*}

Set $\deg
A^{i, i+1}=-\deg A^{i+1, i}=1$, and setting the degree of the lowest
(resp.\ highest) weight vector in the $\fgl(n)$-module of matrices~$B$ (resp.\ $C$) in \eqref{matrgen} be equal to 1 (resp.\ $-1$). We have
the following cases to consider:
\begin{equation}\label{Table}
\renewcommand{\arraystretch}{1.2}
\begin{tabular}{|l|l|}
\hline 1& $\fo_{\rm gen}(2n)$ and $\mathfrak{pe}_{\rm gen}(n)$ (in \eqref{matrgen}
both $B$
and $C$ symmetric);\\
1a&$\fp(\fo)_{\rm gen}(2n):=\fo_{\rm gen}(2n)/\Kee\ 1_{2n}$ and
$\fp(\mathfrak{pe})_{\rm gen}(n):=\mathfrak{pe}_{\rm gen}(n)/\Kee\ 1_{2n}$ \\
\hline
2&$\fo\fp(2n)$ and $\mathfrak{pe}(n)$ (in \eqref{matrgen} both $B$ and $C$
are symmetric, $C\in {\rm ZD}$);\\
2a&$\fp(\fo\fp)(2n):=\fo_{\rm gen}(2n)/\Kee\ 1_{2n}$ and
$\fp(\mathfrak{pe})(n):=\mathfrak{pe}_{\rm gen}(n)/\Kee\ 1_{2n}$;\\
2b&$\fo^{(1)}(2n)$ (in \eqref{matrgen} both $B\in {\rm ZD}$ and $C\in
{\rm ZD}$);\\
2c&$\fp(\fo^{(1)})(2n):=\fo^{(1)}(2n)/\Kee\ 1_{2n}$ \\
\hline
3&$\mathfrak{sop}(2n)$ and $\mathfrak{spe}(n)$ (in \eqref{matrgen} both $B$ and
$C$ are symmetric,\\
& $C\in {\rm ZD}$, and $\operatorname{tr} A=0$);\\
3a&$\fp(\mathfrak{sop})(2n):=\mathfrak{sop}(2n)/\Kee\ 1_{2n}$ and
$\fp(\mathfrak{spe})(n):=\mathfrak{spe}(n)/\Kee\ 1_{2n}$ \\
\hline
4&$\fo^{(2)}(2n)$ (in \eqref{matrgen} both $B\in {\rm ZD}$ and $C\in
{\rm ZD}$, and $\operatorname{tr} A=0$);\\
4a&$\fp(\fo^{(2)})(2n):=\fo^{(2)}(2n)/\Kee\ 1_{2n}$ \\
\hline
\end{tabular}
\end{equation}

\begin{Lemma}\label{3.2.1} \quad
\begin{enumerate}\itemsep=0pt
\item[$(a)$]Let $\fg=\fo_{\rm gen}(2n)$.
For a~basis of $H^2(\fg)$ we
can take the following cocycles $($recall convention~\eqref{conv}$)$
\begin{gather}
c_{-2}^{i,j}= \big(\widehat{E_{n+i,i}}\big)\wedge \big(\widehat{E_{n+j,j}}\big) \qquad \text{for~~} 1\leq i< j \leq n,\nonumber\\
c_{-1}^j=\big(\widehat{E_{j,j}+E_{n+j,n+j}}\big)\wedge \big(\widehat{E_{n+j,j}}\big)\nonumber\\
\hphantom{c_{-1}^j=}{}
+
\sum\limits_{i\not = j} \big(\widehat{E_{i,j}+E_{j+n,i+n}}\big)\wedge \big(\widehat{E_{j+n,i}+E_{i+n,j}}\big)
\qquad \text{for~~}1\leq j \leq n,\nonumber\\
c_0^{i,j}=\big(\widehat{E_{i,i+n}}\big)\wedge \big(\widehat{E_{n+j,j}}\big)\qquad \text{for~~}1\leq i,j\leq n \ \text{~~except~}i=j=n. \label{pogen}
\end{gather}

\item[$(b)$] For a~basis of $H^1(\fg;\fg)$ we can take the
derivation $D_n$, see formula~\eqref{notat}, for any $n>2$ and also the derivations $($recall convention~\eqref{conv}$)$
\begin{gather*} 
c_{-1}^j= \sum\limits_{1\leq i\leq n}
E_{i,i} \otimes \big(\widehat{E_{n+j,j}}\big) \qquad \text{for~~}1\leq j \leq n.
\end{gather*}
\end{enumerate}
\end{Lemma}

\begin{Lemma}\label{3.2.2} \quad
\begin{enumerate}\itemsep=0pt
\item[$(a)$]Let $\fg=\mathfrak{pe}_{\rm gen}(2n)$.
 For a~basis of $H^2(\fg)$
we can take the cocycles \eqref{pogen} $($except for the cocycles
$c_{-1}^j$ and symmetric to them $c_1^j$ that become trivial$)$, and
also the following cocycles $($recall convention~\eqref{conv}$)$
\begin{gather}
c_{-2}^i= \big(\widehat{E_{n+i,i}}\big)^{\wedge^2} \qquad \text{for~~}1\leq i \leq n,\nonumber\\
\widetilde c_{-2}^{i,j}= \big(\widehat{E_{n+i,j} + E_{n+j,i}}\big)^{\wedge^2} \qquad \text{for~~}1\leq i< j \leq n.\label{pegen}
\end{gather}
\item[$(b)$] For a~basis of $H^1(\fg;\fg)$ we can take the
derivation $D_n$ for any $n>2$ and also the derivations $($recall convention~\eqref{conv}$)$
\begin{gather*} 
c_{-1}^j= \sum\limits_{1\leq i\leq n}
E_{i,i} \otimes \big(\widehat{E_{n+j,j}}\big) \qquad \text{for~~}1\leq j \leq n.
\end{gather*}
\end{enumerate}
\end{Lemma}

\begin{Lemma}\label{3.2.3}\quad
\begin{enumerate}\itemsep=0pt
\item[$(a)$] Let $\fg=\fo_{\rm gen}(2n)/\Kee
1_{2n}$. For a~basis of $H^2(\fg)$ we can take the following cocycles $($recall convention~\eqref{conv}$)$
\begin{alignat*}{3}
& \deg=-2\colon \quad && \text{cocycles as in \eqref{pogen}},&\nonumber\\
& \deg=-1\colon \quad &&\text{none},& \nonumber\\
& \deg=0\colon \quad && \text{cocycles as in \eqref{pogen} and the cocycle } \big(\widehat{E_{n,2n}}\big)\wedge \big(\widehat{E_{2n,n}}\big).& 
\end{alignat*}

\item[$(b)$] For a~basis of $H^1(\fg;\fg)$ we can take the derivation $D_n$ for any $n>2$.
\end{enumerate}
\end{Lemma}

\begin{Lemma}\label{3.2.4} \quad
\begin{enumerate}\itemsep=0pt
\item[$(a)$] Let $\fg=\mathfrak{pe}_{\rm gen}(2n)/\Kee
1_{2n}$. For a~basis of
$H^2(\fg)$ we can take the following
cocycles $($recall convention~\eqref{conv}$)$
\begin{alignat*}{3}
& \deg=-2\colon \quad && \text{cocycles as in \eqref{pegen} and \eqref{pogen}},& \nonumber\\
& \deg=-1\colon \quad && \text{none},&\nonumber\\
& \deg=0\colon \quad && \text{cocycles as in \eqref{pogen} and the cocycle},& \nonumber\\
& \text{for $n\neq 4$}\colon \quad &&
\big(\widehat{E_{n,2n}}\big)\wedge \big(\widehat{E_{2n,n}}\big),&\nonumber\\
& \text{for $n= 4$}\colon \quad && \sum \limits_{i=2, 3,
4}\big(\widehat{E_{1,i}+E_{i+4,5}}\big)\wedge \big(\widehat{E_{i,1}+E_{5,i+4}}\big)&\nonumber\\
&&& {}+ \sum \limits_{i=6, 7, 8}\big(\widehat{E_{1,i}+E_{i-4,5}}\big)\wedge \big(\widehat{E_{i,1}+E_{5,i-4}}\big).& 
\end{alignat*}
\item[$(b)$] For a~basis of $H^1(\fg;\fg)$ we
can take the derivation $D_n$ for any $n>2$.
\end{enumerate}
\end{Lemma}
\newpage
\begin{Lemma}\label{3.2.28}\quad
\begin{enumerate}\itemsep=0pt
\item[$(a)$] Let $\fg=\fo^{(2)}(2n)$. For a~basis of
$H^1(\fg;\fg)$ we can take the following
derivations $($recall convention~\eqref{conv}$)$
\begin{alignat*}{3}
& \deg=-1\colon \quad && \text{multiplication by the matrix $d_i$, see Section~{\rm \ref{section7.1}},}& \nonumber\\
&&& \text{for $i=1,\ldots,n$ and any $n>2$}, &\nonumber \\
& \deg=0\colon && \text{multiplication by the matrix } d=\operatorname{diag}(0_n,1_n) \text{ for any $n>2$;}&\nonumber \\
&&& \text{if $n$ is even, there is one more derivation:}&\nonumber\\
&&& \text{multiplication by the matrix $d_{1,1}$}.& 
\end{alignat*}

\item[$(b)$] For a~basis of $H^2(\fg)$
we can take the following cocycles $($recall convention~\eqref{conv}$)$
\begin{gather*}
\deg=-2\colon \quad \text{none for $n\neq 4$}.
\end{gather*}
 For $n=4$, we have (the class of) the cocycle
 \begin{gather*}
c_{-2}= \big(\widehat{E_{7,4}+E_{8,3}}\big)\wedge \big(\widehat{E_{5,2}+E_{6,1}}\big)+ \big(\widehat{E_{6,3}+E_{7,2}}\big)\wedge \big(\widehat{E_{5,4}+E_{8,1}}\big)\\
\hphantom{c_{-2}=}{}
+\big(\widehat{E_{6,4}+
E_{8,2}}\big)\wedge \big(\widehat{E_{5,3}+E_{7,1}}\big),\\
\deg=-1\colon \quad c^j_{-1}= \sum \limits_{i\not =j}
\big(\widehat{E_{i,j}+
E_{j+n,i+n}}\big)\wedge \big(\widehat{E_{j+n,i}+E_{i+n,j}}\big) \qquad \text{for~~}1 \leq j \leq n.
\end{gather*}
For $n=4$, there are also the following $4$ cocycles:
\begin{gather*}
c_{-1}^1= \big(\widehat{E_{1,2}+E_{6,5}}\big)\wedge \big(\widehat{E_{7,4}+E_{8,3}}\big)+ \big(\widehat{E_{1,3}+E_{7,5}}\big)\wedge \big(\widehat{E_{6,4}
+E_{8,2}}\big)\\
\hphantom{c_{-1}^1=}{}
+\big(\widehat{E_{1,4}+
E_{8,5}}\big)\wedge \big(\widehat{E_{6,3}+E_{7,2}}\big), \\
c_{-1}^2= \big(\widehat{E_{2,3}+E_{7,6}}\big)\wedge \big(\widehat{E_{5,4}+E_{8,1}}\big)+ \big(\widehat{E_{2,4}+E_{8,6}}\big)\wedge \big(\widehat{E_{5,3}+E_{7,1}}\big)\\
\hphantom{c_{-1}^2=}{}
+\big(\widehat{E_{2,1}+
E_{5,6}}\big)\wedge \big(\widehat{E_{7,4}+E_{8,3}}\big),\\
c_{-1}^3= \big(\widehat{E_{3,4}+E_{8,7}}\big)\wedge \big(\widehat{E_{5,2}+E_{6,1}}\big)+ \big(\widehat{E_{3,2}+E_{6,7}}\big)\wedge \big(\widehat{E_{5,4}+E_{8,1}}\big)\\
\hphantom{c_{-1}^3=}{}
+\big(\widehat{E_{3,1}+
E_{5,7}}\big)\wedge \big(\widehat{E_{6,4}+E_{8,2}}\big),\\
c_{-1}^4= \big(\widehat{E_{4,3}+E_{7,8}}\big)\wedge \big(\widehat{E_{5,2}+E_{6,1}}\big)+ \big(\widehat{E_{4,2}+E_{6,8}}\big)\wedge \big(\widehat{E_{5,3}+E_{7,1}}\big)\\
\hphantom{c_{-1}^4=}{}
+\big(\widehat{E_{4,1}+
E_{5,8}}\big)\wedge \big(\widehat{E_{6,3}+E_{7,2}}\big),\\
\deg=0\colon\\
 c_0= \sum \limits_{1\leq i\leq n} \big(\widehat{E_{1,i}+
E_{n+i,n+1}}\big)\wedge \big(\widehat{E_{i,1}+E_{n+1,n+i}}\big)\\
\hphantom{c_0=}{}
+
 \sum \limits_{2\leq i <j\leq n} \big(\widehat{E_{i,n+j}+
E_{j,n+i}}\big)\wedge \big(\widehat{E_{n+j,i}+E_{n+i,j}}\big).
\end{gather*}
For $n=4$, there are also the following $6$ cocycles:
\begin{gather*}
c_0^1= \big(\widehat{E_{3,2}-E_{6,7}}\big)\wedge
 \big(\widehat{E_{4,1}-E_{5,8}}\big) +\big(\widehat{E_{4,2}-E_{6,8}}\big)\wedge
 \big(\widehat{E_{3,1}-E_{5,7}}\big)\nonumber\\
\hphantom{c_0^1=}{}
 +\big(\widehat{E_{3,8}-E_{4,7}}\big)\wedge
 \big(\widehat{E_{5,2}+E_{6,1}}\big),\nonumber\\
c_0^2= \big(\widehat{E_{3,4}-E_{8,7}}\big)\wedge
 \big(\widehat{E_{2,1}-E_{5,6}}\big) +\big(\widehat{E_{2,4}-E_{8,6}}\big)\wedge
 \big(\widehat{E_{3,1}-E_{5,7}}\big)\nonumber\\
\hphantom{c_0^2=}{}
 +\big(\widehat{E_{2,7}-E_{3,6}}\big)\wedge
 \big(\widehat{E_{5,4}+E_{8,1}}\big),\nonumber\\
c_0^3= \big(\widehat{E_{4,3}-E_{7,8}}\big)\wedge
 \big(\widehat{E_{2,1}-E_{5,6}}\big)+
\big(\widehat{E_{2,3}-E_{7,6}}\big)\wedge
 \big(\widehat{E_{4,1}-E_{5,8}}\big)\nonumber\\
\hphantom{c_0^3=}{}
 +\big(\widehat{E_{2,8}-E_{4,6}}\big)\wedge
 \big(\widehat{E_{5,3}+E_{7,1}}\big),\nonumber\\
c_0^4= \big(\widehat{E_{2,3}-E_{7,6}}\big)\wedge
 \big(\widehat{E_{1,4}-E_{8,5}}\big) +\big(\widehat{E_{2,4}-E_{8,6}}\big)\wedge
 \big(\widehat{E_{1,3}-E_{7,5}}\big)\nonumber\\
\hphantom{c_0^4=}{}
 +\big(\widehat{E_{1,6}-E_{2,5}}\big)\wedge
 \big(\widehat{E_{7,4}+E_{8,3}}\big),\nonumber\\
c_0^5= \big(\widehat{E_{3,4}-E_{8,7}}\big)\wedge
 \big(\widehat{E_{1,2}-E_{6,5}}\big)+
\big(\widehat{E_{3,2}-E_{6,7}}\big)\wedge
 \big(\widehat{E_{1,4}-E_{8,5}}\big)\nonumber\\
\hphantom{c_0^5=}{}
 + \big(\widehat{E_{1,7}-E_{3,5}}\big)\wedge
 \big(\widehat{E_{6,4}+E_{8,2}}\big),\nonumber\\
c_0^6= \big(\widehat{E_{4,3}-E_{7,8}}\big)\wedge
 \big(\widehat{E_{1,2}-E_{6,5}}\big)+\big(\widehat{E_{4,2}-E_{6,8}}\big)\wedge
 \big(\widehat{E_{1,3}-E_{7,5}}\big)\nonumber\\
\hphantom{c_0^6=}{}
 + \big(\widehat{E_{1,8}-E_{4,5}}\big)\wedge
 \big(\widehat{E_{6,3}+E_{7,2}}\big).
\end{gather*}
\end{enumerate}
\end{Lemma}

\begin{Lemma}\label{3.2.29}\quad
\begin{enumerate}\itemsep=0pt
\item[$(a)$] Let $\fg=\fo^{(2)}(4n)/\Kee 1_{4n}$. For a~basis of
$H^2(\fg)$ we can take the following cocycles $($recall convention~\eqref{conv}$)$
\begin{alignat*}{3}
& \deg=-2\colon \quad && \text{as in the case of $\fo^{(2)}(4n)$},&\nonumber\\
& \deg=-1\colon \quad && \text{as in the case of $\fo^{(2)}(4n)$},&\nonumber\\
& \deg=0\colon \quad && \text{as in the case of $\fo^{(2)}(4n)$ and one more cocycle},& \nonumber\\
&&& \widetilde c_0= \sum \limits_{1\leq i\leq n} \big(\widehat{E_{1,i}+E_{i+n,n+1}}\big)\wedge
 \big(\widehat{E_{i,1}+E_{n+1,n+i}}\big).& 
\end{alignat*}

\item[$(b)$] For a~basis of $H^1(\fg;\fg)$ we can take the following derivations $($recall convention~\eqref{conv}$)$
\begin{alignat*}{3}
& \deg=-2\colon \quad &&
 E_{3,8} \otimes \big(\widehat{E_{5,2}+E_{6,1} }\big)-E_{4,7} \otimes \big(\widehat{E_{5,2}+
 E_{6,1}}\big)+E_{2,7} \otimes \big(\widehat{E_{5,4}+E_{8,1}}\big)\\
 &&&{} -E_{3,6} \otimes \big(\widehat{E_{5,4}+E_{8,1}} \big)
 +
 E_{2,8} \otimes \big(\widehat{E_{5,3}+
E_{7,1} }\big)-E_{4,6}\otimes \big(\widehat{E_{5,3}+E_{7,1}}\big)&\\
&&&{} +E_{1,6}\otimes
\big(\widehat{E_{7,4}+E_{8,3}}\big)-E_{2,5}\otimes \big(\widehat{E_{7,4}+E_{8,3}}\big)
 +E_{1,7}\otimes \big(\widehat{E_{6,4}+E_{8,2}}\big)&\\
 &&&{} -E_{3,5} \otimes \big(\widehat{E_{6,4}+E_{8,2}} \big)+E_{1,8}\otimes \big(\widehat{E_{6,3}+
E_{7,2}}\big)-E_{4,5}\otimes \big(\widehat{E_{6,3}+E_{7,2}}\big),& \\
& \deg=-1 \colon \quad && \text{multiplication by the matrix
$d_i$, see Section~{\rm \ref{section7.1}}},&\\
&&& \text{for $i=1,\ldots,n$ for any $n>2$}.&
\end{alignat*}
 For $n=2$ , there are also four derivations given by
\begin{gather*}
c_{-1}^1= E_{2,1}\otimes \big(\widehat{E_{7,4}+E_{8,3}}\big) +E_{5,6}\otimes \big(\widehat{E_{7,4}+E_{8,3} }\big)+
 E_{3,1}\otimes \big(\widehat{E_{6,4}+E_{8,2}}\big) \nonumber\\
\hphantom{c_{-1}^1=}{}
 +E_{5,7}\otimes \big(\widehat{E_{6,4}+E_{8,2} }\big) +
 E_{4,1}\otimes \big(\widehat{E_{6,3}+E_{7,2}}\big)
 +E_{5,8}\otimes \big(\widehat{E_{6,3}+E_{7,2} }\big)\nonumber\\
\hphantom{c_{-1}^1=}{}
 +
 E_{3,8}\otimes \big(\widehat{E_{1,2}-E_{6,5}}\big)
 +E_{4,7}\otimes \big(\widehat{E_{1,2}-E_{6,5} }\big) +
 E_{2,7}\otimes \big(\widehat{E_{1,4}-E_{8,5}}\big) \nonumber\\
\hphantom{c_{-1}^1=}{}
 +E_{3,6}\otimes \big(\widehat{E_{1,4}-E_{8,5}}\big)+
 E_{2,8}\otimes \big(\widehat{E_{1,3}-E_{7,5}}\big)
 +E_{4,6}\otimes \big(\widehat{E_{1,3}-E_{7,5} }\big),\nonumber\\
c_{-1}^2= E_{1,4}\otimes \big(\widehat{E_{6,3}+E_{7,2} }\big)
 +E_{8,5}\otimes \big(\widehat{E_{6,3}+E_{7,2} }\big) +
 E_{2,4}\otimes \big(\widehat{E_{5,3}+E_{7,1} }\big)\nonumber\\
\hphantom{c_{-1}^2=}{}
 +E_{8,6}\otimes \big(\widehat{E_{5,3}+E_{7,1}}\big) +
 E_{3,4}\otimes \big(\widehat{E_{5,2}+E_{6,1} }\big)
 +E_{8,7}\otimes \big(\widehat{E_{5,2}+E_{6,1} }\big) \nonumber\\
\hphantom{c_{-1}^2=}{}
 +
 E_{2,7}\otimes \big(\widehat{E_{4,1}-E_{5,8} }\big)+
E_{3,6}\otimes \big(\widehat{E_{4,1}-E_{5,8} }\big) +
 E_{1,6}\otimes \big(\widehat{E_{4,3}-E_{7,8} }\big)
 \nonumber\\
\hphantom{c_{-1}^2=}{}
 +E_{2,5}\otimes \big(\widehat{E_{4,3}-E_{7,8} }\big) +
 E_{1,7}\otimes \big(\widehat{E_{4,2}-E_{6,8}}\big)
+E_{3,5}\otimes \big(\widehat{E_{4,2}-E_{6,8}} \big), \nonumber\\
c_{-1}^3= E_{1,2}\otimes \big(\widehat{E_{7,4}+E_{8,3} }\big)
 +E_{6,5}\otimes \big(\widehat{E_{7,4}+E_{8,3} }\big) +
 E_{3,2}\otimes \big(\widehat{E_{5,4}+E_{8,1} }\big)\nonumber\\
\hphantom{c_{-1}^3=}{}
 +E_{6,7}\otimes \big(\widehat{E_{5,4}+E_{8,1} }\big) +
 E_{4,2}\otimes \big(\widehat{E_{5,3}+E_{7,1} }\big)+
E_{6,8}\otimes \big(\widehat{E_{5,3}+E_{7,1} }\big)
\nonumber\\
\hphantom{c_{-1}^3=}{}
+
 E_{3,8}\otimes \big(\widehat{E_{2,1}-E_{5,6}} \big)+
E_{4,7}\otimes \big(\widehat{E_{2,1}-E_{5,6} }\big) +
 E_{1,7}\otimes \big(\widehat{E_{2,4}-E_{8,6} }\big)
\nonumber\\
\hphantom{c_{-1}^3=}{}
 +E_{3,5}\otimes \big(\widehat{E_{2,4}-E_{8,6} }\big) +
 E_{1,8}\otimes \big(\widehat{E_{2,3}-E_{7,6} }\big)+
E_{4,5}\otimes \big(\widehat{E_{2,3}-E_{7,6} }\big),\nonumber\\
c_{-1}^4= E_{1,3}\otimes \big(\widehat{E_{6,4}+E_{8,2} }\big)+E_{7,5}\otimes \big(\widehat{E_{6,4}+E_{8,2} }\big) +
 E_{2,3}\otimes \big(\widehat{E_{5,4}+E_{8,1} }\big)\nonumber\\
\hphantom{c_{-1}^4=}{}
 +E_{7,6}\otimes \big(\widehat{E_{5,4}+E_{8,1}} \big) +
 E_{4,3}\otimes \big(\widehat{E_{5,2}+E_{6,1} }\big)+
 E_{7,8}\otimes \big(\widehat{E_{5,2}+E_{6,1}}\big)
 \nonumber\\
\hphantom{c_{-1}^4=}{}
 +
 E_{2,8}\otimes \big(\widehat{E_{3,1}-E_{5,7} }\big)
+E_{4,6}\otimes \big(\widehat{E_{3,1}-E_{5,7}} \big) +
 E_{1,6}\otimes \big(\widehat{E_{3,4}-E_{8,7}} \big)
 \nonumber\\
\hphantom{c_{-1}^4=}{}
 +E_{2,5}\otimes \big(\widehat{E_{3,4}-E_{8,7}} \big)+
 E_{1,8}\otimes \big(\widehat{E_{3,2}-E_{6,7} }\big)
+E_{4,5}\otimes \big(\widehat{E_{3,2}-E_{6,7} }\big).
\end{gather*}
For $\deg=0$, there is multiplication by the matrix $D_{2n}$ for any $n>1$,
and the multiplication by the matrix $d_{1,1}^{(2n)}$, see Section~{\rm \ref{section7.1}}. Besides, if $n=2$, there are
six derivations given by the following cocycles:
\begin{gather*}
c_0^1= E_{3,1}\otimes \big(\widehat{E_{2,4}-E_{8,6} }\big) +E_{5,7}\otimes \big(\widehat{E_{2,4}-E_{8,6} }\big) +
 E_{4,1}\otimes \big(\widehat{E_{2,3}-E_{7,6} }\big)\nonumber\\
\hphantom{c_0^1=}{}
 +E_{5,8}\otimes \big(\widehat{E_{2,3}-E_{7,6}} \big) +
 E_{3,2}\otimes \big(\widehat{E_{1,4}-E_{8,5} }\big)
 +E_{6,7}\otimes \big(\widehat{E_{1,4}-E_{8,5} }\big)
\nonumber\\
\hphantom{c_0^1=}{}
 +
 E_{4,2}\otimes \big(\widehat{E_{1,3}-E_{7,5}}\big)
+E_{6,8}\otimes \big(\widehat{E_{1,3}-E_{7,5}} \big)
+E_{3,8}\otimes \big(\widehat{E_{1,6}-E_{2,5}}\big)
\nonumber\\
\hphantom{c_0^1=}{}
+E_{4,7}\otimes \big(\widehat{E_{1,6}-E_{2,5} }\big) +
 +E_{5,2}\otimes \big(\widehat{E_{7,4}+E_{8,3} }\big) +
 E_{6,1}\otimes \big(\widehat{E_{7,4}+E_{8,3} }\big),\nonumber\\
c_0^2= E_{2,4}\otimes \big(\widehat{E_{1,3}-E_{7,5}} \big)
 +E_{8,6}\otimes \big(\widehat{E_{1,3}-E_{7,5} }\big) +
 E_{3,4}\otimes \big(\widehat{E_{1,2}-E_{6,5} }\big)\nonumber\\
 \hphantom{c_0^2=}{}
 +
 E_{8,7}\otimes \big(\widehat{E_{1,2}-E_{6,5}} \big) +
 E_{2,1}\otimes \big(\widehat{E_{4,3}-E_{7,8} }\big)
 +E_{5,6}\otimes \big(\widehat{E_{4,3}-E_{7,8} }\big)
 \nonumber\\
 \hphantom{c_0^2=}{}
 +
 E_{3,1}\otimes \big(\widehat{E_{4,2}-E_{6,8} }\big)
+E_{5,7}\otimes \big(\widehat{E_{4,2}-E_{6,8} }\big) +
 E_{2,7}\otimes \big(\widehat{E_{1,8}-E_{4,5}} \big)
 \nonumber\\
 \hphantom{c_0^2=}{}
 +E_{3,6}\otimes \big(\widehat{E_{1,8}-E_{4,5} }\big) +
 E_{5,4}\otimes \big(\widehat{E_{6,3}+E_{7,2} }\big) +
 E_{8,1}\otimes \big(\widehat{E_{6,3}+E_{7,2} }\big),\nonumber\\
c_0^3= E_{2,3}\otimes \big(\widehat{E_{1,4}-E_{8,5}} \big)
 +E_{7,6}\otimes \big(\widehat{E_{1,4}-E_{8,5} }\big) +
 E_{2,1}\otimes \big(\widehat{E_{3,4}-E_{8,7}} \big)\nonumber\\
 \hphantom{c_0^3=}{}
 +E_{5,6}\otimes \big(\widehat{E_{3,4}-E_{8,7} }\big) +
 E_{4,1}\otimes \big(\widehat{E_{3,2}-E_{6,7} }\big)
 +E_{5,8}\otimes \big(\widehat{E_{3,2}-E_{6,7} }\big)
 \nonumber\\
 \hphantom{c_0^3=}{}
 +
 E_{4,3}\otimes \big(\widehat{E_{1,2}-E_{6,5} }\big)
 +E_{7,8}\otimes \big(\widehat{E_{1,2}-E_{6,5} }\big) +
 E_{2,8}\otimes \big(\widehat{E_{1,7}-E_{3,5} }\big)
 \nonumber\\
 \hphantom{c_0^3=}{}
 +E_{4,6}\otimes \big(\widehat{E_{1,7}-E_{3,5} }\big) +
 E_{5,3}\otimes \big(\widehat{E_{6,4}+E_{8,2}}\big) +
 E_{7,1}\otimes \big(\widehat{E_{6,4}+E_{8,2}} \big),\nonumber\\
c_0^4= (E_{1,3}\otimes \big(\widehat{E_{4,2}-E_{6,8} }\big)
 +E_{7,5}\otimes \big(\widehat{E_{4,2}-E_{6,8} }\big) +
 E_{1,4}\otimes \big(\widehat{E_{3,2}-E_{6,7}} \big)\nonumber\\
 \hphantom{c_0^4=}{}
 +E_{8,5}\otimes \big(\widehat{E_{3,2}-E_{6,7} }\big) +
 E_{2,3}\otimes \big(\widehat{E_{4,1}-E_{5,8} }\big)
 +E_{7,6}\otimes \big(\widehat{E_{4,1}-E_{5,8} }\big)
 \nonumber\\
 \hphantom{c_0^4=}{}
 +
 E_{2,4}\otimes \big(\widehat{E_{3,1}-E_{5,7} }\big)
+E_{8,6}\otimes \big(\widehat{E_{3,1}-E_{5,7} }\big) +
 E_{1,6}\otimes \big(\widehat{E_{3,8}-E_{4,7} }\big)
 \nonumber\\
 \hphantom{c_0^4=}{}
 +E_{2,5}\otimes \big(\widehat{E_{3,8}-E_{4,7} }\big) +
 E_{7,4}\otimes \big(\widehat{E_{5,2}+E_{6,1} }\big) +
 E_{8,3}\otimes \big(\widehat{E_{5,2}+E_{6,1} }\big),\nonumber\\
c_0^5= E_{1,2}\otimes \big(\widehat{E_{4,3}-E_{7,8} }\big)
 +E_{6,5}\otimes \big(\widehat{E_{4,3}-E_{7,8} }\big)+
 E_{1,4}\otimes \big(\widehat{E_{2,3}-E_{7,6} }\big)\nonumber\\
 \hphantom{c_0^5=}{}
 +E_{8,5}\otimes \big(\widehat{E_{2,3}-E_{7,6} }\big) +
 E_{3,4}\otimes \big(\widehat{E_{2,1}-E_{5,6}} \big)
 +E_{8,7}\otimes \big(\widehat{E_{2,1}-E_{5,6} }\big)
 \nonumber\\
 \hphantom{c_0^5=}{}
 +
 E_{3,2}\otimes \big(\widehat{E_{4,1}-E_{5,8} }\big)
 +E_{6,7}\otimes \big(\widehat{E_{4,1}-E_{5,8} }\big) +
 E_{1,7}\otimes \big(\widehat{E_{2,8}-E_{4,6} }\big)
 \nonumber\\
 \hphantom{c_0^5=}{}
 +E_{3,5}\otimes \big(\widehat{E_{2,8}-E_{4,6} }\big) +
 E_{6,4}\otimes \big(\widehat{E_{5,3}+E_{7,1} }\big) +
 E_{8,2}\otimes \big(\widehat{E_{5,3}+E_{7,1} }\big),\nonumber\\
c_0^6= E_{1,2}\otimes \big(\widehat{E_{3,4}-E_{8,7} }\big)
 +E_{6,5}\otimes \big(\widehat{E_{3,4}-E_{8,7} }\big) +
 E_{1,3}\otimes \big(\widehat{E_{2,4}-E_{8,6} }\big)\nonumber\\
\hphantom{c_0^6=}{}
 +E_{7,5}\otimes \big(\widehat{E_{2,4}-E_{8,6} }\big) +
 E_{4,2}\otimes \big(\widehat{E_{3,1}-E_{5,7} }\big)
 +E_{6,8}\otimes \big(\widehat{E_{3,1}-E_{5,7} }\big)
\nonumber\\
\hphantom{c_0^6=}{}
 +
 E_{4,3}\otimes \big(\widehat{E_{2,1}-E_{5,6} }\big)
 +E_{7,8}\otimes \big(\widehat{E_{2,1}-E_{5,6} }\big) +
 E_{1,8}\otimes \big(\widehat{E_{2,7}-E_{3,6} }\big)
\nonumber\\
\hphantom{c_0^6=}{}
 +E_{4,5}\otimes \big(\widehat{E_{2,7}-E_{3,6}} \big) +
 E_{6,3}\otimes \big(\widehat{E_{5,4}+E_{8,1} }\big) +
 E_{7,2}\otimes \big(\widehat{E_{5,4}+E_{8,1} }\big).
\end{gather*}
\end{enumerate}
\end{Lemma}

\begin{Lemma}\label{3.2.7}\quad
\begin{enumerate}\itemsep=0pt
\item[$(a)$] Let $\fg=\fo \fp(2n)$. For a~basis of $H^2(\fg)$ we
can take the cocycles
\begin{alignat*}{3}
& \deg=-2, 0\colon \quad &&\text{none},&\nonumber\\
& \deg=1 \colon \quad && c_{-1}^{j}= \sum \limits_{i\not = j}
\big(\widehat{E_{i,j}+E_{j+n,i+n}}\big)\wedge \big(\widehat{E_{j+n,i}+E_{i+n,j}}\big)\qquad \text{for~~}1\leq j \leq n,&\nonumber\\
&&&c_1^j= \sum \limits_{1\leq i\leq n}
\big(\widehat{E_{i,i}+E_{i+n,i+n}}\big)\wedge \big(\widehat{E_{j,n+j}}\big)\qquad \text{for~~}1\leq j \leq
n, &\nonumber\\
&&& \widetilde c_1^j= \sum \limits_{i\not=j}
\big(\widehat{E_{j,i}+E_{i+n,j+n}}\big)\wedge \big(\widehat{E_{i,n+j}+E_{j,n+i}}\big)&\nonumber\\
&&& \hphantom{\widetilde c_1^j=}{} +
 \big(\widehat{E_{j,j}+E_{n+j,n+j}}\big)\wedge \big(\widehat{E_{j,n+j}}\big) \qquad \text{for~~}1\leq j \leq n,& \nonumber\\
& \deg=2 \colon \quad && \text{cocycles as in \eqref{pogen}}& 
\end{alignat*}

\item[$(b)$] For a~basis of $H^1(\fg;\fg)$ we can take the following derivations
\begin{alignat*}{3}
&\deg=-1\colon \quad && \text{none},&\nonumber\\
& \deg=0\colon \quad && \text{multiplication by the matrix $D_n$, see Section~{\rm \ref{section7.1}}},&\nonumber\\
&&& \text{for any $n>2$ and also the cocycle},&\nonumber\\
&&& c_0^2= \sum \limits_{1\leq i\leq 2n}
E_{i,i}\otimes \big(\big(\widehat{E_{1,1}+E_{n+1,n+1}}\big)+\cdots + \big(\widehat{E_{n,n}+E_{2n,2n}}\big) \big),\nonumber\\
& \deg=1\colon \quad && c_1^j= \sum \limits_{1\leq i\leq 2n} E_{i,i}\big(\widehat{E_{j,j+n}}\big), \qquad \text{for~~} 1\leq j \leq n\text{ and derivations},\nonumber\\
&\widetilde c_1^j\colon \quad && \text{multiplication by the matrix $d_j^{\rm t}$, see Section~{\rm \ref{section7.1}}},&\nonumber\\
& && \text{for any $j=1,\dots, n$ and $n>2$.}& 
\end{alignat*}
\end{enumerate}
\end{Lemma}

\begin{Lemma}\label{3.2.8}\quad
\begin{enumerate}\itemsep=0pt
\item[$(a)$] Let $\fg=\mathfrak{pe}(2n)$. For a~basis of $H^2(\fg)$ we can
take $($the classes of$)$ the cocycles
\begin{alignat*}{3}
& \deg=-2 \colon \quad &&
\widetilde c_{-2}^{i,j}= \big(\big(\widehat{E_{n+i,j} + E_{n+j,i}}\big)\big)^{\wedge^2} \qquad \text{for~~}1\leq i< j \leq n, &\nonumber\\
& \deg=-1,0 \colon \quad && \text{none},&\nonumber\\
& \deg=1 \colon \quad && c_1^j=\mathop{\sum}\limits_{1\leq i\leq n}
\big(\widehat{E_{i,i}+E_{i+n,i+n}}\big)\wedge \big(\widehat{E_{j,n+j}}\big) \qquad \text{for~~}1\leq j \leq
n,& \nonumber\\
&&& \widetilde c_1^j= \sum \limits_{i\not=j}^n \big(\widehat{E_{j,i}+
E_{i+n,j+n}}\big)\wedge \big(\widehat{E_{i,n+j}+E_{j,n+i}}\big)& \nonumber\\
&&& \hphantom{\widetilde c_1^j=}{} + \big(\widehat{E_{j,j}+
E_{n+j,n+j}}\big)\wedge \big(\widehat{E_{j,n+j}}\big)\qquad \text{for~~}1\leq j \leq n,&\nonumber\\
& \deg=2 \colon \quad &&
c_{2}=\text{cocycles as in \eqref{pogen} and \eqref{pegen}}.&
\end{alignat*}

\item[$(b)$] For a~basis of $H^1(\fg;\fg)$ we can take the following derivations
\begin{alignat*}{3}
& \deg=-1\colon \quad && \text{none},& \nonumber\\
& \deg=0\colon \quad && \text{multiplication by the matrix $D_n$, see Section~{\rm \ref{section7.1}}}, &\nonumber\\
&&& \text{for any $n>2$ and also the cocycle},&\nonumber\\
&&&c_0^2=
 \sum \limits_{1\leq i\leq 2n}
E_{i,i}\otimes \big(\big(\widehat{E_{1,1}+E_{n+1,n+1}}\big)+\cdots + \big(\widehat{E_{n,n}+E_{2n,2n}}\big) \big), &\nonumber\\
& \deg=1 \colon \quad && \widetilde c_1^j\colon \text{multiplication by the matrix
$d_j^{\rm t}$, see Section~{\rm \ref{section7.1}},} &\nonumber\\
&&& \text{for any $j=1,\dots, n$ and $n>2$. }& 
\end{alignat*}
\end{enumerate}
\end{Lemma}
\newpage

\begin{Lemma}\label{3.2.11}\quad
\begin{enumerate}\itemsep=0pt
\item[$(a)$] Let $\fg=\fp(\fo\fp)(2n)$. For a~basis of $H^2(\fg)$
we can take $($the classes of$)$ the following cocycles
\begin{alignat*}{3}
& \deg=-2, -1\colon \, && \text{none},&\nonumber\\
&\deg=0\colon \quad && \text{none for $n\not =4$. For $n=4$, the cocycle is}&\nonumber\\
&&& c_0= \big(\widehat{E_{1,2}+E_{6,5}}\big) \wedge \big(\widehat{E_{2,1}-E_{5,6}}\big)+
\big(\widehat{E_{2,3}+E_{7,6}}\big) \wedge \big(\widehat{E_{3,2}-E_{6,7}}\big) &\nonumber\\
&&& \hphantom{c_0=}{} + \big(\widehat{E_{2,4}+E_{8,6}}\big) \wedge \big(\widehat{E_{4,2}-E_{6,8}}\big)+
\big(\widehat{E_{2,7}+E_{3,6}}\big) \wedge \big(\widehat{E_{7,2}-E_{6,3}}\big)&\nonumber\\
&&& \hphantom{c_0=}{}+\big(\widehat{E_{2,8}+E_{4,6}}\big) \wedge \big(\widehat{E_{6,4}-E_{8,2}}\big)+
\big(\widehat{E_{1,6}+E_{2,5}}\big) \wedge \big(\widehat{E_{5,2}-E_{6,1}}\big), &\nonumber\\
&&& c_{1}^j= \sum \limits_{2\leq i\leq n} \big(\widehat{E_{i,i}+E_{i+n,i+n}}\big)\wedge \big(\widehat{E_{j,j+n}}\big),
\text{where~} 1\leq j \leq n, &\nonumber\\
& \deg=2\colon \quad && \text{cocycles as in \eqref{pogen}}.& 
\end{alignat*}

\item[$(b)$] For a~basis of $H^1(\fg;\fg)$ we can take the following derivations
\begin{alignat*}{3}
& \deg=-1\colon\ && \text{none},&\\ 
& \deg=0\colon \quad && \text{multiplication by the matrix $D_n$, see Section~{\rm \ref{section7.1}}, for any $n>2$},&\nonumber\\
& \deg=1\colon \quad && \widetilde c_1^j\colon \text{multiplication by the matrix
$d_j^{\rm t}$ for any $j = 1,{\dots},n$ and $n>2$}. &\nonumber
\end{alignat*}
\end{enumerate}
\end{Lemma}

\begin{Lemma}\label{3.2.13} \quad
\begin{enumerate}\itemsep=0pt
\item[$(a)$]
Let $\fg=\fp(\mathfrak{pe})(2n)$. For a~basis of $H^2(\fg)$
we can take $($the classes of$)$ the following cocycles:
\begin{alignat*}{3}
& \deg=-2\colon \quad && c_{-2}^{i,j}= \big(\widehat{E_{n+i,j} +E_{n+j,i}}\big)^{\wedge^2},\qquad
\text{where~~}1\leq i< j \leq n, &\nonumber \\
& \deg=-1\colon \quad && \text{none},&\nonumber\\
& \deg=0\colon \quad && \text{none for $n\not =4$. For $n=4$, the cocycle is} \nonumber\\
&&& c_0= \big(\widehat{E_{1,2}+E_{6,5}}\big) \wedge \big(\widehat{E_{2,1}+E_{5,6}}\big)+
\big(\widehat{E_{2,3}+E_{7,6}}\big) \wedge \big(\widehat{E_{3,2}+E_{6,7}}\big)&\nonumber\\
&&& \hphantom{c_0=}{} + \big(\widehat{E_{2,4}+E_{8,6}}\big)
\wedge \big(\widehat{E_{4,2}+E_{6,8}}\big)+\big(\widehat{E_{2,7}+E_{3,6}}\big) \wedge \big(\widehat{E_{7,2}+E_{6,3}}\big) & \nonumber\\
&&&
 \hphantom{c_0=}{}
+\big(\widehat{E_{2,8}+E_{4,6}}\big)
\wedge \big(\widehat{E_{6,4}+E_{8,2}}\big)+
\big(\widehat{E_{1,6}+E_{2,5}}\big) \wedge \big(\widehat{E_{5,2}+E_{6,1}}\big), & \nonumber \\
& \deg=1\colon \quad && \sum \limits_{2\leq i\leq n}
\big(\widehat{E_{i,i}+E_{i+n,i+n}}\big)\wedge \big(\widehat{E_{j,j+n}}\big), \qquad \text{where~~} 1\leq j \leq n, &\nonumber \\
& \deg=2\colon \quad && \text{cocycles as in \eqref{pogen} and \eqref{pegen}.}& 
\end{alignat*}

\item[$(b)$] For a~basis of $H^1(\fg;\fg)$ we can take the
following derivations
\begin{alignat*}{3}
& \deg=-1 \colon \quad && \text{none},&\nonumber\\
& \deg=0\colon \quad && \text{multiplication by the matrix $D_n$, see Section~{\rm \ref{section7.1}}, for any $n>2$},&\nonumber\\ 
& \deg=1\colon \quad && \widetilde c_1^j\colon \text{multiplication by the matrix $d_j^{\rm t}$, see Section~{\rm \ref{section7.1}},}& \nonumber\\
&&& \text{for any $j=1,\dots, n$ and $n>2$.} & 
\end{alignat*}
\end{enumerate}
\end{Lemma}

\begin{Lemma}\label{3.2.15}\quad
\begin{enumerate}\itemsep=0pt
\item[$(a)$] Let $\fg=\fo^{(1)}(2n)$. For a~basis of $H^2(\fg)$
we can take $($the classes of$)$ the following cocycles:
\begin{alignat*}{3}
 & \deg=-2\colon \quad && \text{none}&\nonumber\\
&&& c_{-1}^j= \sum \limits_{i\not = j}
\big(\widehat{E_{i,j}+E_{j+n,i+n}}\big)\wedge \big(\widehat{E_{j+n,i}+E_{i+n,j}}\big) \qquad \text{for~~}1\leq j \leq n,&\nonumber\\
&&& c_0= \sum \limits_{2\leq i\leq n}
\big(\widehat{E_{1,i}+E_{i+n,n+1}}\big)\wedge \big(\widehat{E_{i,1}+E_{n+1,i+n}}\big)&\nonumber\\
&&& \hphantom{c_0=}{}
+
 \sum \limits_{2 \leq i < j \leq n}
\big(\widehat{E_{i,n+j}+E_{j,n+i}}\big)\wedge \big(\widehat{E_{n+j,i}+E_{n+i,j}}\big).& 
\end{alignat*}

\item[$(b)$] For a~basis of $H^1\big(\fo^{(1)}
(2n);\fo^{(1)} (2n)\big)$ we can take the following derivations
\begin{alignat*}{3}
 & \deg=-1\colon \quad && \widetilde c_{-1}^j\colon \text{multiplication by the
matrix $d_j$, see Section~{\rm \ref{section7.1}}},&\nonumber\\
&&& \text{for any $j=1,\dots, n$ and $n>2$}, &\nonumber\\
& \deg=0\colon \quad && \text{multiplication by the matrix
 $D_n$, see Section~{\rm \ref{section7.1}}},&\nonumber\\
 &&& \text{for any $n>2$ and also the derivation}, &\nonumber\\
&&& c_0^2=
 \sum \limits_{1\leq i\leq 2n}
E_{i,i}\otimes \big(\big(\widehat{E_{1,1}+E_{n+1,n+1}}\big)+
\cdots + \big(\widehat{E_{n,n}+E_{2n,2n}}\big) \big).&
\end{alignat*}
\end{enumerate}
\end{Lemma}

\begin{Lemma}\label{3.2.17}\quad
\begin{enumerate}\itemsep=0pt
\item[$(a)$] Let $\fg=\fp(\fo^{(1)}(2n))$. For a~basis of $H^2(\fg)$ we can take $($the classes of$)$ the following cocycles:
\begin{alignat*}{3}
& \deg=-2\colon \quad && \text{none},&\nonumber\\
&&& c_{-1}^j= \sum\limits_{i\not = j}
\big(\widehat{E_{i,j}+
E_{j+n,i+n}}\big)\wedge \big(\widehat{E_{j+n,i}+E_{i+n,j}}\big) \qquad \text{for~~}1\leq j \leq n,&\nonumber\\
&&& c_0= \sum\limits_{2\leq i\leq n}
\big(\widehat{E_{1,i}+
E_{i+n,n+1}}\big)\wedge \big(\widehat{E_{i,1}+E_{n+1,i+n}}\big) &\nonumber\\
&&& \hphantom{c_0=}{}
+ \sum\limits_{2 \leq i < j \leq n} \big(\widehat{E_{i,n+j}+
E_{j,n+i}}\big)\wedge \big(\widehat{E_{n+j,i}+E_{n+i,j}}\big),&\\ 
& \text{for $n=4$}, \quad && \text{we have one more cocycle}&\nonumber\\
&&& \widetilde c_0= \big(\widehat{E_{1,3}+E_{7,5}}\big)\wedge \big(\widehat{E_{3,1}+E_{5,7}}\big)+\big(\widehat{E_{1,4}+
E_{8,5}}\big)\wedge \big(\widehat{E_{4,1}+E_{5,8}}\big)&\nonumber\\
&&& \hphantom{\widetilde c_0=}{}
+ \big(\widehat{E_{2,3}+E_{7,6}}\big)\wedge
\big(\widehat{E_{3,2}+E_{6,7}}\big)
 +\big(\widehat{E_{2,4}+E_{8,6}}\big)\wedge \big(\widehat{E_{4,2}+
E_{6,8}}\big)&\nonumber\\
&&& \hphantom{\widetilde c_0=}{}
+\big(\widehat{E_{3,8}+E_{4,7}}\big)\wedge \big(\widehat{E_{7,4}+E_{8,3}}\big)+\big(\widehat{E_{1,6}+E_{2,5}}\big)\wedge \big(\widehat{E_{5,2}+
E_{6,1}}\big).\nonumber
\end{alignat*}

\item[$(b)$] For a~basis of $H^1(\fg;\fg)$ we can take the following derivations
\begin{alignat*}{3}
& \deg=-1\colon \quad && \widetilde c_{-1}^j\colon \text{multiplication by the
matrix $d_j$, see Section~{\rm \ref{section7.1}}},&\nonumber\\
&&& \text{for any $j=1,\dots, n$ and $n>2$},& \\ 
&\deg=0\colon\quad && \text{multiplication by the matrix $D_n$, see Section~{\rm \ref{section7.1}}, for any $n>2$.}& \nonumber
\end{alignat*}
\end{enumerate}
\end{Lemma}

\begin{Lemma}\label{3.2.19}\quad
\begin{enumerate}\itemsep=0pt
\item[$(a)$] Let $\fg=\mathfrak{so} \fp(2n)$. For a~basis of $H^2(\fg)$
we can take $($the classes of$)$ the following cocycles:
\begin{alignat*}{3}
& \deg=-2\colon\quad && \text{none for $n\neq 4$},& \nonumber\\
&&& \text{for $n=4$, we have the cocycle}&\nonumber\\
&&& c_{-2}= \big(\widehat{E_{7,4}+E_{8,3}}\big)\wedge \big(\widehat{E_{5,2}+E_{6,1}}\big)+ \big(\widehat{E_{6,3}+
E_{7,2}}\big)\wedge \big(\widehat{E_{5,4}+E_{8,1}}\big)&\nonumber\\
&&& \hphantom{c_{-2}=}{}
+\big(\widehat{E_{6,4}+E_{8,2}}\big)\wedge \big(\widehat{E_{5,3}+E_{7,1}}\big),&\nonumber\\
&&& c_{-1}^{j}= \sum \limits_{i\not = j}
\big(\widehat{E_{i,j}+
E_{j+n,i+n}}\big)\wedge \big(\widehat{E_{j+n,i}+E_{i+n,j}}\big) \qquad \text{for~~}1\leq j \leq n, &\nonumber\\
& \deg=0\colon\quad && \text{none},&\nonumber\\
&&& c_{1}^{j}= \sum \limits_{i\not = j} \big(\widehat{E_{j,i}+
E_{i+n,j+n}}\big)\wedge \big(\widehat{E_{i, n+j}+E_{j,n+i}}\big)&\nonumber\\
&&& \hphantom{c_{1}^{j}=}{}
+\begin{cases}
\big(\widehat{E_{1,1}+E_{2,2}+E_{1+n,1+n}+E_{2+n,2+n}}\big)\wedge \big(\widehat{E_{1,n}}\big)\\
\quad \text{for } i=1,\\
\big(\widehat{E_{n-1,n-1}+E_{n,n}+E_{2n-1,2n-1}+E_{2n,2n}}\big)\wedge \big(\widehat{E_{n,2n}}\big)\\
\quad \text{for } i=n,\\
\big(\widehat{E_{j-1,j-1}+E_{j+1,j+1}+E_{j-1+n,j-1+n}+E_{j+1+n,j+1+n}}\big)\\
\quad \wedge \big(\widehat{E_{j,j+n}}\big) \ \ \text{for } j\not=1,n,
\end{cases}&\nonumber\\
& \deg=2\colon\quad && \text{cocycles as in \eqref{pogen}}.& 
\end{alignat*}

\item[$(b)$] For a~basis of $H^1(\fg;\fg)$ we can take the following derivations
\begin{alignat*}{3}
& \deg=-1\colon\quad && \text{none},& \nonumber\\
& \deg=0\colon\quad && \text{multiplication by the matrix
$D_n$, see Section~{\rm \ref{section7.1}}, for any $n>2$},&\nonumber\\
&&& \text{if $n$ is even, we have the additional derivation:}&\nonumber\\
&&& \text{multiplication by \smash{$d_{1,1}^{(n)}$}, see Section~{\rm \ref{section7.1}}},& \nonumber\\
& \deg=1\colon\quad && c_1^j\colon \text{multiplication by the matrix $d_j^{\rm t}$, see Section~{\rm \ref{section7.1}}},&\nonumber\\
&&& \text{for any $j=1,\dots, n$ and $n>2$}.& 
\end{alignat*}
\end{enumerate}
\end{Lemma}

\begin{Lemma}\label{3.2.21}\quad
\begin{enumerate}\itemsep=0pt
\item[$(a)$] Let $\fg=\mathfrak{so}
\fp(4n)/\Kee 1_{4n}$. For a~basis of $H^2(\fg)$ we can take $($the classes of$)$ the following cocycles:
\begin{alignat*}{3}
& \deg=-2\colon\quad && \text{none for $n\neq 2$. For $n=2$, we have the cocycle}, & \nonumber\\
&&& c_{-2}=\big(\widehat{E_{7,4}+E_{8,3}}\big)\wedge \big(\widehat{E_{5,2}+E_{6,1}}\big)+
\big(\widehat{E_{6,3}+E_{7,2}}\big)\wedge \big(\widehat{E_{5,4}+E_{8,1}}\big)&\nonumber\\
&&& \hphantom{c_{-2}=}{}
+\big(\widehat{E_{6,4}+E_{8,2}}\big)\wedge \big(\widehat{E_{5,3}+E_{7,1}}\big),& \nonumber\\
&&&
c_{-1}^{j}= \sum \limits_{i\not = j} \big(\widehat{E_{i,j}+
E_{j+2n,i+2n}}\big)\wedge \big(\widehat{E_{j+2n,i}+E_{i+2n,j}}\big) \quad \text{for~~}1\leq j \leq 2n,& \nonumber\\
&&& c_0= \sum \limits_{1\leq i\leq 2n}
\big(\widehat{E_{1,i}+E_{2n+i,2n+1}}\big)\wedge \big(\widehat{E_{i,1}+E_{2n+1,2n+i}}\big),&\nonumber\\
& \deg=1\colon\quad && \text{none},& \nonumber\\
& \deg=2\colon\quad && \text{cocycles as in \eqref{pogen}}.& 
\end{alignat*}

\item[$(b)$] For a~basis of $H^1(\fg;\fg)$ we can take the following
derivations
\begin{alignat*}{3}
& \deg=-2\colon\quad && \text{none for $n\neq 2$}, & \nonumber\\
&&& \text{for $n=2$, we have the cocycle}&\nonumber\\
&&& c _{-2}= (E_{3,8}+E_{4,7})\otimes \big(\widehat{E_{5,2}+E_{6,1}}\big)+
(E_{2,7}+E_{3,6})\otimes \big(\widehat{E_{5,4}+E_{8,1}}\big)&\nonumber\\
&&& \hphantom{c _{-2}=}{}
+(E_{2,8}+E_{4,6})\otimes \big(\widehat{E_{5,3}+E_{7,1}}\big)+
 (E_{1,6}+ E_{2,5})\otimes \big(\widehat{E_{7,4}+E_{8,3}}\big)&\nonumber\\
&&& \hphantom{c _{-2}=}{}
 +(E_{1,7}+E_{3,5})\otimes \big(\widehat{E_{6,4}+E_{8,2}}\big)+(E_{1,8}+
E_{4,5})\otimes \big(\widehat{E_{6,3}+E_{7,2}}\big), &\nonumber\\
& \deg=\pm1\colon\quad && \text{none}, &\nonumber\\
& \deg=0 \colon \quad && \text{multiplication by the matrix
$D_{2n}$, see Section~{\rm \ref{section7.1}}, for any $n>1$}, &\nonumber\\
&&& \text{and another derivation: multiplication by the matrix $d_{1,1}^{(n)}$},&\nonumber\\
&&& \text{see Section~{\rm \ref{section7.1}}.}& 
\end{alignat*}
\end{enumerate}
\end{Lemma}

\begin{Lemma}\label{3.2.23}\quad
\begin{enumerate}\itemsep=0pt
\item[$(a)$] For a~basis of $H^2(\mathfrak{spe}(2n))$
we can take $($the classes of$)$ the following cocycles:
\begin{alignat*}{3}
& \deg=-2\colon\quad && \widetilde c_{-2}^{i,j}=\big(\widehat{E_{n+i,j}+E_{n+j,i}}\big)^{\wedge^2} \qquad \text{for~~} 1\leq i <j\leq n, &\nonumber \\
&&& \text{for $n=4$, we have one more cocycle}&\nonumber\\
&& & c_{-2}=\big(\widehat{E_{7,4}+E_{8,3}}\big)\wedge \big(\widehat{E_{5,2}+E_{6,1}}\big)+ \big(\widehat{E_{6,3}+
E_{7,2}}\big)\wedge \big(\widehat{E_{5,4}+E_{8,1}}\big)&\nonumber\\
&&& \hphantom{c_{-2}=}{}
+\big(\widehat{E_{6,4}+E_{8,2}}\big)\wedge \big(\widehat{E_{5,3}+E_{7,1}}\big), & \nonumber\\
& \deg=-1, 0\colon\quad && \text{none}, &\nonumber\\
& \deg=1\colon\quad && c_{1}^{j}= \sum \limits_{i\not
= j} \big(\widehat{E_{j,i}+
E_{i+n,j+n}}\big)\wedge \big(\widehat{E_{i, n+j}+E_{j,n+i}}\big)& \nonumber\\
&&& \hphantom{c_{1}^{j}=}{} +\begin{cases}
\big(\widehat{E_{1,1}+E_{2,2}+E_{1+n,1+n}+E_{2+n,2+n}}\big)\wedge \big(\widehat{E_{1,n}}\big)\\
\quad \text{for } i=1,\\
\big(\widehat{E_{n-1,n-1}+E_{n,n}+E_{2n-1,2n-1}+E_{2n,2n}}\big)\wedge \big(\widehat{E_{n,2n}}\big)\\
\quad \text{for } i=n,\\
\big(\widehat{E_{j-1,j-1}+E_{j+1,j+1}+E_{j-1+n,j-1+n}+E_{j+1+n,j+1+n}}\big)\\
\quad \wedge \big(\widehat{E_{j,j+n}}\big) \ \ \text{for } j\not=1,n,
\end{cases}& \nonumber
\\
& \deg=2\colon\quad && \text{cocycles as in \eqref{pogen} and \eqref{pegen}}.& 
\end{alignat*}

\item[$(b)$] Let $\fg=\mathfrak{spe} (2n)$. For a
basis of $H^1(\fg;\fg)$ we can take the following derivations
\begin{alignat*}{3}
& \deg=-1\colon\quad && \text{none}& \nonumber\\
& \deg=0\colon\quad && \text{multiplication by the matrix
$D_n$, see Section~{\rm \ref{section7.1}}, for any $n>2$}, & \nonumber\\
&&& \text{if $n$ is even, we have one more derivation: multiplication}&\nonumber\\
&&& \text{by the matrix $d_{1,1}$, see Section~{\rm \ref{section7.1}}},& \nonumber\\
& \deg=1\colon\quad && c_1^j\colon \text{multiplication by the matrix $d_j^{\rm t}$, see Section~{\rm \ref{section7.1}}}, &\nonumber\\
&&& \text{for any $j=1,\dots, n$ and $n>2$}, & 
\end{alignat*}
\end{enumerate}
\end{Lemma}

\begin{Lemma}\label{3.2.25}\quad
\begin{enumerate}\itemsep=0pt
\item[$(a)$] For a~basis of $H^2(\mathfrak{spe}(4n)/
\Kee 1_{4n})$ we can take (the classes of) the following cocycles:
\begin{alignat*}{3}
& \deg=-2\colon\quad && \widetilde c_{-2}^{i,j}=\big(\widehat{E_{2n+i,j}+E_{2n+j,i}} \big)^{\wedge^2}
\qquad \text{for~~}1\leq i <j\leq 2n, &\nonumber\\
&&& \text{for $n=2$, we have one more cocycle} & \nonumber\\
&&& c_{-2}=\big(\widehat{E_{7,4}+E_{8,3}}\big)\wedge \big(\widehat{E_{5,2}+E_{6,1}}\big)+ \big(\widehat{E_{6,3}+
E_{7,2}}\big)\wedge (\widehat{E_{5,4}+E_{8,1}}\big)&\nonumber\\
&&& \hphantom{c_{-2}=}{}
+\big(\widehat{E_{6,4}+E_{8,2}}\big)\wedge \big(\widehat{E_{5,3}+E_{7,1}}\big),& \nonumber\\
& \deg=\pm 1\colon\quad && \text{none}, & \nonumber\\
&&& c_0= \sum \limits_{1\leq i\leq 2n}
\big(\widehat{E_{1,i}+E_{2n+i,2n+1}}\big)\wedge \big(\widehat{E_{i,1}+
E_{2n+1,2n+i}}\big), &\nonumber\\
& \deg=2\colon\quad && \text{cocycles as in \eqref{pogen} and
\eqref{pegen}}.& 
\end{alignat*}

\item[$(b)$] Let $\fg=\mathfrak{spe} (4n)/\Kee
1_{4n}$. For a~basis of $H^1(\fg;\fg)$ we can take the following
derivations
\begin{alignat*}{3}
& \deg=-2\colon\quad && \text{none for $n\neq 2$}, &\nonumber\\
&&& \text{for $n=2$, we have the cocycle}& \nonumber\\
&&& c _{-2}= (E_{3,8}+E_{4,7})\otimes \big(\widehat{E_{5,2}+E_{6,1}}\big)+(E_{2,7}+E_{3,6})\otimes \big(\widehat{E_{5,4}+
E_{8,1}}\big)&\nonumber\\
&&& \hphantom{c _{-2}=}{} +(E_{2,8}+E_{4,6})\otimes \big(\widehat{E_{5,3}+E_{7,1}}\big)+
 (E_{1,6}+ E_{2,5})\otimes \big(\widehat{E_{7,4}+E_{8,3}}\big)&\nonumber\\
&&& \hphantom{c _{-2}=}{}
 +(E_{1,7}+E_{3,5})\otimes \big(\widehat{E_{6,4}+E_{8,2}}\big)+(E_{1,8}+
E_{4,5})\otimes \big(\widehat{E_{6,3}+E_{7,2}}\big),\nonumber\\
& \deg=\pm 1\colon\quad && \text{none}, &\nonumber\\
& \deg=0\colon\quad && \text{multiplication by the matrix
$D_{2n}$, see Section~{\rm \ref{section7.1}}, for any $n>1$} &\nonumber\\
&&& \text{and one more derivation: multiplication by the matrix $d_{1,1}^{(2n)}$}, &\nonumber\\
&&& \text{see Section~{\rm \ref{section7.1}}.}& 
\end{alignat*}
\end{enumerate}
\end{Lemma}

\appendix\renewcommand{\thefootnote}{}

\section[Derivations and central extensions of the true deforms of symmetric simple\\ Lie algebras and superalgebras.
Appendix by Andrey Krutov]{Derivations and central extensions of the true deforms\\ of symmetric simple Lie algebras and superalgebras.\\
Appendix by Andrey Krutov\footnote{Institute of Mathematics, Czech Academy of Sciences, \v{Z}itn\'a 25, 115 67 Prague, Czech Republic}
\footnote{E-mail: \href{mailto:krutov@math.cas.cz}{krutov@math.cas.cz}}}

In this appendix, I consider derivations and central extensions of the following Lie algebras and superalgebras: true deforms with even parameter~$\eps$ (classified in~\cite{BGL2}) of symmetric Lie algebras and superalgebras with indecomposable Cartan matrix (classified in~\cite{BGL1}), and their simple relatives.

I do not consider deforms with odd parameter, although this might be the most interesting task from the ``super'' point of view.

In this appendix, $\fg_{c_i}$ designates the deform with parameter $\eps$ of Lie (super)algebras $\fg$ determined by the cocycle $c_i$ explicitly given in
\cite{BGL2}. The results below are obtained with the help of \textsc{SuperLie} package~\cite{Gr}.

\begin{Lemma} \label{AKwk4} Let $p=2$ and $\fg=\fwk(4;\alpha)$ for $\alpha\neq 0,1$. Then
$H^1(\fg_{c_i};\fg_{c_i}) = 0$ and $H^2(\fg_{c_i}) = 0$ for all $i\neq 0$.
For $\fg_{c_0}$, we have
\begin{enumerate}\itemsep=0pt
\item[$1)$] if $\eps \neq \alpha$, then $H^1(\fg_{c_i};\fg_{c_i}) = 0$ and $H^2(\fg_{c_i}) = 0$;
\item[$2)$] if $\eps=\alpha$, then $H^2(\fg)=0$ and $H^1(\fg_{c_0}; \fg_{c_0})$ is spanned by the cocycle
\begin{gather*}
x_1\otimes \widehat{x}_1 + x_2\otimes \widehat{x}_2 + x_4\otimes \widehat{x}_4 +
x_6\otimes \widehat{x}_6 + x_7\otimes \widehat{x}_7 + x_{10}\otimes \widehat{x}_{10}
+ x_{11}\otimes \widehat{x}_{11} \\
\qquad{} + x_{14}\otimes \widehat{x}_{14} +y_1\otimes
\widehat{y}_1 + y_2\otimes \widehat{y}_2 + y_4\otimes \widehat{y}_4 + y_6\otimes
\widehat{y}_6 + y_7\otimes \widehat{y}_7 + y_{10}\otimes \widehat{y}_{10}\\
\qquad{} + y_{11}\otimes \widehat{y}_{11} + y_{14}\otimes \widehat{y}_{14}.
\end{gather*}
\end{enumerate}
\end{Lemma}

\begin{Lemma} \label{AKwk3big} Let $p=2$ and $\fg=\fwk(3;\alpha)$ for~$\alpha\neq 0,1$.
Then, the space $H^1(\fg_{c_i};\fg_{c_i})$ for~$i\neq0$ is spanned by
\[
 (h_1 + \alpha h_3)\otimes\widehat{h}_4.
\]
The space $H^1(\fg_{c_0};\fg_{c_0})$ is spanned by
\[
 (h_1 + \alpha h_3)\otimes\widehat{h}_4 + \eps(h_1 + \alpha h_3)\otimes\widehat{h}_2.
\]
We have $H^2(\fg_{c_i})=0$ for all cocycles~$c_i$.
\end{Lemma}

\begin{Lemma} \label{AKwk3der} Let $p=2$ and $\fg=\fwk^{(1)}(3;\alpha)$ for $\alpha\neq 0,1$. Then, the
space $H^1(\fg_{c_i};\fg_{c_i})$ is
spanned by one and the same derivation for any cocycle $c_i$ for $i\neq0$:
\[
 x_1\otimes\widehat x_1 + x_3\otimes\widehat x_3 + x_4\otimes\widehat
x_4 + x_5\otimes\widehat x_5 + y_1\otimes\widehat y_1+y_3\otimes\widehat y_3
 + y_4\otimes\widehat y_4 + y_5\otimes\widehat y_5.
\]
We have $H^2(\fg_{c_i}) = 0$ for all cocycles~$c_i$.

For cocycle~$c_0$, we have
\begin{enumerate}\itemsep=0pt
\item[$1)$] if $\eps\notin\{1,\alpha\}$, then the space $H^1(\fg_{c_0};\fg_{c_0})$ is
spanned by
\[
x_2\otimes\widehat{x}_2+x_3\otimes\widehat{x}_3+x_4\otimes\widehat{x}_4+x_7\otimes\widehat{x}_7+y_2\otimes\widehat{y}_2+y_3\otimes\widehat{y}_3+y_4\otimes\widehat{y}_4+y_7\otimes\widehat{y}_7;
\]
\item[$2)$] if $\eps = 1$, then the space $H^1(\fg_{c_0};\fg_{c_0})$ is spanned by
\begin{gather*}
 d_{-2} = y_3\otimes\widehat{x}_3 \qquad \text{(but there is no $d_2$!)},\\
 d_{0,1}= \alpha h_3\otimes\widehat{h}_2+h_1\otimes\widehat{h}_2,\\
 d_{0,2}= \alpha
h_3\otimes\widehat{h}_3+h_1\otimes\widehat{h}_3+\alpha
x_2\otimes\widehat{x}_2+\alpha x_4\otimes\widehat{x}_4+\alpha
x_5\otimes\widehat{x}_5\\
\hphantom{d_{0,2}=}{} +\alpha x_6\otimes\widehat{x}_6+\alpha
y_3\otimes\widehat{y}_3,\\
 d_{0,3}=
x_1\otimes\widehat{x}_1+x_3\otimes\widehat{x}_3+x_4\otimes\widehat{x}_4+x_5\otimes\widehat{x}_5+y_1\otimes\widehat{y}_1+y_3\otimes\widehat{y}_3\\
\hphantom{d_{0,3}=}{}
+y_4\otimes\widehat{y}_4+y_5\otimes\widehat{y}_5.
\end{gather*}
\item[$3)$] if $\eps = \alpha$, then the space $H^1(\fg_{c_0};\fg_{c_0})$ is spanned by
\[
x_1\otimes\widehat{x}_1+x_3\otimes\widehat{x}_3+x_4\otimes\widehat{x}_4+x_5\otimes\widehat{x}_5+y_1\otimes\widehat{y}_1
+y_3\otimes\widehat{y}_3+y_4\otimes\widehat{y}_4+y_5\otimes\widehat{y}_5.
\]
\end{enumerate}
\end{Lemma}

\begin{Remark} The weights of basis elements in the deform of~$\fwk(3;\alpha)$ that correspond to the
cocycle~$c_0$ remain unchanged, whereas it is not so in the case of the deform of~$\fwk^{(1)}(3;\alpha)$. In this case,
the weights depend on the deformation parameter~$\eps$, which causes a~\lq\lq spontaneous breach of symmetry'', as a~physicist might say, see Lemma~\ref{AKwk3der} for $\eps=1$: no cocycles of weight~2.
\end{Remark}

\begin{Lemma} \label{AKwk3simple} Let $p=2$ and $\fg=\fwk^{(1)}(3;\alpha)/\fc$ for $\alpha\neq0,1$, where $\fc$ is
a~$1$-dimensional center spanned by $h_1+\alpha h_3$.
Then, the space $H^1(\fg_{c_i};\fg_{c_i})$ for all $i\neq 0$ is spanned by
\[
 x_1\otimes\widehat{x}_1+x_3\otimes\widehat{x}_3+x_4\otimes\widehat{x}_4+x_5\otimes\widehat{x}_5+y_1\otimes\widehat{y}_1+y_3\otimes\widehat{y}_3+y_4\otimes\widehat{y}_4+y_5\otimes\widehat{y}_5,
\]
the space $H^2(\fg_{c_i})$ for all $i\neq0$ is spanned by
\[
\big(\alpha^2+\alpha\big) \widehat{x}_7\wedge \widehat{y}_7+\alpha \widehat{x}_4\wedge \widehat{y}_4+\alpha \widehat{x}_6\wedge \widehat{y}_6+\widehat{x}_1\wedge \widehat{y}_1.
\]
For the cocycle~$c_0$, we have:
\begin{enumerate}\itemsep=0pt
\item[$1)$] if $\eps\notin\{1,\alpha\}$, then the space $H^1(\fg_{c_0};\fg_{c_0})$ is spanned by
\[
x_2\otimes\widehat{x}_2+x_3\otimes\widehat{x}_3+x_4\otimes\widehat{x}_4+x_7\otimes\widehat{x}_7+y_2\otimes\widehat{y}_2+y_3\otimes\widehat{y}_3
+y_4\otimes\widehat{y}_4+y_7\otimes\widehat{y}_7,
\]
and the space $H^2(\fg_{c_0})$ is spanned by
\begin{gather*}
 \big(\alpha^2\eps^2 + \alpha^2 + \alpha\eps^3 + \alpha\eps^2 + \alpha\eps + \alpha + \eps^3 + \eps\big) \widehat{x}_7\wedge \widehat{y}_7 \\
\qquad{} + \big(\alpha\eps + \alpha + \eps^2 + \eps\big) \widehat{x}_4\wedge \widehat{y}_4 + \big(\alpha\eps + \alpha + \eps^2 + \eps\big) \widehat{x}_6\wedge \widehat{y}_6 + \widehat{x}_1\wedge \widehat{y}_1,
\end{gather*}
\item[$2)$] if $\eps=1$, then the space $H^1(\fg_{c_0};\fg_{c_0})$ is spanned by
\begin{gather*}
 d_{0,1} = x_1\otimes\widehat{x}_1+x_3\otimes\widehat{x}_3+x_4\otimes\widehat{x}_4+x_5\otimes\widehat{x}_5+y_1\otimes\widehat{y}_1+y_3\otimes\widehat{y}_3\\
\hphantom{d_{0,1} =}{}
 +y_4\otimes\widehat{y}_4+y_5\otimes\widehat{y}_5,\\
 d_{0,2} = x_2\otimes\widehat{x}_2+x_4\otimes\widehat{x}_4+x_5\otimes\widehat{x}_5+x_6\otimes\widehat{x}_6,
\end{gather*}
and the space $H^2(\fg_{c_0})$ is spanned by
\begin{gather*}
 z_{0,1} = \widehat{x}_1\wedge\widehat{y}_1,\\
 z_{0,2} = \alpha \widehat{x}_4\wedge \widehat{y}_4+\alpha \widehat{x}_6\wedge \widehat{y}_6+\widehat{x}_2\wedge \widehat{y}_2+\widehat{x}_5\wedge \widehat{y}_5,
\end{gather*}
\item[$3)$] if $\eps=\alpha$, then the space $H^1(\fg_{c_0};\fg_{c_0})$ is spanned by
\[
x_1\otimes\widehat{x}_1+x_3\otimes\widehat{x}_3+x_4\otimes\widehat{x}_4+x_5\otimes\widehat{x}_5+y_1\otimes\widehat{y}_1+y_3\otimes\widehat{y}_3+y_4\otimes\widehat{y}_4+y_5\otimes\widehat{y}_5,
\]
and $H^2(\fg_{c_0})=\operatorname{Span}(\widehat{x}_1\wedge\widehat{y}_1)$.
\end{enumerate}
\end{Lemma}

\begin{Lemma} \label{AKbr3} Let $p=3$ and $\fg=\mathfrak{br}(3)$. Then, $H^1(\fg_{c_i};\fg_{c_i}) = 0$ and $H^2(\fg_{c_i}) = 0$ for any $c_i$.
\end{Lemma}

\begin{Lemma} Let $p=3$ and $\fg=\mathfrak{L}(2,2)$, the explicitly described in {\rm \cite{BLW}} exceptional deform~of~$\fo(5)$.
Then, $H^1(\fg;\fg) = 0$ and $H^2(\fg) = 0$.
\end{Lemma}

\begin{Lemma} Let $p=3$ and $\fg=\mathfrak{brj}(2;3)$. Then, the
space~$H^1(\fg_{c};\fg_{c})$ is spanned by the following cocycles=derivations, where $\eps$ is the parameter of deformation that defines $\fg_{c}$
 for $c=c_{12, i}$ and $c=c_{6, i}$, where $i=1,2\colon$
\begin{gather*}
c_{12,1} = 2 \eps x_5\otimes \widehat{y}_8+\eps x_8\otimes
\widehat{y}_5 +2 y_2\otimes \widehat{x}_4
 +y_4\otimes \widehat{x}_2+x_1\otimes \widehat{x}_6+x_3\otimes \widehat{x}_7\\
 \hphantom{c_{12,1} =}{}
 + 2 y_6\otimes \widehat{y}_1+y_7\otimes \widehat{y}_3,\\
c_{12,2} = x_2\otimes \widehat{y}_4+x_4\otimes \widehat{y}_2+x_6\otimes\widehat{x}_1+2 x_7\otimes \widehat{x}_3
 +y_1\otimes \widehat{y}_6+y_3\otimes \widehat{y}_7,
\\
c_{6,1} =
 2 y_2\otimes \widehat{x}_4+y_4\otimes \widehat{x}_2+x_1\otimes
\widehat{x}_6+x_3\otimes \widehat{x}_7+2y_6\otimes
\widehat{y}_1+y_7\otimes \widehat{y}_3,\\
c_{6,2} =
 2 \eps x_5\otimes \widehat{y}_8+\eps x_8\otimes
\widehat{y}_5+x_2\otimes \widehat{y}_4+x_4\otimes
\widehat{y}_2+x_6\otimes \widehat{x}_1+2 x_7\otimes
\widehat{x}_3+y_1\otimes \widehat{y}_6+y_3\otimes \widehat{y}_7.
\end{gather*}

We have~$H^2(\fg_{c_i})=0$ for any cocycle $c_i$.
\end{Lemma}

A direct corollary of the above lemma is that the deforms of
$\fwk^{(1)}(3;\alpha)/\fc$ and $\fwk(4;\alpha)$
corresponding to the cocyles~$c_0$ when the deformation parameter $\eps$ is
equal to~$\alpha$ are not isomorphic
to the original algebras. This fact agrees with results of~\cite{BKLS}, where
it was shown that these deforms are
not simple Lie algebras.

\subsection*{Acknowledgements} We are thankful to A.~Krutov for his help, e.g., for computing several examples (Section~\ref{sssAK}). We
are thankful to N.~Chebochko and M.~Kuznetsov for helpful
discussions of their unpublished results pertaining to this paper; to A.~Dzhumadildaev, P.~Zusmanovich,
and Sh.~Ibraev for helpful discussions. We are very thankful to the referees,
carefully selected by SIGMA, especially one of them, for very constructive criticism.

S.B.~and D.L.~were supported by the grant AD 065 NYUAD. D.L.~is thankful to MPIMiS, Leipzig, where he was Sophus-Lie-Professor
(2004-07), when a~part of the ideas of this paper were conceived, for financial
support and most creative environment. The authors of the main text and A.~Krutov, who wrote the Appendix, are grateful to M.~Al Barwani, Director of the High Performance Computing resources at New York University Abu Dhabi for the possibility to perform the difficult computations of this research.

Andrey Krutov was supported by the GA\v{C}R project 20-17488Y and RVO: 67985840.

\LastPageEnding

\end{document}